 \documentclass[a4paper,12pt]{book}

 \usepackage{amsfonts,amsmath,mathrsfs,amssymb}
 \usepackage{times,helvet,courier,type1cm}

 \usepackage{makeidx}
 \usepackage{harvard}

 \usepackage[dvips]{color}

 \allowdisplaybreaks

 \newtheorem{thm0}{Theorem}[section]
 \newtheorem{exa0}{Theorem}[section]

 \newtheorem{def1}[thm0]{Definition}
 \newtheorem{lem1}[thm0]{Lemma}
 \newtheorem{thm1}[thm0]{Theorem}
 \newtheorem{cor1}[thm0]{Corollary}
 \newtheorem{pro1}[thm0]{Proposition}
 \newtheorem{con1}[thm0]{Condition}
 \newtheorem{exa1}[exa0]{\it{Example}}

 \def\bglemma{\begin{lem1}}\def\edlemma{\end{lem1}}
 \def\bgtheorem{\begin{thm1}}\def\edtheorem{\end{thm1}}
 \def\bgcorollary{\begin{cor1}}\def\edcorollary{\end{cor1}}
 \def\bgproposition{\begin{pro1}}\def\edproposition{\end{pro1}}
 \def\bgcondition{\begin{con1}}\def\edcondition{\end{con1}}
 \def\bexample{\begin{exa1}\rm{}\def\eexample{\end{exa1}}}

 \def\benumerate{\begin{enumerate}}
 \def\eenumerate{\end{enumerate}}
 \def\itm{\item}

 \def\beqlb{\begin{eqnarray}}\def\eeqlb{\end{eqnarray}}
 \def\beqnn{\begin{eqnarray*}}\def\eeqnn{\end{eqnarray*}}

 \def\eqref#1{{\rm(\ref{#1})}}

 \def\ar{\!\!\!&}\def\nnm{\nonumber}\def\ccr{\nnm\\}

 \def\mcr{\mathscr}\def\mbb{\mathbb}\def\mbf{\mathbf}
 \def\mrm{\mathrm}

 \def\qed{\hfill$\square$\smallskip}

 \def\<{\langle}\def\>{\rangle}

 \def\ulim{\uparrow\!\!\lim}\def\dlim{\downarrow\!\!\lim}

 \def\d{\mrm{d}}\def\e{\mrm{e}}
 \def\b{\mrm{b}}

 \def\x{x}

 \def\const{\mrm{const}}

 \makeindex

\begin{document}

\


\vskip6.0cm

\centerline{\LARGE\sf{\textbf{CONTINUOUS-STATE BRANCHING}}}

\bigskip

\centerline{\LARGE\sf{\textbf{PROCESSES}}}

\vskip2.0cm

\centerline{\Large\sf{Zenghu Li}}

\vskip3.0cm

\centerline{\large\rm{(February 15, 2012)}}

\vskip5.0cm

\centerline{\large\rm{Beijing Normal University}}

\newpage


\

\vskip3.0cm

Professor Zenghu Li

School of Mathematical Sciences

Beijing Normal University

Beijing 100875, China

E-mail: lizh@bnu.edu.cn

URL: http://math.bnu.edu.cn/\~{}lizh/

\vskip9.0cm

Title: Continuous-state branching processes

\bigskip



\chapter*{Preface}

\pagenumbering{roman}

\bigskip

These notes were used in a short graduate course on branching processes
the author gave in Beijing Normal University. The following main topics
are covered: scaling limits of Galton--Watson processes, continuous-state
branching processes, extinction probabilities, conditional limit theorems,
decompositions of sample paths, martingale problems, stochastic equations,
Lamperti's transformations, independent and dependent immigration
processes. Some of the results are simplified versions of those in the
author's book ``Measure-valued branching Markov processes'' (Springer,
2011). We hope these simplified results will set out the main ideas in an
easy way and lead the reader to a quick access of the subject.

\vskip1.0cm

\noindent Zenghu Li

\medskip

\noindent Beijing, China

\tableofcontents


\chapter{Preliminaries}\label{ch01}

\pagenumbering{arabic}

In this chapter, we discuss the basic properties of Laplace transforms of
finite measures on the positive half line. In particular, we give some
characterizations of the weak convergence of those measures in terms of
their Laplace transforms. Based on these results, a general
representation for infinitely divisible distributions on the positive
half line is established. We also give some characterizations of
continuous functions on the positive half line with L\'{e}vy--Khintchine
type representations.

\section{Laplace transforms of measures}\label{sch01.1}

 \setcounter{equation}{0}

In this section, we discuss the basic properties of Laplace transforms of
finite measures on the positive half line $\mbb{R}_+ := [0,\infty)$. Let
$B(\mbb{R}_+) = \b\mcr{B}(\mbb{R}_+)$ be the set of bounded Borel
functions on $\mbb{R}_+$. Given a finite measure $\mu$ on $\mbb{R}_+$, we
define the \index{Laplace functional} \textit{Laplace transform} $L_\mu$
of $\mu$ by
 \beqlb\label{ch01.1.1}
L_\mu(\lambda)
 =
\int_0^\infty \e^{-\lambda\x} \mu(\d\x), \qquad \lambda\ge 0.
 \eeqlb

 \bgtheorem\label{tch01.1.1}
A finite measure on $\mbb{R}_+$ is uniquely determined by its Laplace
transform. \edtheorem

 \begin{proof}
Suppose that $\mu_1$ and $\mu_2$ are finite measures on $\mbb{R}_+$ and
$L_{\mu_1}(\lambda) = L_{\mu_2}(\lambda)$ for all $\lambda\ge 0$. Let
$\mcr{K} = \{x\mapsto \e^{-\lambda\x}: \lambda\ge 0\}$ and let $\mcr{L}$
be the class of functions $F\in B(\mbb{R}_+)$ so that
 \beqnn
\int_0^\infty F(x)\mu_1(\d x) = \int_0^\infty F(x)\mu_2(\d x).
 \eeqnn
Then $\mcr{K}$ is closed under multiplication and $\mcr{L}$ is a monotone
vector space containing $\mcr{K}$. It is easy to see $\sigma(\mcr{K}) =
\mcr{B}(\mbb{R}_+)$. Then the monotone class theorem implies
$\mcr{L}\supset \b\sigma(\mcr{K}) = B(\mbb{R}_+)$. That proves the
desired result. \qed\end{proof}

 \bgtheorem\label{tch01.1.2}
Let $\{\mu_n\}$ be finite measures on $\mbb{R}_+$ and let $\lambda\mapsto
L(\lambda)$ be a continuous function on $[0,\infty)$. If there is a dense
subset $D$ of $(0,\infty)$ to that $\lim_{n\to \infty} L_{\mu_n}(\lambda)
= L(\lambda)$ for every $\lambda\in D$, then there is a finite measure
$\mu$ on $\mbb{R}_+$ such that $L_\mu=L$ and $\lim_{n\to \infty} \mu_n =
\mu$ by weak convergence. \edtheorem

 \begin{proof}
We can regard each $\mu_n$ as a finite measure on $\bar{\mbb{R}}_+ :=
[0,\infty]$, the one-point compactification of $[0,\infty)$. Let $F_n$
denote the distribution function of $\mu_n$. By applying Helly's theorem
one can see that any subsequence of $\{F_n\}$ contains a weakly
convergent subsequence $\{F_{n_k}\}$. Then the corresponding subsequence
$\{\mu_{n_k}\}$ converges weakly on $\bar{\mbb{R}}_+$ to a finite measure
$\mu$. It follows that
 \beqnn
\mu(\bar{\mbb{R}}_+) = \lim_{k\to \infty} \mu_{n_k}(\mbb{R}_+)
 =
\lim_{k\to \infty} L_{\mu_{n_k}}(0) = L(0).
 \eeqnn
Moreover, for $\lambda\in D$ we have
 \beqlb\label{ch01.1.2}
\int_{\bar{\mbb{R}}_+} \e^{-\lambda\x} \mu(\d\x)
 =
\lim_{k\to \infty} \int_0^\infty \e^{-\lambda\x} \mu_{n_k}(\d\x)
 =
L(\lambda),
 \eeqlb
where $\e^{-\lambda\cdot\infty} = 0$ by convention. By letting
$\lambda\to 0+$ along $D$ in \eqref{ch01.1.2} and using the continuity of
$L$ at $\lambda=0$ we find $\mu(\mbb{R}_+) = L(0)$, so $\mu$ is supported
by $\mbb{R}_+$. Then $\lim_{n\to \infty} \mu_{n_k} = \mu$ weakly on
$\mbb{R}_+$. It is easy to see that \eqref{ch01.1.2} in fact holds for
all $\lambda\ge 0$, so we have $L_\mu=L$. By a standard argument one sees
$\lim_{n\to \infty} \mu_n = \mu$ weakly on $\mbb{R}_+$. \qed\end{proof}

 \bgtheorem\label{tch01.1.3}
Let $\mu_1, \mu_2, \ldots$ and $\mu$ be finite measures on $\mbb{R}_+$.
Then $\mu_n \to \mu$ weakly if and only if $L_{\mu_n}(\lambda)\to
L_\mu(\lambda)$ for every $\lambda\ge 0$. \edtheorem

 \begin{proof}
If $\mu_n \to \mu$ weakly, we have $\lim_{n\to \infty} L_{\mu_n}(\lambda)
= L_\mu(\lambda)$ for every $\lambda\ge 0$ by dominated convergence. The
converse assertion is a consequence of Theorem~\ref{tch01.1.2}.
\qed\end{proof}

We next give a necessary and sufficient condition for a continuous real
function to be the Laplace transform of a finite measure on $\mbb{R}_+$.
For a constant $c\ge 0$ and a function $f$ on an interval $T\subset
\mbb{R}$ we write
 \beqnn
\Delta_c f(\lambda) = f(\lambda+c) - f(\lambda), \qquad
\lambda,\lambda+c\in T.
 \eeqnn
Let $\Delta_c^0$ be the identity and define $\Delta_c^n = \Delta_c^{n-1}
\Delta_c$ for $n\ge 1$ inductively. Then we have
 \beqnn
\Delta_c^m f(\lambda) = (-1)^m \sum_{i=0}^m {m\choose i} (-1)^i
f(\lambda+ic).
 \eeqnn
The \index{Bernstein polynomials} \textit{Bernstein polynomials} of a
function $f$ on $[0,1]$ are given by
 \beqlb\label{ch01.1.3}
B_{f,m}(s) = \sum_{i=0}^m {m \choose i} \Delta_{1/m}^if(0)s^i, \qquad
0\le s\le 1, m=1,2,\ldots.
 \eeqlb
It is well-known that $B_{f,m}(s) \to f(s)$ uniformly as $m\to \infty$;
see, e.g., Feller (1971, p.222). A real function $\theta$ on $[0,\infty)$
is said to be \index{completely monotone function} \textit{completely
monotone} if it satisfies
 \beqlb\label{ch01.1.4}
(-1)^i\Delta_c^i\theta(\lambda)\ge 0, \qquad \lambda\ge 0,c\ge
0,i=0,1,2,\ldots.
 \eeqlb

 \bgtheorem\label{tch01.1.4}
A continuous real function $\theta$ on $[0,\infty)$ is the Laplace
transform of a finite measure $\mu$ on $\mbb{R}_+$ if and only if it is
completely monotone. \edtheorem

 \begin{proof}
If $\theta$ is the Laplace transform of a finite measure on $\mbb{R}_+$
it is clearly a completely monotone function. Conversely, suppose that
\eqref{ch01.1.4} holds. For fixed $a>0$, we let $\gamma_a(s) = \theta
(a-as)$ for $0\le s\le 1$. The complete monotonicity of $\theta$ implies
 \beqnn
\Delta_{1/m}^i\gamma_a(0)\ge 0, \qquad i=0,1,\ldots, m.
 \eeqnn
Then the Bernstein polynomial $B_{\gamma_a,m}(s)$ has positive
coefficients, so $B_{\gamma_a,m}(\e^{-\lambda/a})$ is the Laplace
transform of a finite measure on $\mbb{R}_+$. By Theorem~\ref{tch01.1.2},
 \beqnn
\theta(\lambda) = \lim_{a\to\infty}\lim_{m\to\infty} B_{\gamma_a,m}
(\e^{-\lambda/a}), \qquad \lambda\ge 0,
 \eeqnn
is the Laplace transform of a finite measure on $\mbb{R}_+$.
\qed\end{proof}

We often use a variation of the Laplace transform in dealing with
$\sigma$-finite measures on $(0,\infty)$. A typical case is considered in
the following:

 \bgtheorem\label{tch01.1.5}
Let $\mu_1$ and $\mu_2$ be two $\sigma$-finite measures on $(0,\infty)$.
If for every $\lambda\ge 0$,
 \beqlb\label{ch01.1.5}
\int_0^\infty (1-\e^{-\lambda\x}) \mu_1(\d\x)
 =
\int_0^\infty (1-\e^{-\lambda\x}) \mu_2(\d\x)
 \eeqlb
and the value is finite, then we have $\mu_1 = \mu_2$. \edtheorem

 \begin{proof}
By setting $\mu_1(\{0\}) = \mu_2(\{0\}) = 0$ we extend $\mu_1$ and
$\mu_2$ to $\sigma$-finite measures on $[0,\infty)$. Taking the
difference of \eqref{ch01.1.5} for $\lambda$ and $\lambda+1$ we obtain
 \beqnn
\int_0^\infty \e^{-\lambda\x}(1-\e^{-\x}) \mu_1(\d\x)
 =
\int_0^\infty \e^{-\lambda\x}(1-\e^{-\x}) \mu_2(\d\x).
 \eeqnn
Then the result of Theorem~\ref{tch01.1.1} implies that
 \beqnn
(1-\e^{-\x}) \mu_1(\d\x)
 =
(1-\e^{-\x}) \mu_2(\d\x)
 \eeqnn
as finite measures on $[0,\infty)$. Since $1-\e^{-\x}$ is strictly
positive on $(0,\infty)$, it follows that $\mu_1 = \mu_2$ as
$\sigma$-finite measures on $(0,\infty)$. \qed\end{proof}

Now let us consider a complete separable metric space $E$ with the Borel
$\sigma$-algebra denoted by $\mcr{B}(E)$. Suppose that $h$ is a strictly
positive bounded Borel function on $E$. Let $B_h(E)$ be the set of Borel
functions $f$ on $E$ such that $|f|\le \const\cdot h$. Let $M_h$ be the
set of Borel measures $\mu$ on $E$ such that $\int_E h\d\mu< \infty$. Let
$\mcr{M}_h$ be the $\sigma$-algebra on $M_h$ generated by the mappings
 \beqnn
\mu\mapsto \mu(f) := \int_E f(x)\mu(\d x), \qquad f\in B_h(E).
 \eeqnn
Given a finite measure $Q$ on $(M_h,\mcr{M}_h)$, we define the
\index{Laplace functional} \textit{Laplace functional} $L_Q$ of $Q$ by
 \beqlb\label{ch01.1.6}
L_Q(f)=\int_{M_h} \e^{-\nu(f)} Q(\d\nu ), \qquad f\in B_h(E)^+.
 \eeqlb
A random element $X$ taking values on $(M_h,\mcr{M}_h)$ is called a
\index{random measure} \textit{random measure} on $E$. The \index{Laplace
functional} \textit{Laplace functional} of a random measure means the
Laplace functional of its distribution on $(M_h,\mcr{M}_h)$. The reader
may refer to Kallenberg (1975) or Li (2011) for the basic theory of
random measure. In particular, the proofs of the following results can be
found in the two references:

 \bgtheorem\label{tch01.1.6}
A finite measure on $(M_h,\mcr{M}_h)$ is uniquely determined by its
Laplace functional. \edtheorem

Suppose that $\lambda$ is a $\sigma$-finite measure on $(E,\mcr{B}(E))$.
A random measure $X$ on $E$ is called a \index{Poisson random measure}
\textit{Poisson random measure} with \index{intensity of a Poisson random
measure} \textit{intensity} $\lambda$ provided:
 \benumerate

\itm[{\rm(1)}] for each $B\in\mcr{B}(E)$ with $\lambda(B)< \infty$, the
random variable $X(B)$ has the Poisson distribution with parameter
$\lambda(B)$, that is,
 \beqnn
\mbf{P}\{X(B)=n\} = \frac{\lambda(B)^n}{n!}\e^{-\lambda(B)}, \qquad
n=0,1,2,\ldots;
 \eeqnn

\itm[{\rm(2)}] if $B_1, \ldots, B_n\in\mcr{B}(E)$ are disjoint and
$\lambda(B_i)< \infty$ for each $i=1, \ldots, n$, then $X(B_1),$
$\ldots,$ $X(B_n)$ are mutually independent random variables.

 \eenumerate

 \bgtheorem\label{tch01.1.7}
A random measure $X$ on $E$ is Poissonian with intensity $\lambda \in
M_h(E)$ if and only if its Laplace functional is given by
 \beqlb\label{ch01.1.7}
\mbf{E}\exp\{-X(f)\}
 =
\exp\Big\{-\int_E(1-\e^{-f(x)})\lambda(\d x)\Big\}, \quad f\in B_h(E)^+.
 \eeqlb
\edtheorem

 \begin{proof}
Suppose that $X$ is a Poisson random measure on $E$ with intensity
$\lambda$. Let $B_1, \ldots, B_n\in\mcr{B}(E)$ be disjoint sets
satisfying $\lambda(B_i)< \infty$ for each $i=1, \ldots, n$. For any
constants $\alpha_1, \ldots, \alpha_n\ge 0$ we can use the above two
properties to see
 \beqlb\label{ch01.1.8}
\mbf{E}\exp\Big\{- \sum_{i=1}^n\alpha_iX(B_i)\Big\}
 =
\exp\Big\{- \sum_{i=1}^n (1-\e^{-\alpha_i})\lambda(B_i)\Big\}.
 \eeqlb
Then we get \eqref{ch01.1.7} by approximating $f\in B_h(E)^+$ by simple
functions and using dominated convergence. Conversely, if the Laplace
functional of $X$ is given by \eqref{ch01.1.7}, we may apply the equality
to the simple function $f = \sum_{i=1}^n \alpha_i1_{B_i}$ to get
\eqref{ch01.1.8}. Then $X$ satisfies the above two properties in the
definition of a Poisson random measure on $E$ with intensity $\lambda$.
\end{proof}

\section{Infinitely divisible distributions}\label{sch01.2}

 \setcounter{equation}{0}

For probability measures $\mu_1$ and $\mu_2$ on $\mbb{R}_+$, the product
$\mu_1 \times \mu_2$ is a probability measure on $\mbb{R}_+^2$. The image
of $\mu_1 \times \mu_2$ under the mapping $(\x_1,\x_2)\mapsto \x_1+\x_2$
is called the \index{convolution} \textit{convolution} of $\mu_1$ and
$\mu_2$ and is denoted by $\mu_1*\mu_2$, which is a probability measure
on $\mbb{R}_+$. According to the definition, for any $F\in B(\mbb{R}_+)$
we have
 \beqlb\label{ch01.2.1}
\int_0^\infty F(\x) (\mu_1*\mu_2)(\d\x)
 =
\int_0^\infty \mu_1(\d\x_1)\int_0^\infty F(\x_1+\x_2) \mu_2(\d\x_2).
 \eeqlb
Clearly, if $\xi_1$ and $\xi_2$ are independent random variables with
distributions $\mu_1$ and $\mu_2$ on $\mbb{R}_+$, respectively, then the
random variable $\xi_1+\xi_2$ has distribution $\mu_1*\mu_2$. It is easy
to show that
 \beqlb\label{ch01.2.2}
L_{\mu_1*\mu_2}(\lambda) = L_{\mu_1}(\lambda)L_{\mu_2}(\lambda), \qquad
\lambda\ge 0.
 \eeqlb
Let $\mu^{*0} = \delta_0$ and define $\mu^{*n} = \mu^{*(n-1)}*\mu$
inductively for integers $n\ge 1$. We say a probability distribution
$\mu$ on $\mbb{R}_+$ is \index{infinitely divisible distribution}
\textit{infinitely divisible} if for each integer $n\ge1$, there is a
probability $\mu_n$ such that $\mu = \mu_n^{*n}$. In this case, we call
$\mu_n$ the \index{the $n$-th root of a probability} \textit{$n$-th root}
of $\mu$. A positive random variable $\xi$ is said to be
\index{infinitely divisible random measure} \textit{infinitely divisible}
if it has infinitely divisible distribution on $\mbb{R}_+$.

We next give a characterization for the class of infinitely divisible
probability measures on $\mbb{R}_+$. Write $\psi\in \mcr{I}$ if
$\lambda\mapsto \psi(\lambda)$ is a positive function on $[0,\infty)$
with the representation
 \beqlb\label{ch01.2.3}
\psi(\lambda) = h\lambda + \int_0^\infty (1-\e^{-\lambda u}) l(\d u),
 \eeqlb
where $h\ge 0$ and $(1\land u)l(\d u)$ is a finite measure on
$(0,\infty)$.

 \bgproposition\label{tch01.2.1}
The pair $(h,l)$ in \eqref{ch01.2.3} is uniquely determined by the
function $\psi\in \mcr{I}$. \edproposition

 \begin{proof}
Suppose that $\psi$ can also be represented by \eqref{ch01.2.3} with
$(h,l)$ replaced by $(h^\prime,l^\prime)$. For $\lambda> 0$ and
$\theta\ge 0$, we can evaluate $\psi(\lambda+\theta) - \psi(\theta)$ with
the two representations and get
 \beqnn
h\lambda + \int_0^\infty\big(1-\e^{-\lambda u}\big) \e^{-\theta u} l(\d
u)
 =
h^\prime\lambda + \int_0^\infty\big(1-\e^{-\lambda u}\big) \e^{-\theta u}
l^\prime(\d u).
 \eeqnn
By letting $\theta \to \infty$ we get $h = h^\prime$, and so $l(\d u) =
l^\prime(\d u)$ by Theorem~\ref{tch01.1.5}. \qed\end{proof}

 \bgtheorem\label{tch01.2.2}
Suppose that $\psi$ is a continuous function on $[0,\infty)$. If there is
a sequence $\{\psi_n\}\subset \mcr{I}$ such that $\psi(\lambda) =
\lim_{n\to \infty} \psi_n(\lambda)$ for all $\lambda\ge 0$, then $\psi\in
\mcr{I}$. \edtheorem

 \begin{proof}
Suppose that $\psi_n\in \mcr{I}$ is given by \eqref{ch01.2.3} with
$(h,l)$ replaced by $(h_n,l_n)$. We can define a finite measure $F_n$ on
$\bar{\mbb{R}}_+$ by setting $F_n(\{0\}) = h_n$, $F_n(\{\infty\}) = 0$
and $F_n(\d u) = (1-\e^{-u}) l_n(\d u)$ for $0<u<\infty$. For $\lambda>0$
let
 \beqlb\label{ch01.2.4}
\xi(u,\lambda)
 =
\left\{\begin{array}{ll}
 (1-\e^{-u})^{-1}(1-\e^{-u\lambda}) \quad \ar\mbox{if $0<u<\infty$,}\cr
 \lambda \quad \ar\mbox{if $u=0$,} \cr
 1 \quad \ar\mbox{if $u=\infty$.}
\end{array}\right.
 \eeqlb
Then we have
 \beqnn
\psi_n(\lambda) = \int_{\bar{\mbb{R}}_+} \xi(u,\lambda) F_n(\d u), \qquad
\lambda>0.
 \eeqnn
It is evident that $\{F_n(\bar{\mbb{R}}_+)\}$ is a bounded sequence. Take
any subsequence $\{F_{n_k}\}\subset \{F_n\}$ such that $\lim_{k\to
\infty} F_{n_k} = F$ weakly for a finite measure $F$ on
$\bar{\mbb{R}}_+$. Since $u\mapsto \xi(u,\lambda)$ is continuous on
$\bar{\mbb{R}}_+$, we have
 \beqnn
\psi(\lambda) = \int_{\bar{\mbb{R}}_+} \xi(u,\lambda) F(\d u), \qquad
\lambda>0.
 \eeqnn
Observe also that $\lim_{n\to \infty} \psi(1/n) = \psi(0) = 0$ implies
$F(\{\infty\}) =0$. Then the desired conclusion follows by a change of
the integration variable. \qed\end{proof}

 \bgtheorem\label{tch01.2.3}
The relation $\psi = -\log L_\mu$ establishes a one-to-one correspondence
between the functions $\psi\in \mcr{I}$ and infinitely divisible
probability measures $\mu$ on $\mbb{R}_+$. \edtheorem

 \begin{proof}
Suppose that $\psi\in \mcr{I}$ is given by \eqref{ch01.2.3}. Let $N$ be a
Poisson random measure on $(0,\infty)$ with intensity $l(\d u)$ and let
 \beqnn
\xi = h + \int_0^\infty x N(\d x).
 \eeqnn
By Theorem~\ref{tch01.1.7} for any $\lambda\ge 0$ we have
 \beqnn
\mbf{E}\,\e^{-\lambda\xi}
 =
\exp\Big\{-h\lambda - \int_0^\infty \big(1-\e^{-\lambda u}\big) l(\d
u)\Big\}.
 \eeqnn
Then $\psi = -\log L_\mu$ for a probability measure $\mu$ on $\mbb{R}_+$.
Similarly, for each integer $n\ge 1$ there is a probability measure
$\mu_n$ on $\mbb{R}_+$ so that $\psi/n = -\log L_{\mu_n}$. It is easy to
see that $\mu_n^{*n} = \mu$. That gives the infinite divisibility of
$\mu$. Conversely, suppose that $\psi = -\log L_\mu$ for an infinitely
divisible probability measure $\mu$ on $\mbb{R}_+$. For $n\ge 1$ let
$\mu_n$ be the $n$-th root of $\mu$. Then
 \beqnn
\psi(\lambda) = \lim_{n\to \infty} n[1-\e^{-n^{-1}\psi(\lambda)}]
 =
\lim_{n\to \infty} \int_0^\infty \big(1-\e^{-\lambda\x}\big)
n\mu_n(\d\x).
 \eeqnn
By Theorem~\ref{tch01.2.2} we have $\psi\in \mcr{I}$. \qed\end{proof}

The above theorem gives a complete characterization of infinitely
divisible probability measures on $\mbb{R}_+$. We write $\mu = I(h,l)$ if
$\mu$ is an infinitely divisible probability measure on $\mbb{R}_+$ with
$\psi := -\log L_\mu$ given by \eqref{ch01.2.3}.

 \bgtheorem\label{tch01.2.4}
If $\psi_1,\psi_2\in \mcr{I}$, then $\psi_1\circ \psi_2\in \mcr{I}$.
\edtheorem

 \begin{proof}
For every $x\ge 0$ we clearly have $x\psi_2\in\mcr{I}$, so there is an
infinitely divisible probability measure $\nu_x$ on $\mbb{R}_+$
satisfying $-\log L_{\nu_x} = x\psi_2$. By a monotone class argument one
can see $\nu_x(\d y)$ is a probability kernel on $\mbb{R}_+$. Let $\mu$
be the infinitely divisible probability measure on $\mbb{R}_+$ with
$-\log L_\mu = \psi_1$ and define
 \beqnn
\eta(\d y) = \int_0^\infty \mu(\d x)\nu_x(\d y), \qquad y\ge 0.
 \eeqnn
It is not hard to show that $-\log L_\eta= \psi_1\circ \psi_2$. By the
same reasoning, for each integer $n\ge1$ there is a probability measure
$\eta_n$ such that $-\log L_{\eta_n} = n^{-1}\psi_1\circ \psi_2$. Then
$\eta = \eta_n^{*n}$ and hence $\eta$ is infinitely divisible. By
Theorem~\ref{tch01.2.3} we conclude that $\psi_1\circ \psi_2\in \mcr{I}$.
\qed\end{proof}

 \bexample\label{ech01.2.1}
Let $b>0$ and $\alpha>0$. The \index{Gamma distribution} \textit{Gamma
distribution} $\gamma$ on $\mbb{R}_+$ with parameters $(b,\alpha)$ is
defined by
 \beqnn
\gamma(B)
 =
\frac{\alpha^b}{\Gamma(b)}\int_Bx^{b-1}\e^{-\alpha x}\d x, \qquad B\in
\mcr{B}(\mbb{R}_+),
 \eeqnn
which reduces to the \index{exponential distribution} \textit{exponential
distribution} when $b=1$. The Laplace transform of $\gamma$ is
 \beqnn
L_\gamma(\lambda)
 =
\Big(\frac{\alpha}{\alpha+\lambda}\Big)^b, \qquad \lambda\ge 0.
 \eeqnn
It is easily seen that $\gamma$ is infinitely divisible and its $n$-th
root is the Gamma distribution with parameters $(b/n,\alpha)$. \eexample

 \bexample\label{ech01.2.2}
For $c>0$ and $0<\alpha<1$ the function $\lambda\mapsto c\lambda^\alpha$
admits the representation \eqref{ch01.2.3}. Indeed, it is simple to show
 \beqlb\label{ch01.2.5}
\lambda^\alpha
 =
\frac{\alpha}{\Gamma(1-\alpha)}\int_0^\infty (1-\e^{-\lambda u}) \frac{\d
u}{u^{1+\alpha}}, \qquad \lambda\ge 0.
 \eeqlb
The infinitely divisible probability measure $\nu$ on $\mbb{R}_+$
satisfying $-\log L_\nu(\lambda) = c\lambda^\alpha$ is known as the
\index{one-sided stable distribution (with index $0<\alpha<1$)}
\textit{one-sided stable distribution} with index $0<\alpha<1$. This
distribution does not charge zero and is absolutely continuous with
respect to the Lebesgue measure on $(0,\infty)$ with continuous density.
For $\alpha = 1/2$ it has density
 \beqnn
q(x)
 :=
\frac{c}{2\sqrt{\pi}}x^{-3/2}\e^{-c^2/4x}, \qquad x>0.
 \eeqnn
For a general index the density can be given using an infinite series;
see, e.g., Sato (1999, p.88). \eexample

\section{L\'{e}vy--Khintchine type representations}\label{sch01.3}

 \setcounter{equation}{0}

In this section, we give some criteria for continuous functions on
$[0,\infty)$ to have L\'{e}vy--Khintchine type representations. The
results are useful in the study of high-density limits of discrete
branching processes. For $u\ge 0$ and $\lambda\ge 0$ let
 \beqnn
\xi_n(u, \lambda) = \e^{-\lambda u} -1 - (1+u^n)^{-1}
\sum_{i=1}^{n-1}\frac{(-\lambda u)^i}{i!}, \qquad n=1,2,\ldots.
 \eeqnn
We are interested in functions $\phi$ on $[0,\infty)$ with the
representation
 \beqlb\label{ch01.3.1}
\phi(\lambda) = \sum_{i=0}^{n-1}a_i\lambda^i + \int_0^\infty \xi_n(u,
\lambda)(1-\e^{-u})^{-n} G(\d u), \quad \lambda\ge 0,
 \eeqlb
where $n\ge 1$ is an integer, $\{a_0,\ldots,a_{n-1}\}$ is a set of
constants and $G(\d u)$ is a finite measure on $\mbb{R}_+$. The value at
$u=0$ of the integrand in \eqref{ch01.3.1} is defined by continuity as
$(-\lambda)^n/n!$. The following theorem was proved in Li (1991, 2011):

 \bgtheorem\label{tch01.3.1}
A continuous real function $\phi$ on $[0,\infty)$ has the representation
\eqref{ch01.3.1} if and only if for every $c\ge 0$ the function
 \beqlb\label{ch01.3.2}
\theta_c(\lambda) := (-1)^n\Delta_c^n\phi(\lambda), \qquad \lambda\ge 0
 \eeqlb
is the Laplace transform of a finite measure on $\mbb{R}_+$. \edtheorem

Based on the above theorem we can give canonical representations for the
limit functions of some sequences involving probability generating
functions. Let $\{\alpha_k\}$ be a sequence of positive numbers and let
$\{g_k\}$ be a sequence of probability generating functions, that is,
 \beqnn
g_k(z) = \sum_{i=0}^\infty p_{ki} z^i, \qquad |z|\le 1,
 \eeqnn
where $p_{ki}\ge 0$ and $\sum_{i=0}^\infty p_{ki} = 1$. We first consider
the sequence of functions $\{\psi_k\}$ defined by
 \beqlb\label{ch01.3.3}
\psi_k(\lambda) = \alpha_k [1-g_k(1-\lambda/k)], \qquad 0\le \lambda\le
k.
 \eeqlb

 \bgtheorem\label{tch01.3.2}
If the sequence $\{\psi_k\}$ defined by \eqref{ch01.3.3} converges to a
continuous real function $\psi$ on $[0,\infty)$, then the limit function
belongs to the class $\mcr{I}$ defined by \eqref{ch01.2.3}. \edtheorem

 \begin{proof}
For any $c,\lambda\ge 0$ and sufficiently large $k\ge 1$ we have
 \beqnn
\Delta_c\psi_k(\lambda)
 =
-\alpha_k\Delta_c g_k(1-\cdot/k)(\lambda).
 \eeqnn
Since for each integer $i\ge 1$ the $i$-th derivative $g_k^{(i)}$ is a
power series with positive coefficients, we have
 \beqnn
(-1)^i\frac{\d^i}{\d\lambda^i}\Delta_c\psi_k(\lambda)
 =
- k^{-i}\alpha_k\Delta_cg_k^{(i)}(1-\cdot/k)(\lambda)\ge 0.
 \eeqnn
By the mean-value theorem, one sees inductively $(-1)^i \Delta_h^i
\Delta_c \psi_k(\lambda)\ge 0$. Letting $k\to \infty$ we obtain $(-1)^i
\Delta_h^i \Delta_c \psi(\lambda)\ge 0$. Then $\Delta_c\psi(\lambda)$ is
a completely monotone function of $\lambda\ge 0$, so by
Theorem~\ref{tch01.1.4} it is the Laplace transform of a finite measure
on $\mbb{R}_+$. Since $\psi(0) = \lim_{k\to \infty}\psi_k(0) = 0$, by
Theorem~\ref{tch01.3.1} there is a finite measure $F$ on $\mbb{R}_+$ so
that
 \beqnn
\psi(\lambda) = \int_0^\infty (1 - \e^{-\lambda u}) (1-\e^{-u})^{-1} F(\d
u),
 \eeqnn
where the value of the integrand at $u=0$ is defined as $\lambda$ by
continuity. Then \eqref{ch01.2.3} follows with $\beta = F(\{0\})$ and
$n(\d u) = (1-\e^{-u})^{-1} F(\d u)$ for $u>0$. \qed\end{proof}

 \bexample\label{ech01.3.1}
Suppose that $g$ is a probability generating function so that $\beta :=
g^\prime(1-)< \infty$. Let $\alpha_k = k$ and $g_k(z) = g(z)$. Then the
sequence $\psi_k(\lambda)$ defined by \eqref{ch01.3.3} converges to
$\beta\lambda$ as $k\to \infty$. \eexample

 \bexample\label{ech01.3.2}
For any $0< \alpha\le 1$ the function $\psi(\lambda) = \lambda^\alpha$
has the representation \eqref{ch01.2.3}. For $\alpha=1$ that is trivial,
and for $0< \alpha< 1$ that follows from \eqref{ch01.2.5}. Let
$\psi_k(\lambda)$ be defined by \eqref{ch01.3.3} with $\alpha_k =
k^\alpha$ and $g_k(z) = 1-(1-z)^\alpha$. Then $\psi_k(\lambda) =
\lambda^\alpha$ for $0\le \lambda\le k$. \eexample

In the study of limit theorems of branching models, we shall also need to
consider the limit of another function sequence defined as follows. Let
$\{\alpha_k\}$ and $\{g_k\}$ be given as above and let
 \beqlb\label{ch01.3.4}
\phi_k(\lambda) = \alpha_k [g_k(1-\lambda/k) - (1-\lambda/k)], \qquad
0\le \lambda\le k.
 \eeqlb

 \bgtheorem\label{tch01.3.3}
If the sequence $\{\phi_k\}$ defined by \eqref{ch01.3.4} converges to a
continuous real function $\phi$ on $[0,\infty)$, then the limit function
has the representation
 \beqlb\label{ch01.3.5}
\phi(\lambda) = a\lambda + c\lambda^2 + \int_0^\infty \Big(\e^{-\lambda
u} - 1 + \frac{\lambda u}{1+u^2}\Big) m(\d u),
 \eeqlb
where $c\ge 0$ and $a$ are constants, and $m(\d u)$ is a $\sigma$-finite
measure on $(0,\infty)$ satisfying
 \beqlb\label{ch01.3.6}
\int_0^\infty (1\land u^2) m(\d u)<\infty.
 \eeqlb
\edtheorem

 \begin{proof}
Since $\phi(0) = \lim_{k\to \infty}\phi_k(0) = 0$, arguing as in the
proof of Theorem~\ref{tch01.3.2} we see that $\phi$ has the
representation \eqref{ch01.3.1} with $n=2$ and $a_0=0$, which can be
rewritten into the equivalent form \eqref{ch01.3.5}. \qed\end{proof}

Now let us consider a special case of the function $\phi$ given by
\eqref{ch01.3.5}. Observe that if the measure $m(\d u)$ satisfies the
integrability condition
 \beqlb\label{ch01.3.7}
\int_0^\infty (u\land u^2)\, m(\d u)<\infty,
 \eeqlb
we have
 \beqlb\label{ch01.3.8}
\phi(\lambda) = b\lambda + c\lambda^2 + \int_0^\infty \big(\e^{-\lambda
u} -1 + \lambda u\big) m(\d u),
 \eeqlb
where
 \beqnn
b = a - \int_0^\infty \frac{u^3}{1+u^2} m(\d u).
 \eeqnn

 \bgproposition\label{tch01.3.4}
A function $\phi$ with the representation \eqref{ch01.3.5} is locally
Lipschitz if and only if \eqref{ch01.3.7} holds. \edproposition

 \begin{proof}
For computational convenience we first rewrite \eqref{ch01.3.5} as
 \beqlb\label{ch01.3.9}
\phi(\lambda) = b_1\lambda + c\lambda^2 + \int_0^\infty \big(\e^{-\lambda
u} - 1 + \lambda u1_{\{u\le 1\}}\big) m(\d u),
 \eeqlb
where
 \beqnn
b_1 := a + \int_0^\infty \Big(\frac{u}{1+u^2} - u1_{\{u\le 1\}}\Big) m(\d
u).
 \eeqnn
By applying dominated convergence to \eqref{ch01.3.9}, for each
$\lambda>0$ we have
 \beqnn
\phi^\prime(\lambda) = b_1 + 2c\lambda + \int_{(0,1]} u \big(1 -
\e^{-\lambda u}\big) m(\d u) - \int_{(1,\infty)} u\e^{-\lambda u} m(\d
u).
 \eeqnn
Then we use monotone convergence to the two integrals to get
 \beqnn
\phi^\prime(0+) = b_1 - \int_{(1,\infty)} u m(\d u).
 \eeqnn
If $\phi$ is locally Lipschitz, we have $\phi^\prime(0+)> -\infty$ and
the integral on the right-hand side is finite. This together with
\eqref{ch01.3.6} implies \eqref{ch01.3.7}. Conversely, if
\eqref{ch01.3.7} holds, then $\phi^\prime$ is bounded on each bounded
interval and so $\phi$ is locally Lipschitz. \qed\end{proof}

 \bgcorollary\label{tch01.3.5}
If the sequence $\{\phi_k\}$ defined by \eqref{ch01.3.4} is uniformly
Lipschitz on each bounded interval and $\phi_k(\lambda)\to \phi(\lambda)$
for all $\lambda\ge 0$ as $k\to \infty$, then the limit function has the
representation \eqref{ch01.3.8}. \edcorollary

 \bexample\label{ech01.3.3}
Suppose that $g$ is a probability generating function so that
$g^\prime(1-) = 1$ and $c := g^{\prime\prime}(1-)/2< \infty$. Let
$\alpha_k = k^2$ and $g_k(z) = g(z)$. By Taylor's expansion it is easy to
show that the sequence $\phi_k(\lambda)$ defined by \eqref{ch01.3.4}
converges to $c\lambda^2$ as $k\to\infty$. \eexample

 \bexample\label{ech01.3.4}
For $0< \alpha< 1$ the function $\phi(\lambda) = -\lambda^\alpha$ has the
representation \eqref{ch01.3.5}. That follows from \eqref{ch01.2.5} as we
notice
 \beqnn
\int_0^\infty \Big(\frac{u}{1+u^2}\Big)\frac{\d u}{u^{1+\alpha}}
 =
\int_0^\infty \Big(\frac{1}{1+u^2}\Big)\frac{\d u}{u^\alpha} <\infty.
 \eeqnn
The function is the limit of the sequence $\phi_k(\lambda)$ defined by
\eqref{ch01.3.4} with $\alpha_k = k^\alpha$ and $g_k(z) =
1-(1-z)^\alpha$. \eexample

 \bexample\label{ech01.3.5}
For any $1\le \alpha\le 2$ the function $\phi(\lambda) = \lambda^\alpha$
can be represented in the form of \eqref{ch01.3.8}. In particular, for
$1<\alpha<2$ we have
 \beqnn
\lambda^\alpha
 =
\frac{\alpha(\alpha-1)}{\Gamma(2-\alpha)}\int_0^\infty (\e^{-\lambda u} -
1 + \lambda u) \frac{\d u}{u^{1+\alpha}}, \qquad \lambda\ge 0.
 \eeqnn
Let $\phi_k(\lambda)$ be defined by \eqref{ch01.3.4} with $\alpha_k =
\alpha k^\alpha$ and $g_k(z) = z + \alpha^{-1} (1-z)^\alpha$. Then
$\phi_k(\lambda) = \lambda^\alpha$ for $0\le \lambda\le k$. \eexample


\chapter{Continuous-state branching processes}\label{ch02}

In this chapter, we first give a construction of CB-processes as the
scaling limits of discrete Galton--Watson branching processes. This
approach also gives the interpretations of the CB-processes. We shall
study some basic properties of the CB-processes. In particular, some
conditional limit theorems will be given. We also give a reconstruction
of the sample paths of the CB-processes in terms of excursions.

\section{Construction by scaling limits}\label{sch02.1}

 \setcounter{equation}{0}

Suppose that $\{\xi_{n,i}: n,i=1,2,\ldots\}$ is a family of positive
integer-valued i.i.d.\ random variables with distribution given by the
probability generating function $g$. Given the positive integer $x(0)=m$,
we define inductively
 \beqlb\label{ch02.1.1}
x(n) = \sum_{i=1}^{x(n-1)}\xi_{n,i}, \qquad n=1,2,\ldots.
 \eeqlb
It is easy to show that $\{x(n): n = 0,1,2,\ldots\}$ is a discrete-time
positive integer-valued Markov chain with transition matrix $P(i,j)$
defined by
 \beqlb\label{ch02.1.2}
\sum^\infty_{j=0} P(i,j)z^j = g(z)^i, \qquad i=0,1,2,\ldots,~ |z|\le 1.
 \eeqlb
The random variable $x(n)$ can be thought of as the number of individuals
in generation $n\ge 0$ of an evolving particle system. After one unit
time, each of the $x(n)$ particles splits independently of others into a
random number of offspring according to the distribution given by $g$;
see, e.g., Athreya and Ney (1972). For $n\ge 0$ the $n$-step transition
matrix $P^n(i,j)$ is determined by
 \beqlb\label{ch02.1.3}
\sum^\infty_{j=0} P^n(i,j)z^j = g^n(z)^i, \qquad i=0,1,2,\ldots,~ |z|\le
1,
 \eeqlb
where $g^n(z)$ is defined by $g^n(z) = g(g^{n-1}(z))$ successively with
$g^0(z) = z$. We call any positive integer-valued Markov chain with
transition matrix given by \eqref{ch02.1.2} or \eqref{ch02.1.3} a
\index{Galton--Watson branching process} \textit{Galton--Watson branching
process} \index{GW-process} (GW-process). If $g^\prime(1-)<\infty$, the
first moment of the discrete distribution $\{P^n(i,j); j=0,1,2,\ldots\}$
is given by
 \beqlb\label{ch02.1.4}
\sum^\infty_{j=1} jP^n(i,j) = ig^\prime(1-)^n,
 \eeqlb
which can be obtained by differentiating both sides of \eqref{ch02.1.3}.

Now suppose we have a sequence of GW-processes $\{x_k(n): n\ge 0\}$ with
offspring distribution given by the sequence of probability generating
functions $\{g_k\}$. Let $z_k(n) = x_k(n)/k$ for $n\ge 0$. Then
$\{z_k(n): n\ge 0\}$ is a Markov chain with state space $E_k :=
\{0,1/k,2/k,\ldots\}$ and $n$-step transition probability $P_k^n(x,dy)$
determined by
 \beqlb\label{ch02.1.5}
\int_{E_k}\e^{-\lambda y}P_k^n(x,\d y)
 =
g_k^n(\e^{-\lambda/k})^{kx}, \qquad \lambda\ge 0.
 \eeqlb
Suppose that $\{\gamma_k\}$ is a positive sequence so that $\gamma_k\to
\infty$ as $k\to \infty$. Let $[\gamma_kt]$ denote the integer part of
$\gamma_kt\ge 0$. We are interested in the asymptotic behavior of the
sequence of continuous time processes $\{z_k([\gamma_kt]): t\ge 0\}$. By
\eqref{ch02.1.5} we have
 \beqlb\label{ch02.1.6}
\int_{E_k}\e^{-\lambda y}P_k^{[\gamma_kt]}(x,\d y)
 =
\exp\{-xv_k(t,\lambda)\},
 \eeqlb
where
 \beqlb\label{ch02.1.7}
v_k(t,\lambda)
 =
-k\log g_k^{[\gamma_kt]}(\e^{-\lambda/k}), \qquad \lambda\ge 0.
 \eeqlb
Clearly, if $z_k(0) = x\in E_k$, then the probability
$P_k^{[\gamma_kt]}(x,\cdot)$ gives the distribution of $z_k([\gamma_kt])$
on $\mbb{R}_+$. Let us consider the function sequences
 \beqlb\label{ch02.1.8}
G_k(z)
 =
k\gamma_k[g_k(\e^{-z/k})-\e^{-z/k}], \qquad z\ge 0,
 \eeqlb
and
 \beqlb\label{ch02.1.9}
\phi_k(z) = k\gamma_k[g_k(1-z/k) - (1-z/k)], \quad 0\le z\le k.
 \eeqlb

 \bgproposition\label{tch02.1.1}
The sequence $\{G_k\}$ is uniformly Lipschitz on each bounded interval if
and only if so is $\{\phi_k\}$. In this case, we have $\lim_{k\to\infty}
|\phi_k(z)-G_k(z)| = 0$ uniformly on each bounded interval.
\edproposition

 \begin{proof}
From \eqref{ch02.1.8} and \eqref{ch02.1.9} it is simple to see that
 \beqlb\label{ch02.1.10}
G^\prime_k(z)
 =
\gamma_k\e^{-z/k}[1-g^\prime_k(\e^{-z/k})], \qquad z\ge 0,
 \eeqlb
and
 \beqlb\label{ch02.1.11}
\phi^\prime_k(z)
 =
\gamma_k[1-g^\prime_k(1-z/k)], \qquad 0\le z\le k.
 \eeqlb
Clearly, the sequence $\{G_k^\prime\}$ is uniformly bounded on each
bounded interval if and only if so is $\{\phi_k^\prime\}$. Then the first
assertion is immediate. We next assume $\{G_k\}$ is uniformly Lipschitz
on each bounded interval. Let $a\ge 0$. By the mean-value theorem, for
$k\ge a$ and $0\le z\le a$ we have
 \beqlb\label{ch02.1.12}
G_k(z) - \phi_k(z)
 \ar=\ar
k\gamma_k\big[g_k(\e^{-z/k}) - g_k(1-z/k) - \e^{-z/k} + (1 - z/k)\big]
\quad \cr
 \ar=\ar
k\gamma_k[g^\prime_k(\eta_k)-1](\e^{-z/k}-1+z/k),
 \eeqlb
where
 \beqnn
1-a/k\le 1-z/k\le \eta_k\le \e^{-z/k}\le 1.
 \eeqnn
Choose $k_0\ge a$ so that $\e^{-2a/k_0}\le 1-a/k_0$. Then $\e^{-2a/k}\le
1-a/k$ for $k\ge k_0$ and hence
 \beqnn
\gamma_k |g^\prime_k(\eta_k)-1|
 \le
\sup_{0\le z\le 2a} \gamma_k|g^\prime_k(\e^{-z/k})-1|, \qquad k\ge k_0.
 \eeqnn
Since $\{G_k\}$ is uniformly Lipschitz on each bounded interval, the
sequence \eqref{ch02.1.10} is uniformly bounded on $[0,2a]$. Then
$\{\gamma_k |g^\prime_k(\eta_k)-1|: k\ge k_0\}$ is a bounded sequence.
Now the desired result follows from \eqref{ch02.1.12}. \qed\end{proof}

By the above proposition, if either $\{G_k\}$ or $\{\phi_k\}$ is
uniformly Lipschitz on each bounded interval, then they converge or
diverge simultaneously and in the convergent case they have the same
limit. For the convenience of statement of the results, we formulate the
following conditions:

 \bgcondition\label{tch02.1.2}
The sequence $\{G_k\}$ is uniformly Lipschitz on $[0,a]$ for every $a\ge
0$ and there is a function $\phi$ on $[0,\infty)$ so that $G_k(z)\to
\phi(z)$ uniformly on $[0,a]$ for every $a\ge 0$ as $k\to \infty$.
\edcondition

 \bgproposition\label{tch02.1.3}
Suppose that Condition~\ref{tch02.1.2} is satisfied. Then the function
$\phi$ has representation
 \beqlb\label{ch02.1.13}
\phi(z) = bz + cz^2 + \int_0^\infty\big(\e^{-zu}-1+zu\big)m(\d u), \quad
z\ge 0,
 \eeqlb
where $c\ge 0$ and $b$ are constants and $(u\land u^2)m(\d u)$ is a
finite measure on $(0,\infty)$. \edproposition

 \begin{proof}
By Proposition~\ref{tch02.1.1}, the sequence $\{\phi_k\}$ is uniformly
Lipschitz on $[0,a]$ and $\phi_k(z)\to \phi(z)$ uniformly on $[0,a]$ for
every $a\ge 0$. Then the result follows by Corollary~\ref{tch01.3.5}.
\qed\end{proof}

 \bgproposition\label{tch02.1.4}
For any function $\phi$ with representation \eqref{ch02.1.13} there is a
sequence $\{G_k\}$ in the form of \eqref{ch02.1.8} satisfying
Condition~\ref{tch02.1.2}. \edproposition

 \begin{proof}
By Proposition~\ref{tch02.1.1} it suffices to construct a sequence
$\{\phi_k\}$ in the form of \eqref{ch02.1.9} uniformly Lipschitz on
$[0,a]$ and $\phi_k(z)\to \phi(z)$ uniformly on $[0,a]$ for every $a\ge
0$. To simplify the formulations we decompose the function $\phi$ into
two parts. Let $\phi_0(z) = \phi(z) - bz$. We first define
 \beqnn
\gamma_{0,k} = (1+2c)k + \int_0^\infty u(1-\e^{-ku}) m(\d u)
 \eeqnn
and
 \beqnn
g_{0,k}(z) = z + k^{-1}\gamma_{0,k}^{-1} \phi_0(k(1-z)), \qquad |z|\le 1.
 \eeqnn
It is easy to see that $z\mapsto g_{0,k}(z)$ is an analytic function in
$(-1,1)$ satisfying $g_{0,k}(1) = 1$ and
 \beqnn
\frac{\d^n}{\d z^n}g_{0,k}(0)\ge 0, \qquad n\ge 0.
 \eeqnn
Therefore $g_{0,k}(\cdot)$ is a probability generating function. Let
$\phi_{0,k}$ be defined by \eqref{ch02.1.9} with $(\gamma_k,g_k)$
replaced by $(\gamma_{0,k}, g_{0,k})$. Then $\phi_{0,k}(z) = \phi_0(z)$
for $0\le z\le k$. That completes the proof if $b=0$. In the case $b\neq
0$, we set
 \beqnn
g_{1,k}(z) = \frac{1}{2}\Big(1+\frac{b}{|b|}\Big) +
\frac{1}{2}\Big(1-\frac{b}{|b|}\Big)z^2.
 \eeqnn
Let $\gamma_{1,k} = |b|$ and let $\phi_{1,k}(z)$ be defined by
\eqref{ch02.1.9} with $(\gamma_k,g_k)$ replaced by $(\gamma_{1,k},
g_{1,k})$. Thus we have
 \beqnn
\phi_{1,k}(z) = bz + \frac{1}{2k}(|b|-b)z^2.
 \eeqnn
Finally, let $\gamma_k = \gamma_{0,k} + \gamma_{1,k}$ and $g_k =
\gamma_k^{-1}(\gamma_{0,k}g_{0,k} + \gamma_{1,k}g_{1,k})$. Then the
sequence $\phi_k(z)$ defined by \eqref{ch02.1.9} is equal to
$\phi_{0,k}(z) + \phi_{1,k}(z)$ which satisfies the required condition.
\qed\end{proof}

 \bglemma\label{tch02.1.5}
Suppose that the sequence $\{G_k\}$ defined by \eqref{ch02.1.8} is
uniformly Lipschitz on $[0,1]$. Then there are constants $B,N\ge 0$ such
that $v_k(t,\lambda)\le \lambda\e^{Bt}$ for every $t,\lambda\ge 0$ and
$k\ge N$. \edlemma

 \begin{proof}
Let $b_k := G_k^\prime(0+)$ for $k\ge 1$. Since $\{G_k\}$ is uniformly
Lipschitz on $[0,1]$, the sequence $\{b_k\}$ is bounded. {From}
\eqref{ch02.1.8} we have $b_k = \gamma_k[1-g_k^\prime(1-)]$. By
\eqref{ch02.1.4} it is not hard to obtain
 \beqnn
\int_{E_k} y P_k^{[\gamma_kt]}(x,\d y)
 =
xg_k^\prime(1-)^{[\gamma_kt]}
 =
x\Big(1-\frac{b_k}{\gamma_k}\Big)^{[\gamma_kt]}.
 \eeqnn
Let $B\ge 0$ be a constant such that $2|b_k|\le B$ for all $k\ge 1$.
Since $\gamma_k\to \infty$ as $k\to \infty$, there is $N\ge 1$ so that
 \beqnn
0\le \Big(1-\frac{b_k}{\gamma_k}\Big)^{\frac{\gamma_k}{B}}
 \le
\Big(1+\frac{B}{2\gamma_k}\Big)^{\frac{\gamma_k}{B}}
 \le
\e, \qquad k\ge N.
 \eeqnn
It follows that, for $t\ge 0$ and $k\ge N$,
 \beqlb\label{ch02.1.14}
\int_{E_k} y P_k^{[\gamma_kt]}(x,\d y)
 \le
x\exp\Big\{\frac{B}{\gamma_k}[\gamma_kt]\Big\}
 \le
x\e^{Bt}.
 \eeqlb
Then the desired estimate follows from \eqref{ch02.1.6},
\eqref{ch02.1.14} and Jensen's inequality.
 \qed\end{proof}

 \bgtheorem\label{tch02.1.6}
Suppose that Condition~\ref{tch02.1.2} holds. Then for every $a\ge 0$ we
have $v_k(t,\lambda)\to$ some $v_t(\lambda)$ uniformly on $[0,a]^2$ as
$k\to \infty$ and the limit function solves the integral equation
 \beqlb\label{ch02.1.15}
v_t(\lambda) = \lambda - \int_0^t\phi(v_s(\lambda))\d s, \qquad
\lambda,t\ge 0.
 \eeqlb
\edtheorem

 \begin{proof}
The following argument is a modification of that of Aliev and Shchurenkov
(1982) and Aliev (1985). For any $n\ge 0$ we may write
 \beqnn
\log g_k^{n+1}(\e^{-\lambda/k})
 \ar=\ar
\log\big[g_k(g_k^n(\e^{-\lambda/k})) g_k^n(\e^{-\lambda/k})^{-1}\big] +
\log g_k^n(\e^{-\lambda/k}) \cr
 \ar=\ar
(k\gamma_k)^{-1}\bar{G}_k\big(-k\log g_k^n(\e^{-\lambda/k})\big) + \log
g_k^n(\e^{-\lambda/k}),
 \eeqnn
where
 \beqnn
\bar{G}_k(z)=k\gamma_k\log\big[g_k(\e^{-z/k})\e^{z/k}\big].
 \eeqnn
From this and \eqref{ch02.1.7} it follows that
 \beqnn
v_k(t+\gamma_k^{-1},\lambda)
 =
v_k(t,\lambda)-\gamma_k^{-1}\bar{G}_k(v_k(t,\lambda)).
 \eeqnn
By applying the above equation to $t=0,1/\gamma_k, \ldots,
([\gamma_kt]-1)/\gamma_k$ and adding the resulting equations we obtain
 \beqnn
v_k(t,\lambda)
 =
\lambda-\sum_{i=1}^{[\gamma_kt]}\gamma_k^{-1}\bar{G}_k(v_k(\gamma_k^{-1}(i-1),\lambda)).
 \eeqnn
Then we can write
 \beqlb\label{ch02.1.16}
v_k(t,\lambda)
 =
\lambda + \varepsilon_k(t,\lambda) - \int_0^t
\bar{G}_k(v_k(s,\lambda))ds,
 \eeqlb
where
 \beqnn
\varepsilon_k(t,\lambda)
 =
\big(t-\gamma_k^{-1}[\gamma_kt]\big)
\bar{G}_k\big(v_k(\gamma_k^{-1}[\gamma_kt],\lambda)\big).
 \eeqnn
It is not hard to see
 \beqnn
\bar{G}_k(z) = k\gamma_k\log\big[1+(k\gamma_k)^{-1}G_k(z)\e^{z/k}\big].
 \eeqnn
By Condition~\ref{tch02.1.2}, for any $0<\varepsilon\le 1$ we can enlarge
$N\ge 1$ so that
 \beqlb\label{ch02.1.17}
|\bar{G}_k(z)-\phi(z)|\le \varepsilon, \qquad 0\le z\le a\e^{Ba}, k\ge N.
 \eeqlb
It then follows that
 \beqlb\label{ch02.1.18}
|\varepsilon_k(t,\lambda)|
 \le
\gamma_k^{-1}M, \qquad 0\le t,\lambda\le a,
 \eeqlb
where
 \beqnn
M = 1+\sup_{0\le z\le a\e^{Ba}}|\phi(z)|.
 \eeqnn
For $n\ge k\ge N$ let
 \beqnn
K_{k,n}(t,\lambda)
 =
\sup_{0\le s\le t}|v_n(s,\lambda)-v_k(s,\lambda)|.
 \eeqnn
By \eqref{ch02.1.16}, \eqref{ch02.1.17} and \eqref{ch02.1.18} we obtain
 \beqnn
K_{k,n}(t,\lambda)
 \le
2(\gamma_k^{-1}M + \varepsilon a) + L\int_0^t K_{k,n}(s,\lambda) ds,
\qquad 0\le t,\lambda\le a,
 \eeqnn
where $L = \sup_{0\le s\le a\e^{Ba}}|\phi^\prime(z)|$. By Gronwall's
inequality,
 \beqnn
K_{k,n}(t,\lambda)
 \le
2(\gamma_k^{-1}M + \varepsilon a)\exp\{Lt\}, \qquad 0\le t,\lambda\le a.
 \eeqnn
Then $v_k(t,\lambda)\to$ some $v_t(\lambda)$ uniformly on $[0,a]^2$ as
$k\to \infty$ for every $a\ge 0$. From \eqref{ch02.1.16} we get
\eqref{ch02.1.15}. \qed\end{proof}

 \bgtheorem\label{tch02.1.7}
Suppose that $\phi$ is given by \eqref{ch02.1.13}. Then for any
$\lambda\ge 0$ there is a unique locally bounded positive solution
$t\mapsto v_t(\lambda)$ to \eqref{ch02.1.15}. Moreover, the solution
satisfies the semigroup property
 \beqlb\label{ch02.1.19}
v_{r+t}(\lambda) = v_r\circ v_t(\lambda) = v_r(v_t(\lambda)), \qquad
r,t,\lambda\ge 0.
 \eeqlb
\edtheorem

 \begin{proof}
By Propositions~\ref{tch02.1.4} and \ref{ch02.1.6} there is a locally
bounded positive solution to \eqref{ch02.1.15}. The proof of the
uniqueness of the solution is a standard application of Gronwall's
inequality. The relation \eqref{ch02.1.19} follows from the uniqueness of
the solution to \eqref{ch02.1.15}.
 \qed\end{proof}

 \bgtheorem\label{tch02.1.8}
Suppose that $\phi$ is given by \eqref{ch02.1.13}. Then there is a Feller
transition semigroup $(Q_t)_{t\ge 0}$ on $\mbb{R}_+$ defined by
 \beqlb\label{ch02.1.20}
\int_0^\infty \e^{-\lambda y} Q_t(x,\d y)
 =
\e^{-xv_t(\lambda)}, \qquad \lambda\ge 0,x\ge 0.
 \eeqlb
Moreover, if $E_k\ni x_k\to x\ge 0$, we have $P_k^{[\gamma_kt]}
(x_k,\cdot)\to Q_t(x,\cdot)$ weakly. \edtheorem

 \begin{proof}
By Proposition~\ref{tch02.1.4} and Theorems~\ref{tch01.1.2}
and~\ref{tch02.1.6} there is a probability kernel $Q_t(x,\d y)$ on
$\mbb{R}_+$ defined by \eqref{ch02.1.20}. Moreover, we have
$P_k^{[\gamma_kt]} (x_k,\cdot)\to Q_t(x,\cdot)$ weakly if $x_k\to x$. The
semigroup property of the family of kernels $(Q_t)_{t\ge 0}$ follows from
\eqref{ch02.1.19} and \eqref{ch02.1.20}. For $\lambda>0$ and $x\ge 0$ set
$e_\lambda(x) = \e^{-\lambda x}$. We denote by $D_1$ the linear span of
$\{e_\lambda: \lambda>0\}$. Clearly, the operator $Q_t$ preserves $D_1$
for every $t\ge 0$. By the continuity of $t\mapsto v_t(\lambda)$ it is
easy to show that $t\mapsto Q_te_\lambda(x)$ is continuous for
$\lambda>0$ and $x\ge 0$. Then $t\mapsto Q_tf(x)$ is continuous for every
$f\in D_1$ and $x\ge 0$. Let $C_0(\mbb{R}_+)$ be the space of continuous
functions on $\mbb{R}_+$ vanishing at infinity. By the Stone--Weierstrass
theorem, the set $D_1$ is uniformly dense in $C_0(\mbb{R}_+)$; see, e.g.,
Hewitt and Stromberg (1965, pp.98-99). Then each operator $Q_t$ preserves
$C_0(\mbb{R}_+)$ and $t\mapsto Q_tf(x)$ is continuous for every $f\in
C_0(\mbb{R}_+)$ and $x\ge 0$. That gives the Feller property of the
semigroup $(Q_t)_{t\ge 0}$. \qed\end{proof}

A Markov process is called a \index{continuous-state branching process}
\textit{continuous-state branching process} \index{CB-process}
(CB-process) with \index{branching mechanism} \textit{branching
mechanism} $\phi$ if it has transition semigroup $(Q_t)_{t\ge 0}$ defined
by \eqref{ch02.1.20}. It is simple to see that
 \beqlb\label{ch02.1.21}
Q_t(x_1+x_2,\cdot)
 =
Q_t(x_1,\cdot)*Q_t(x_2,\cdot), \qquad t,x_1,x_2\ge 0,
 \eeqlb
which is called the \index{branching property} \textit{branching
property} of $(Q_t)_{t\ge 0}$. The family of functions $(v_t)_{t\ge 0}$
is called the \index{cumulant semigroup} \textit{cumulant semigroup} of
the CB-process. By Theorem~\ref{tch02.1.8} the process has a
c\`{a}dl\`{a}g realization. Let $\Omega = D([0,\infty), \mbb{R}_+)$
denote the space of c\`{a}dl\`{a}g paths from $[0,\infty)$ to $\mbb{R}_+$
furnished with the Skorokhod topology. The following theorem gives an
interpretation of the CB-process as the approximation of the GW-process.

 \bgtheorem\label{tch02.1.9}
Suppose that Condition~\ref{tch02.1.2} holds. Let $\{x(t): t\ge 0\}$ be a
CB-process with transition semigroup $(Q_t)_{t\ge 0}$ defined by
\eqref{ch02.1.20}. If $z_k(0)$ converges to $x(0)$ in distribution, then
$\{z_k([\gamma_kt]): t\ge 0\}$ converges to $\{x(t): t\ge 0\}$ in
distribution on $D([0,\infty), \mbb{R}_+)$. \edtheorem

 \begin{proof}
For $\lambda>0$ and $x\ge 0$ set $e_\lambda(x) = \e^{-\lambda x}$. Let
$C_0(\mbb{R}_+)$ be the space of continuous functions on $\mbb{R}_+$
vanishing at infinity. By Theorem~\ref{tch02.1.6} it is easy to show
 \beqnn
\lim_{k\to\infty} \sup_{x\in E_k}\big|P_k^{[\gamma_kt]}e_\lambda(x) -
Q_te_\lambda(x)\big| = 0, \qquad \lambda>0.
 \eeqnn
Then the Stone--Weierstrass theorem implies
 \beqnn
\lim_{k\to\infty} \sup_{x\in E_k}\big|P_k^{[\gamma_kt]}f(x) -
Q_tf(x)\big| = 0, \qquad f\in C_0(\mbb{R}_+).
 \eeqnn
By Ethier and Kurtz (1986, p.226 and pp.233-234) we conclude that
$\{z_k([\gamma_kt]): t\ge 0\}$ converges to the CB-process $\{x(t): t\ge
0\}$ in distribution on $D([0,\infty), \mbb{R}_+)$. \qed\end{proof}

The convergence of rescaled Galton--Watson branching processes to
diffusion processes was first studied by Feller (1951). Ji\v{r}ina (1958)
introduced CB-processes in both discrete and continuous times. Lamperti
(1967a) showed that the continuous-time processes are weak limits of
rescaled Galton--Watson branching processes. We have followed Aliev and
Shchurenkov (1982) and Li (2006) in some of the above calculations; see
also Li (2011).

\section{Simple properties of CB-processes}\label{sch02.2}

 \setcounter{equation}{0}

In this section we prove some basic properties of CB-processes. Most of
the results presented here can be found in Grey (1974) and Li (2000). We
shall follow the treatments in Li (2011). Suppose that $\phi$ is a
branching mechanism defined by \eqref{ch02.1.13}. Then a CB-process has
transition semigroup $(Q_t)_{t\ge 0}$ defined by \eqref{ch02.1.15} and
\eqref{ch02.1.20}. It is easy to see for each $x\ge 0$, the probability
measure $Q_t(x,\cdot)$ is infinitely divisible. Then $(v_t)_{t\ge 0}$ can
be expressed canonically as
 \beqlb\label{ch02.2.1}
v_t(\lambda) = h_t\lambda + \int_0^\infty (1-\e^{-\lambda u}) l_t(\d u),
\qquad t\ge 0,\lambda\ge 0,
 \eeqlb
where $h_t\ge 0$ and $ul_t(\d u)$ is a finite measure on $(0,\infty)$.
{From} \eqref{ch02.1.15} we see that $t\mapsto v_t(\lambda)$ is first
continuous and then continuously differentiable. Moreover, we have the
backward differential equation:
 \beqlb\label{ch02.2.2}
\frac{\partial}{\partial t}v_t(\lambda)
 =
-\phi(v_t(\lambda)),
 \qquad
v_0(\lambda) = \lambda.
 \eeqlb
By \eqref{ch02.2.2} and the semigroup property $v_r\circ v_t = v_{r+t}$
for $r, t\ge 0$ we also have the forward differential equation
 \beqlb\label{ch02.2.3}
\frac{\partial}{\partial t}v_t(\lambda)
 =
-\phi(\lambda)\frac{\partial}{\partial \lambda}v_t(\lambda),
 \quad
v_0(\lambda)=\lambda.
 \eeqlb
By differentiating both sides of \eqref{ch02.1.15} it is easy to find
 \beqlb\label{ch02.2.4}
\frac{\partial}{\partial \lambda}v_t(0+)
 =
\e^{-bt}, \qquad t\ge 0,
 \eeqlb
which together with \eqref{ch02.1.20} yields
 \beqlb\label{ch02.2.5}
\int_0^\infty y Q_t(x,\d y)
 =
x\e^{-bt}, \qquad t\ge 0,x\ge 0.
 \eeqlb
We say the CB-process is \index{critical CB-process} \textit{critical},
\index{subcritical CB-process} \textit{subcritical} or
\index{supercritical CB-process} \textit{supercritical} according as
$b=0$, $\ge 0$ or $\le 0$.

 \bgproposition\label{tch02.2.1}
For every $t\ge 0$ the function $\lambda\mapsto v_t(\lambda)$ is strictly
increasing on $[0,\infty)$. \edproposition

 \begin{proof}
By the continuity of $t\mapsto v_t(\lambda)$, for any $\lambda_0>0$ there
is $t_0>0$ so that $v_t(\lambda_0)>0$ for $0\le t\le t_0$. Then
\eqref{ch02.1.20} implies $Q_t(x,\{0\})< 1$ for $x>0$ and $0\le t\le
t_0$, and so $\lambda\mapsto v_t(\lambda)$ is strictly increasing for
$0\le t\le t_0$. By the semigroup property of $(v_t)_{t\ge 0}$ we infer
$\lambda\mapsto v_t(\lambda)$ is strictly increasing for all $t\ge 0$.
\qed\end{proof}

 \bgcorollary\label{tch02.2.2}
The transition semigroup $(Q_t)_{t\ge 0}$ defined by \eqref{ch02.1.20} is
a Feller semigroup. \edcorollary

 \begin{proof}
By Proposition~\ref{tch02.2.1} for $t\ge 0$ and $\lambda>0$ we have
$v_t(\lambda)>0$. {From} \eqref{ch02.1.20} we see the operator $Q_t$ maps
$\{x\mapsto \e^{-\lambda x}: \lambda>0\}$ to itself. By the
Stone--Weierstrass theorem, the linear span of $\{x\mapsto \e^{-\lambda
x}: \lambda>0\}$ is dense in $C_0(\mbb{R}_+)$ in the supremum norm. Then
$Q_t$ maps $C_0(\mbb{R}_+)$ to itself. The Feller property of
$(Q_t)_{t\ge 0}$ follows by the continuity of $t\mapsto v_t(\lambda)$.
\qed\end{proof}

 \bgproposition\label{tch02.2.3}
Suppose that $\lambda> 0$ and $\phi(\lambda)\neq 0$. Then the equation
$\phi(z) = 0$ has no root between $\lambda$ and $v_t(\lambda)$. Moreover,
we have
 \beqlb\label{ch02.2.6}
\int_{v_t(\lambda)}^\lambda\phi(z)^{-1}\d z =t, \qquad t\ge 0.
 \eeqlb
\edproposition

 \begin{proof}
By \eqref{ch02.1.13} we see $\phi(0)=0$ and $z\mapsto \phi(z)$ is a
convex function. Since $\phi(\lambda)\neq 0$ for some $\lambda> 0$
according to the assumption, the equation $\phi(z) = 0$ has at most one
root in $(0,\infty)$. Suppose that $\lambda_0\ge 0$ is a root of $\phi(z)
= 0$. Then \eqref{ch02.2.3} implies $v_t(\lambda_0) = \lambda_0$ for all
$t\ge 0$. By Proposition~\ref{tch02.2.1} we have $v_t(\lambda)>\lambda_0$
for $\lambda> \lambda_0$ and $0<v_t(\lambda)<\lambda_0$ for $0<\lambda<
\lambda_0$. Then $\lambda>0$ and $\phi(\lambda)\neq 0$ imply there is no
root of $\phi(z) = 0$ between $\lambda$ and $v_t(\lambda)$. {From}
\eqref{ch02.2.2} we get \eqref{ch02.2.6}. \qed\end{proof}

 \bgproposition\label{tch02.2.4}
For any $t\ge 0$ and $\lambda\ge 0$ let $v_t^\prime(\lambda) =
(\partial/\partial \lambda) v_t(\lambda)$. Then we have
 \beqlb\label{ch02.2.7}
v_t^\prime(\lambda)
 =
\exp\Big\{- \int_0^t \phi^\prime(v_s(\lambda))\d s\Big\},
 \eeqlb
where
 \beqlb\label{ch02.2.8}
\phi^\prime(z)
 =
b + 2cz + \int_0^\infty u\big(1-\e^{-zu}\big)m(\d u).
 \eeqlb
\edproposition

 \begin{proof}
Based on \eqref{ch02.1.15} and \eqref{ch02.2.2} it is elementary to see
that
 \beqnn
\frac{\partial}{\partial t} v_t^\prime(\lambda)
 =
- \phi^\prime(v_t(\lambda))v_t^\prime(\lambda)
 =
\frac{\partial}{\partial \lambda}\frac{\partial}{\partial t}
v_t(\lambda).
 \eeqnn
It follows that
 \beqnn
\frac{\partial}{\partial t}\big[\log v_t^\prime(\lambda)\big]
 =
v_t^\prime(\lambda)^{-1} \frac{\partial}{\partial t} v_t^\prime(\lambda)
 =
- \phi^\prime(v_t(\lambda)).
 \eeqnn
Then we have \eqref{ch02.2.7} since $v_0^\prime(\lambda) = 1$.
\qed\end{proof}

Since $(Q_t)_{t\ge 0}$ is a Feller semigroup by
Corollary~\ref{tch02.2.2}, the CB-process has a Hunt realization $X =
(\Omega, \mcr{F}, \mcr{F}_t, x(t), \mbf{Q}_x)$. Let $\tau_0 := \inf\{s\ge
0: x(s)=0\}$ denote the \index{extinction time} \textit{extinction time}
of the CB-process.

 \bgtheorem\label{tch02.2.5}
For every $t\ge 0$ the limit $\bar{v}_t = \ulim_{\lambda\to \infty}
v_t(\lambda)$ exists in $(0,\infty]$. Moreover, the mapping $t\mapsto
\bar{v}_t$ is decreasing and for any $t\ge 0$ and $x>0$ we have
 \beqlb\label{ch02.2.9}
\mbf{Q}_x\{\tau_0\le t\} = \mbf{Q}_x\{x(t)=0\} = \exp\{-x\bar{v}_t\}.
 \eeqlb
\edtheorem

 \begin{proof}
By Proposition~\ref{tch02.2.1} the limit $\bar{v}_t = \ulim_{\lambda\to
\infty} v_t(\lambda)$ exists in $(0,\infty]$ for every $t\ge 0$. For
$t\ge r\ge 0$ we have
 \beqlb\label{ch02.2.10}
\bar{v}_t
 =
\ulim_{\lambda\to \infty} v_r(v_{t-r}(\lambda))
 =
v_r(\bar{v}_{t-r})
 \le
\bar{v}_r.
 \eeqlb
Since zero is a trap for the CB-process, we get \eqref{ch02.2.9} by
letting $\lambda\to \infty$ in \eqref{ch02.1.20}. \qed\end{proof}

For the convenience of statement of the results in the sequel, we
formulate the following condition on the branching mechanism:

 \bgcondition\label{tch02.2.6}
There is some constant $\theta>0$ so that
 \beqnn
\phi(z)>0 ~\mbox{for}~ z\ge \theta ~\mbox{and}~ \int_\theta^\infty
\phi(z)^{-1}\d z< \infty.
 \eeqnn
\edcondition

 \bgtheorem\label{tch02.2.7}
We have $\bar{v}_t< \infty$ for some and hence all $t>0$ if and only if
Condition~\ref{tch02.2.6} holds. \edtheorem

 \begin{proof}
By \eqref{ch02.2.10} it is simple to see that $\bar{v}_t =
\ulim_{\lambda\to \infty} v_t(\lambda)< \infty$ for all $t>0$ if and only
if this holds for some $t>0$. If Condition~\ref{tch02.2.6} holds, we can
let $\lambda\to \infty$ in \eqref{ch02.2.6} to obtain
 \beqlb\label{ch02.2.11}
\int_{\bar{v}_t}^\infty \phi(z)^{-1}\d z =t
 \eeqlb
and hence $\bar{v}_t< \infty$ for $t>0$. For the converse, suppose that
$\bar{v}_t<\infty$ for some $t>0$. By \eqref{ch02.2.2} there exists some
$\theta>0$ so that $\phi(\theta)>0$, for otherwise we would have
$v_t(\lambda)\ge \lambda$, yielding a contradiction. Then $\phi(z)>0$ for
all $z\ge \theta$ by the convexity of the branching mechanism. As in the
above we see that \eqref{ch02.2.11} still holds, so
Condition~\ref{tch02.2.6} is satisfied. \qed\end{proof}

 \bgtheorem\label{tch02.2.8}
Let $\bar{v} = \dlim_{t\to \infty} \bar{v}_t\in [0,\infty]$. Then for any
$x>0$ we have
 \beqlb\label{ch02.2.12}
\mbf{Q}_x\{\tau_0< \infty\} = \exp\{-x\bar{v}\}.
 \eeqlb
Moreover, we have $\bar{v}< \infty$ if and only if
Condition~\ref{tch02.2.6} holds, and in this case $\bar{v}$ is the
largest root of $\phi(z)=0$. \edtheorem

 \begin{proof}
The first assertion follows immediately from Theorem~\ref{tch02.2.5}. By
Theorem~\ref{tch02.2.7} we have $\bar{v}_t< \infty$ for some and hence
all $t>0$ if and only if Condition~\ref{tch02.2.6} holds. This is clearly
equivalent to $\bar{v}< \infty$. {From} \eqref{ch02.2.11} it is easy to
see that $\bar{v}$ is the largest root of $\phi(z)=0$. \qed\end{proof}

 \bgcorollary\label{tch02.2.9}
Suppose that Condition~\ref{tch02.2.6} holds. Then for any $x>0$ we have
$\mbf{Q}_x\{\tau_0< \infty\} = 1$ if and only if $b\ge 0$. \edcorollary

Let $(Q_t^\circ)_{t\ge 0}$ be the restriction to $(0,\infty)$ of the
semigroup $(Q_t)_{t\ge 0}$. A family of $\sigma$-finite measures
$(H_t)_{t>0}$ on $(0,\infty)$ is called an \index{entrance law}
\textit{entrance law} for $(Q_t^\circ)_{t\ge 0}$ if $H_rQ_t^\circ =
H_{r+t}$ for all $r,t>0$. The special case of the canonical
representation \eqref{ch02.2.1} with $h_t=0$ for all $t>0$ is
particularly interesting. In this case, we have
 \beqlb\label{ch02.2.13}
v_t(\lambda) = \int_0^\infty (1-\e^{-\lambda u}) l_t(\d u), \qquad
t>0,\lambda\ge 0.
 \eeqlb

 \bgtheorem\label{tch02.2.10}
The cumulant semigroup admits representation \eqref{ch02.2.13} if and
only if
 \beqlb\label{ch02.2.14}
\phi^\prime(\infty)
 :=
b + 2c\cdot\infty + \int_0^\infty u\, m(\d u)
 =
\infty
 \eeqlb
with $0\cdot\infty = 0$ by convention. If condition \eqref{ch02.2.14} is
satisfied, then $(l_t)_{t>0}$ is an entrance law for the restricted
semigroup $(Q_t^\circ)_{t\ge 0}$. \edtheorem

 \begin{proof}
{From} \eqref{ch02.2.8} it is clear that the limit $\phi^\prime(\infty) =
\lim_{z\to \infty} \phi^\prime(z)$ always exists in $(-\infty,\infty]$.
By \eqref{ch02.2.1} we have
 \beqlb\label{ch02.2.15}
v_t^\prime(\lambda) = h_t + \int_0^\infty u\e^{-\lambda u} l_t(\d u),
\qquad t\ge 0,\lambda\ge 0.
 \eeqlb
{From} \eqref{ch02.2.7} and \eqref{ch02.2.15} it follows that
 \beqlb\label{ch02.2.16}
h_t = v_t^\prime(\infty) = \exp\Big\{- \int_0^t \phi^\prime(\bar{v}_s)\d
s\Big\}.
 \eeqlb
Then $h_t = 0$ for any $t>0$ implies $\phi^\prime(\infty) = \infty$. For
the converse, assume that $\phi^\prime(\infty) = \infty$. If
Condition~\ref{tch02.2.6} holds, by Theorem~\ref{tch02.2.7} for every
$t>0$ we have $\bar{v}_t<\infty$, so $h_t = 0$ by \eqref{ch02.2.1}. If
Condition~\ref{tch02.2.6} does not hold, then $\bar{v}_t = \infty$ for
$t>0$ by Theorem~\ref{tch02.2.7}. Then \eqref{ch02.2.16} implies $h_t =
0$ for $t>0$. That proves the first assertion of the theorem. If
$(v_t)_{t>0}$ admits the representation \eqref{ch02.2.13}, we can use the
semigroup property of $(v_t)_{t\ge 0}$ to see
 \beqnn
\int_0^\infty (1-\e^{-\lambda u}) l_{r+t}(\d u)
 \ar=\ar
\int_0^\infty (1-\e^{- uv_t(\lambda)}) l_r(\d u) \cr
 \ar=\ar
\int_0^\infty l_r(\d x)\int_0^\infty (1-\e^{-\lambda u}) Q_t^\circ(x,\d
u)
 \eeqnn
for $r,t>0$ and $\lambda\ge 0$. Then $(l_t)_{t>0}$ is an entrance law for
$(Q_t^\circ)_{t\ge 0}$. \qed\end{proof}

 \bgcorollary\label{tch02.2.11}
If Condition~\ref{tch02.2.6} holds, the cumulant semigroup admits the
representation \eqref{ch02.2.13} and $t\mapsto \bar{v}_t = l_t(0,\infty)$
is the minimal solution of the differential equation
 \beqlb\label{ch02.2.17}
\frac{\d}{\d t}\bar{v}_t = -\phi(\bar{v}_t), \qquad t>0
 \eeqlb
with singular initial condition $\bar{v}_{0+} = \infty$. \edcorollary

 \begin{proof}
Under Condition~\ref{tch02.2.6}, for every $t>0$ we have $\bar{v}_t<
\infty$ by Theorem~\ref{tch02.2.7}. Moreover, the condition and the
convexity of $z\mapsto \phi(z)$ imply $\phi^\prime(\infty) = \infty$.
Then we have the representation \eqref{ch02.2.13} by
Theorem~\ref{tch02.2.10}. The semigroup property of $(v_t)_{t\ge 0}$
implies $\bar{v}_{t+s} = v_t(\bar{v}_s)$ for $t>0$ and $s>0$. Then
$t\mapsto \bar{v}_t$ satisfies \eqref{ch02.2.17}. {From}
\eqref{ch02.2.11} it is easy to see $\bar{v}_{0+} = \infty$. Using the
relation $\bar{v}_t = \lim_{\lambda\to \infty} v_t(\lambda)$ it is easy
to show that any solution $t\mapsto u_t$ of \eqref{ch02.2.17} with
$u_{0+} = \infty$ satisfies $u_t\ge \bar{v}_t$ for $t>0$. \qed\end{proof}

 \bgcorollary\label{tch02.2.12}
Suppose that Condition~\ref{tch02.2.6} holds. Then for any $t>0$ the
function $\lambda\mapsto v_t(\lambda)$ is strictly increasing and concave
on $[0,\infty)$, and $\bar{v}$ is the largest solution of the equation
$v_t(\lambda) = \lambda$. Moreover, we have $\bar{v} = \ulim_{t\to\infty}
v_t(\lambda)$ for $0<\lambda<\bar{v}$ and $\bar{v} = \dlim_{t\to\infty}
v_t(\lambda)$ for $\lambda>\bar{v}$. \edcorollary

 \begin{proof}
By Corollary~\ref{tch02.2.11} we have the canonical representation
\eqref{ch02.2.13} for every $t>0$. Since $\lambda\mapsto v_t(\lambda)$ is
strictly increasing by Proposition~\ref{tch02.2.1}, the measure $l_t(\d
u)$ is non-trivial, so $\lambda\mapsto v_t(\lambda)$ is strictly concave.
The equality $\bar{v} = v_t(\bar{v})$ follows by letting $s\to \infty$ in
$\bar{v}_{t+s} = v_t(\bar{v}_s)$, where $\bar{v}_{t+s}\le \bar{v}_s$.
Then $\bar{v}$ is clearly the largest solution to $v_t(\lambda) =
\lambda$. When $b\ge 0$, we have $\bar{v}=0$ by Theorem~\ref{tch02.2.8}
and Corollary~\ref{tch02.2.9}. Furthermore, since $\phi(z)\ge 0$, from
\eqref{ch02.2.2} we see $t\mapsto v_t(\lambda)$ is decreasing, and hence
$\dlim_{t\to\infty} v_t(\lambda) = \dlim_{t\to\infty} \bar{v}_t = 0$. If
$b<0$ and $0<\lambda<\bar{v}$, we have $\lambda\le v_t(\lambda)<
v_t(\bar{v}) = \bar{v}$ for all $t\ge 0$. Then the limit
$v_\infty(\lambda) = \ulim_{t\uparrow\infty} v_t(\lambda)$ exists. {From}
the relation $v_t(v_s(\lambda)) = v_{t+s}(\lambda)$ we have
$v_t(v_\infty(\lambda)) = v_\infty(\lambda)$, and hence
$v_\infty(\lambda) =\bar{v}$ since $\bar{v}$ is the unique solution to
$v_t(\lambda)=\lambda$ in $(0,\infty)$. The assertion for $b<0$ and
$\lambda>\bar{v}$ can be proved similarly. \qed\end{proof}

We remark that in Theorem~\ref{tch02.2.10} one usually cannot extend
$(l_t)_{t>0}$ to a $\sigma$-finite entrance law for the semigroup
$(Q_t)_{t\ge 0}$ on $\mbb{R}_+$. For example, let us assume
Condition~\ref{tch02.2.6} holds and $(\bar{l}_t)_{t>0}$ is such an
extension. For any $0< r< \varepsilon< t$ we have
 \beqnn
\bar{l}_t(\{0\})
 \ar\ge\ar
\int_0^\infty Q_{t-r}^\circ(x,\{0\}) l_r(\d x) \ge \int_0^\infty
\e^{-x\bar{v}_{t-\varepsilon}} l_r(\d x) \cr
 \ar=\ar
\bar{v}_r - \int_0^\infty (1-\e^{-u\bar{v}_{t-\varepsilon}}) l_r(\d u)
 =
\bar{v}_r - v_r(\bar{v}_{t-\varepsilon}).
 \eeqnn
The right-hand side tends to infinity as $r\to 0$. Then $\bar{l}_t(\d x)$
cannot be a $\sigma$-finite measure on $\mbb{R}_+$.

 \bexample\label{ech02.2.1}
Suppose that there are constants $c>0$, $0<\alpha\le 1$ and $b$ so that
$\phi(z) = cz^{1+\alpha} + bz$. Then Condition~\ref{tch02.2.6} is
satisfied. Let $q^0_\alpha(t) = \alpha t$ and
 \beqnn
q^b_\alpha(t)
 =
b^{-1}(1-\e^{-\alpha bt}), \qquad b\neq 0.
 \eeqnn
By solving the equation
 \beqnn
\frac{\partial}{\partial t}v_t(\lambda)
 =
- cv_t(\lambda)^{1+\alpha} - bv_t(\lambda),
 \qquad
v_0(\lambda) = \lambda
 \eeqnn
we get
 \beqlb\label{ch02.2.18}
v_t(\lambda)
 =
\frac{\e^{-bt}\lambda}{\big[1 + cq^b_\alpha(t)
\lambda^\alpha\big]^{1/\alpha}}, \qquad t\ge 0, \lambda\ge 0.
 \eeqlb
Thus $\bar{v}_t = c^{-1/\alpha}\e^{-bt}q^b_\alpha(t)^{-1/\alpha}$ for
$t>0$. In particular, if $\alpha = 1$, then \eqref{ch02.2.13} holds with
 \beqnn
l_t(\d u)
 =
\frac{\e^{-bt}}{c^2q^b_1(t)^2}\exp\Big\{-\frac{u}{cq^b_1(t)}\Big\}\d u,
\qquad t>0, u>0.
 \eeqnn
\eexample

\section{Conditional limit theorems}\label{sch02.3}

 \setcounter{equation}{0}

Let $(Q_t)_{t\ge 0}$ denote the transition semigroup of the CB-process
with branching mechanism $\phi$ given by \eqref{ch02.1.13}. Let
$(Q_t^\circ)_{t\ge 0}$ be the restriction of $(Q_t)_{t\ge 0}$ to
$(0,\infty)$. It is easy to check that $Q^b_t(x,\d y) := \e^{bt}x^{-1}y
Q_t^\circ(x,\d y)$ defines a Markov semigroup on $(0,\infty)$. Let
$q_t(\lambda) = \e^{bt} v_t(\lambda)$ and let $q_t^\prime(\lambda) =
(\partial/\partial \lambda) q_t(\lambda)$. Recall that $z\mapsto
\phi^\prime(z)$ is defined by \eqref{ch02.2.8}. {From} \eqref{ch02.2.7}
we have
 \beqlb\label{ch02.3.1}
q_t^\prime(\lambda)
 =
\exp\Big\{- \int_0^t\phi_0^\prime(v_s(\lambda))\d s\Big\},
 \eeqlb
where $\phi_0^\prime(z) = \phi^\prime(z) - b$. By differentiating both
sides of \eqref{ch02.1.20} we see
 \beqlb\label{ch02.3.2}
\int_0^\infty\e^{-\lambda y}Q^b_t(x,\d y)
 =
\exp\{-xv_t(\lambda)\}q_t^\prime(\lambda), \quad \lambda\ge 0.
 \eeqlb
It follows that
 \beqlb\label{ch02.3.3}
\int_0^\infty\e^{-\lambda y}Q^b_t(x,\d y)
 =
\exp\Big\{-xv_t(\lambda) - \int_0^t\phi^\prime_0(v_s(\lambda))\d s\Big\}.
 \eeqlb
Using \eqref{ch02.3.3} it is easy to extend $(Q^b_t)_{t\ge 0}$ to a
Feller semigroup on $\mbb{R}_+$.

 \bgtheorem\label{tch02.3.1}
Suppose that $\phi^\prime(z)\to \infty$ as $z\to \infty$. Let
$(l_t)_{t>0}$ be defined by \eqref{ch02.2.13}. Then for $t>0$ the
probability measure $Q^b_t(0,\cdot)$ is supported by $(0,\infty)$ and
$Q^b_t(0,\d u) = u\e^{bt} l_t(\d u)$. \edtheorem

 \begin{proof}
We first note that $(v_t)_{t>0}$ really has the representation
\eqref{ch02.2.13} by Theorem~\ref{tch02.2.10}. Then
 \beqlb\label{ch02.3.4}
q_t(\lambda) = \int_0^\infty \big(1-\e^{-\lambda u}\big)\e^{bt} l_t(\d
u),
 \eeqlb
and
 \beqlb\label{ch02.3.5}
q_t^\prime(\lambda) = \int_0^\infty u\e^{-\lambda u}\e^{bt}l_t(\d u).
 \eeqlb
{From} \eqref{ch02.3.2} and \eqref{ch02.3.5} we see $Q^b_t(0,\d u) =
u\e^{bt} l_t(\d u)$ for $u>0$. \qed\end{proof}

Now let $X = (\Omega, \mcr{F}, \mcr{F}_t, x(t), \mbf{Q}_x)$ be a Hunt
realization of the CB-process with the augmented natural
$\sigma$-algebras. Let $\tau_0 := \inf\{s\ge 0: x(s)=0\}$ denote the
extinction time of $X$.

 \bgtheorem\label{tch02.3.2}
Suppose that $b\ge 0$ and Condition~\ref{tch02.2.6} holds. Let $t\ge 0$
and $x>0$. Then the distribution of $x(t)$ under $\mbf{Q}_x\{\cdot|r+t<
\tau_0\}$ converges as $r\to \infty$ to $Q^b_t(x,\cdot)$. \edtheorem

 \begin{proof}
Since zero is a trap for the CB-process, for any $r>0$ we can use the
Markov property of $\{x(t): t\ge 0\}$ to see
 \beqlb\label{ch02.3.6}
\mbf{Q}_x\big[\e^{-\lambda x(t)}|r+t<\tau_0\big]
 \ar=\ar
\frac{\mbf{Q}_x\big[\e^{-\lambda x(t)}1_{\{r+t<\tau_0\}}\big]}
{\mbf{Q}_x[1_{\{r+t<\tau_0\}}]} \cr
 \ar=\ar
\lim_{\theta\to \infty}\frac{\mbf{Q}_x\big[\e^{-\lambda x(t)}(1-
\e^{-\theta x(r+t)})\big]} {\mbf{Q}_x\big[(1- \e^{-\theta x(r+t)})\big]}
\cr
 \ar=\ar
\frac{\mbf{Q}_x\big[\e^{-\lambda x(t)}(1- \e^{-x(t)\bar{v}_r})\big]}
{1-\e^{-x\bar{v}_{r+t}}}.
 \eeqlb
Recall that $\bar{v}_{r+t} = v_t(\bar{v}_r)$ and $v_t^\prime(0)
=\e^{-bt}$. By Theorem~\ref{tch02.2.8} and Corollary~\ref{tch02.2.9} we
have $\lim_{r\to \infty} \bar{v}_r = 0$. Then
 \beqnn
\lim_{r\to \infty}\mbf{Q}_x\big[\e^{-\lambda x(t)}| r+t<\tau_0\big]
 \ar=\ar
\lim_{r\to \infty} \frac{\mbf{Q}_x\big[\e^{-\lambda x(t)}
\bar{v}_r^{-1}(1- \e^{-x(t)\bar{v}_r})\big]}
{\bar{v}_r^{-1}(1-\e^{-xv_t(\bar{v}_r)})} \cr
 \ar=\ar
\frac{1}{x}\e^{bt}\mbf{Q}_x[x(t)\e^{-\lambda x(t)}].
 \eeqnn
That gives the desired convergence result. \qed\end{proof}

It is easy to see that $t\mapsto Z(t) := \e^{bt}x(t)$ is a positive
$(\mcr{F}_t)$-martingale. By the theory of Markov processes, for each
$x>0$ there is a unique probability measure $\mbf{Q}^b_x$ on $(\Omega,
\mcr{F})$ so that
 \beqlb\label{ch02.3.7}
\mbf{Q}^b_x(F) = \mbf{Q}_x[Z(t)F]
 \eeqlb
for any $\mcr{F}_t$-measurable bounded random variable $F$. Moreover,
under this new probability measure $\{x(t): t\ge 0\}$ is a Markov process
in $(0,\infty)$ with transition semigroup $(Q^b_t)_{t\ge 0}$. By a
modification of the proof of Theorem~\ref{tch02.3.2} we get the
following:

 \bgtheorem\label{tch02.3.3}
Suppose that $b\ge 0$ and Condition~\ref{tch02.2.6} holds. Let $x>0$ and
$t\ge 0$. Then for any $\mcr{F}_t$-measurable bounded random variable $F$
we have
 \beqlb\label{ch02.3.8}
\mbf{Q}^b_x[F] = \lim_{r\to \infty}\mbf{Q}_x[F|r+t<\tau_0].
 \eeqlb
\edtheorem

By the above theorem, in the critical and subcritical cases,
$\mbf{Q}^b_x$ is intuitively the law of $\{x(t): t\ge 0\}$ conditioned on
large extinction times. See Lambert (2007), Li (2000, 2011) and Pakes
(1999) for more conditional limit theorems.

\section{A reconstruction from excursions}\label{sch02.4}

 \setcounter{equation}{0}

In this section, we give a reconstruction of the sample paths of the
CB-process from excursions. Let $(Q_t)_{t\ge 0}$ denote the transition
semigroup of the process with branching mechanism $\phi$ given by
\eqref{ch02.1.13}. Recall that $(Q_t^\circ)_{t\ge 0}$ is the restriction
of $(Q_t)_{t\ge 0}$ to $(0,\infty)$. Let $D(0,\infty)$ be the space of
c\`{a}dl\`{a}g paths $t\mapsto w_t$ from $(0,\infty)$ to $\mbb{R}_+$
having zero as a trap. Let $(\mcr{A}^0, \mcr{A}^0_t)$ be the natural
$\sigma$-algebras on $D(0,\infty)$ generated by the coordinate process.
For any entrance law $(H_t)_{t>0}$ for $(Q_t^\circ)_{t\ge 0}$ there is a
unique $\sigma$-finite measure $\mbf{Q}_H(\d w)$ on $\mcr{A}^0$ such that
$\mbf{Q}_H(\{0\}) = 0$ and
 \beqlb\label{ch02.4.1}
\ar\ar\mbf{Q}_H(w_{t_1}\in \d x_1, w_{t_2}\in \d x_2, \ldots, w_{t_n}\in
\d x_n) \ccr
 \ar\ar\qquad
=\, H_{t_1}(\d x_1) Q_{t_2-t_1}^\circ(x_1,\d x_2) \cdots
Q_{t_n-t_{n-1}}^\circ (x_{n-1},\d x_n)
 \eeqlb
for every $\{t_1< \cdots< t_n\}\subset (0,\infty)$ and $\{x_1, \ldots,
x_n\}\subset (0,\infty)$. See, e.g., Getoor and Glover (1987) for the
proof of the existence of $\mbf{Q}_H$ in the setting of Borel right
processes. Roughly speaking, the above formula means that $\{w_t: t>0\}$
under $\mbf{Q}_H$ is a Markov process in $(0,\infty)$ with transition
semigroup $(Q_t^\circ)_{t\ge 0}$ and one-dimensional distributions
$(H_t)_{t>0}$.

 \bgtheorem\label{tch02.4.1}
Suppose that the condition \eqref{ch02.2.14} is satisfied. Let
$\mbf{Q}_{(0)}$ be the $\sigma$-finite measure on $D(0,\infty)$
determined by the entrance law $(l_t)_{t>0}$ given by \eqref{ch02.2.13}.
Then for $\mbf{Q}_{(0)}$-a.e.\ $w\in D(0,\infty)$ we have $w_t\to 0$ as
$t\to 0$. \edtheorem

 \begin{proof}
Let $(Q_t^b)_{t\ge 0}$ be the transition semigroup on $\mbb{R}_+$ defined
by \eqref{ch02.3.3}. Then we have $Q^b_t(x,\d y) = x^{-1}\e^{bt}y
Q_t^\circ(x,\d y)$ for $x,y>0$. Observe that
 \beqnn
\int_{D(0,\infty)}\e^{bt}w_t\mbf{Q}_{(0)}(\d w)
 =
\int_0^\infty \e^{bt}y l_t(\d y)
 =
\e^{bt}\frac{\partial}{\partial\lambda}v_t(0+) = 1.
 \eeqnn
Then for fixed $u>0$ we can define a probability measure
$\mbf{Q}_{(0)}^u(\d w) := \e^{bu}w_u\mbf{Q}_{(0)}(\d w)$ on
$D(0,\infty)$. Under this measure, the coordinate process $\{w_t: 0<t\le
u\}$ is an immigration process with transition semigroup $(Q_t^b)_{t\ge
0}$ and one-dimensional distributions
 \beqnn
Q^b_t(0,\d y) = \e^{bt}yl_t(\d y), \qquad 0<t\le u.
 \eeqnn
By the uniqueness of the transition law of the immigration process we
have $w_t\to 0$ as $t\to 0$ for $\mbf{Q}_{(0)}^u$-a.e.\ $w\in
D(0,\infty)$. Note that $\mbf{Q}_{(0)}^u(\d w)$ and $\mbf{Q}_{(0)}(\d w)$
are absolutely continuous with respect to each other on $D_u(0,\infty) :=
\{w\in D(0,\infty): w_u> 0\}$. Since $D(0,\infty) = \{0\}\cup
(\cup_{n=1}^\infty D_{1/n}(0,\infty))$ and $\mbf{Q}_{(0)}(\{0\}) = 0$, we
have $w_t\to 0$ as $t\to 0$ for $\mbf{Q}_{(0)}$-a.e.\ $w\in D(0,\infty)$.
\qed\end{proof}

By Theorem~\ref{tch02.4.1} we can regard $\mbf{Q}_{(0)}$ as a
$\sigma$-finite measure on the space of \index{excursion}
\textit{excursions} $D_0[0,\infty) := \{w\in D(0,\infty): w_0 = 0\}$. We
call $\mbf{Q}_{(0)}$ an \index{excursion law} \textit{excursion law} for
the CB-process. Using the excursions, we can give a reconstruction of the
sample paths of the CB-process. Let $x\ge 0$ and let $N(\d w)$ be a
Poisson random measure on $D_0[0,\infty)$ with intensity
$x\mbf{Q}_{(0)}(\d w)$. We define the process $\{X_t: t\ge 0\}$ by
$X_0=x$ and
 \beqlb\label{ch02.4.2}
X_t = \int_{D_0[0,\infty)} w_t N(\d w), \qquad t>0.
 \eeqlb

 \bgtheorem\label{tch02.4.2}
For $t\ge 0$ let $\mcr{G}_t$ be the $\sigma$-algebra generated by the
collection of random variables $\{N(A): A\in \mcr{A}^0_t\}$. Then
$\{(X_t,\mcr{G}_t): t\ge 0\}$ is a realization of the CB-process.
\edtheorem

 \begin{proof}
We first remark that the random variable $X_t$ has distribution
$Q_t(x,\cdot)$ on $\mbb{R}_+$. In fact, for any $t>0$ and $\lambda\ge 0$
we have
 \beqnn
\mbf{P}\big[\exp\{-\lambda X_t\}\big]
 \ar=\ar
\exp\Big\{-x\int_{D_0[0,\infty)} (1-\e^{-\lambda w_t}) \mbf{Q}_{(0)}(\d
w)\Big\} \cr
 \ar=\ar
\exp\Big\{-x\int_0^\infty (1 - \e^{-\lambda z})l_t(\d z)\Big\}
 =
\exp\{-xv_t(\lambda)\}.
 \eeqnn
Let $t>r>0$ and let $h$ be a bounded positive function on $D_0[0,\infty)$
measurable relative to $\mcr{A}^0_r$. For any $\lambda\ge 0$ we have
 \beqnn
 \ar\ar\mbf{P}\Big[\exp\Big\{-\int_{D_0[0,\infty)}h(w)N(\d w)
-\lambda X_t\Big\}\Big] \cr
 \ar\ar\qquad
=\, \exp\Big\{-x\int_{D_0[0,\infty)}\big(1 - \e^{-h(w)-\lambda w_t}\big)
\mbf{Q}_{(0)}(\d w)\Big\} \cr
 \ar\ar\qquad
=\, \exp\Big\{-x\int_{D_0[0,\infty)} \big(1 - \e^{-h(w)}\big)
\mbf{Q}_{(0)}(\d w)\Big\} \cr
 \ar\ar\qquad\quad
\cdot\exp\Big\{-x\int_{D_0[0,\infty)} \e^{-h(w)}\big(1 - \e^{- \lambda
w_t}\big)\mbf{Q}_{(0)}(\d w)\Big\},
 \eeqnn
where we made the convention $\e^{-\infty}=0$. By the Markov property of
$\mbf{Q}_{(0)}$ we have
 \beqnn
 \ar\ar\int_{D_0[0,\infty)} \e^{-h(w)}\big(1 - \e^{-\lambda w_t}\big)
\mbf{Q}_{(0)}(\d w) \cr
 \ar\ar\qquad
= \int_{D_0[0,\infty)}\e^{-h(w)}\mbf{Q}_{(0)}(\d w)\int_0^\infty (1 -
\e^{-\lambda z}) Q_{t-r}^\circ(w_r,\d z) \cr
 \ar\ar\qquad
= \int_{D_0[0,\infty)}\e^{-h(w)}\big(1 - \e^{-w_rv_{t-r}(\lambda)}\big)
\mbf{Q}_{(0)}(\d w).
 \eeqnn
It follows that
 \beqnn
 \ar\ar\mbf{P}\Big[\exp\Big\{-\int_{D_0[0,\infty)}h(w)N(\d w)
- \lambda X_t\Big\}\Big] \cr
 \ar\ar\qquad
=\, \exp\Big\{-\int_{D_0[0,\infty)}\big(1 - \e^{-h(w)}
\e^{-w_rv_{t-r}(\lambda)}\big) \mbf{Q}_{(0)}(\d w)\Big\} \cr
 \ar\ar\qquad
=\, \mbf{P}\Big[\exp\Big\{-\int_{D_0[0,\infty)}h(w)N(\d w) -
X_rv_{t-r}(\lambda)\Big\}\Big].
 \eeqnn
Then $\{(X_t,\mcr{G}_t): t\ge 0\}$ is a Markov process with transition
semigroup $(Q_t)_{t\ge 0}$. \qed\end{proof}

\chapter{Structures of independent immigration }\label{ch03}

In this chapter we study independent immigration structures associated
with CB-processes. We first give a formulation of the structures using
skew convolution semigroups. Those semigroups are in one-to-one
correspondence with infinitely divisible distributions on $\mbb{R}_+$. We
show the corresponding immigration process arise as scaling limits of
Galton--Watson processes with immigration. We discuss briefly limit
theorems and stationary distributions of the immigration superprocesses.
The trajectories of the immigration processes are constructed using
stochastic integrals with respect to Poisson random measures determined
by entrance laws.

\section{Formulation of immigration processes}\label{sch03.1}

 \setcounter{equation}{0}

In this section, we introduce a generalization of the CB-process. Let
$(Q_t)_{t\ge 0}$ be the transition semigroup defined by \eqref{ch02.1.15}
and \eqref{ch02.1.20}. Let $(\gamma_t)_{t\ge 0}$ be a family of
probability measures on $\mbb{R}_+$. We call $(\gamma_t)_{t\ge 0}$ a
\index{skew convolution semigroup} \textit{skew convolution semigroup}
\index{SC-semigroup} (SC-semigroup) associated with $(Q_t)_{t\ge 0}$
provided
 \beqlb\label{ch03.1.1}
\gamma_{r+t} = (\gamma_rQ_t)*\gamma_t, \qquad r,t\ge 0.
 \eeqlb
It is easy to show that \eqref{ch03.1.1} holds if and only if
 \beqlb\label{ch03.1.2}
u_{r+t}(\lambda) = u_t(\lambda) + u_r(v_t(\lambda)), \qquad
r,t,\lambda\ge 0,
 \eeqlb
where
 \beqlb\label{ch03.1.3}
u_t(\lambda) = -\log\int_0^\infty \e^{-y\lambda} \gamma_t(\d y).
 \eeqlb
The concept of SC-semigroup is of interest because of the following:

 \bgtheorem\label{tch03.1.1}
The family of probability measures $(\gamma_t)_{t\ge 0}$ on $\mbb{R}_+$
is an SC-semigroup if and only if
 \beqlb\label{ch03.1.4}
Q^\gamma_t(x,\cdot) := Q_t(x,\cdot)*\gamma_t, \qquad t,x\ge 0
 \eeqlb
defines a Markov semigroup $(Q^\gamma_t)_{t\ge 0}$ on $\mbb{R}_+$.
\edtheorem

 \begin{proof}
Let $(\gamma_t)_{t\ge 0}$ probability measures on $\mbb{R}_+$ and let
$Q^\gamma_t(x,\cdot)$ be the probability kernel defined by
\eqref{ch03.1.4}. Then we have
 \beqlb\label{ch03.1.5}
\int_0^\infty \e^{-y\lambda} Q^\gamma_t(x,\d y)
 =
\exp\big\{-xv_t(\lambda)-u_t(\lambda)\big\}, \qquad t,x,\lambda\ge 0.
 \eeqlb
Using this relation it is easy to show that $(Q^\gamma_t)_{t\ge 0}$
satisfies the Chapman-Kolmogorov equation if and only if \eqref{ch03.1.2}
is satisfied. That proves the result. \qed\end{proof}

If $\{y(t): t\ge 0\}$ is a positive Markov process with transition
semigroup $(Q^\gamma_t)_{t\ge 0}$ given by \eqref{ch03.1.4}, we call it
an \index{immigration process} \textit{immigration process} or a
\index{CBI-process} \textit{CBI-process} associated with $(Q_t)_{t\ge
0}$. The intuitive meaning of the model is clear in view of
\eqref{ch03.1.1} and \eqref{ch03.1.4}. {From} \eqref{ch03.1.4} we see
that the population at any time $t\ge 0$ is made up of two parts; the
native part generated by the mass $x\ge 0$ has distribution
$Q_t(x,\cdot)$ and the immigration in the time interval $(0,t]$ gives the
distribution $\gamma_t$. In a similar way, the equation \eqref{ch03.1.1}
decomposes the mass immigrating to the population during the time
interval $(0,r+t]$ into two parts; the immigration in the interval
$(r,r+t]$ gives the distribution $\gamma_t$ while the immigration in the
interval $(0,r]$ generates the distribution $\gamma_r$ at time $r$ and
gives the distribution $\gamma_rQ_t$ at time $r+t$. It is not hard to
understand that \eqref{ch03.1.4} gives a general formulation of the
immigration independent of the state of the population.

 \bgtheorem\label{tch03.1.2}
The family of probability measures $(\gamma_t)_{t\ge 0}$ on $\mbb{R}_+$
is an SC-semigroup if and only if there exists $\psi\in \mcr{I}$ so that
 \beqlb\label{ch03.1.6}
\int_0^\infty \e^{-\lambda y} \gamma_t(\d y)
 =
\exp\Big\{-\int_0^t\psi(v_s(\lambda))\d s\Big\}, \qquad t,\lambda\ge 0.
 \eeqlb
\edtheorem

\begin{proof} It is easy to check that for any $\psi\in \mcr{I}$ the family $(\gamma_t)_{t\ge 0}$ defined by \eqref{ch03.1.6} is an SC-semigroup. Conversely, suppose that $(\gamma_t)_{t\ge 0}$ is an SC-semigroup. For $t,\lambda\ge 0$ let $u_t(\lambda)$ be defined by \eqref{ch03.1.3}. Then $t\mapsto u_t(\lambda)$ is increasing. By Lebesgue's theorem, the limit
 \beqnn
\psi_t(\lambda) := \lim_{s\to 0^+}s^{-1}[u_{t+s}(\lambda)-u_t(\lambda)] =
\lim_{s\to 0^+}s^{-1}u_s(v_t(\lambda))
 \eeqnn
exists for almost all $t\ge 0$; see, e.g., Hewitt and Stromberg (1965,
p.264). By the continuity of $t\mapsto v_t(\lambda)$, there is a dense
subset $D$ of $(0,\infty)$ so that the following limits exist:
 \beqlb\label{ch03.1.7}
\psi_0(\lambda) := \lim_{s\to 0^+}s^{-1}u_s(\lambda)
 =
\lim_{s\to 0^+}s^{-1}[1-\e^{-u_s(\lambda)}], \qquad \lambda\in D.
 \eeqlb
For $u\in [0,\infty]$ and $\lambda\in (0,\infty)$ let $\xi(u,\lambda)$ be
defined by \eqref{ch01.2.4}. Then $\xi(u,\lambda)$ is jointly continuous
in $(u,\lambda)$. By \eqref{ch03.1.7} we have
 \beqlb\label{ch03.1.8}
\psi_0(\lambda)
 =
\lim_{s\to 0^+}s^{-1}u_s(\lambda)
 =
\lim_{s\to 0^+}s^{-1} \int_{[0,\infty]} \xi(u,\lambda) G_s(\d u)
 \eeqlb
for $\lambda\in D$, where $G_s(\d u) = (1-\e^{-u})\gamma_s(\d u)$. Then
for some $\delta>0$ we have
 \beqnn
\sup_{0<s<\delta} s^{-1}G_s([0,\infty])< \infty,
 \eeqnn
so the family of finite measures $\{s^{-1}G_s: 0<s<\delta\}$ on
$[0,\infty]$ is relatively compact. Suppose that $s_n\to 0$ and
$s_n^{-1}G_{s_n}\to$ some $G$ weakly as $n\to\infty$. Then
 \beqnn
\lim_{n\to \infty}s_n^{-1}u_{s_n}(\lambda)
 =
\lim_{n\to \infty}s_n^{-1} \int_{[0,\infty]} \xi(u,\lambda) G_{s_n}(\d u)
 =
\int_{[0,\infty]}\xi(u,\lambda) G(du), \quad \lambda>0.
 \eeqnn
Thus we can extend $\psi_0$ to a continuous function on $(0,\infty)$
given by
 \beqnn
\psi_0(\lambda) = \int_{[0,\infty]}\xi(u,\lambda) G(du)
 =
G(\{\infty\}) + h\lambda + \int_{(0,\infty)}(1-\e^{-u\lambda}) l(du),
 \eeqnn
where $h = G(\{0\})$ and $l(\d u) = (1-\e^{-u})^{-1}G_s(\d u)$. By a
standard argument one sees that \eqref{ch03.1.8} holds actually for all
$\lambda>0$. From \eqref{ch03.1.2} and \eqref{ch03.1.8} it follows that
 \beqnn
D^+u_s(\lambda)\big|_{s=t}
 =
D^+u_s(v_t(\lambda))\big|_{s=0}
 =
\psi_0(v_t(\lambda)),
 \eeqnn
where $D^+$ denotes the right derivative relative to $s\ge 0$. Here the
right-hand side is continuous in $t\ge 0$. Thus $t\mapsto u_t(\lambda)$
is continuously differentiable and \eqref{ch03.1.6} holds with
$\psi(\lambda) = \psi_0(\lambda)$ for $\lambda>0$. By letting $\lambda\to
0+$ in \eqref{ch03.1.6} one sees $G(\{\infty\})=0$. Then $\psi=\psi_0\in
\mcr{I}$. \qed\end{proof}

By Theorem~\ref{tch03.1.2} there is a 1-1 correspondence between
SC-semigroups and infinitely divisible distributions on $\mbb{R}_+$. Then
the theorem generalizes the 1-1 correspondence between classical
convolution semigroups and infinitely divisible distributions. In fact,
from \eqref{ch03.1.1} it is easy to see that $(\gamma_t)_{t\ge 0}$
reduces to a classical convolution semigroup if $Q_t$ is the identity
operator for all $t\ge 0$. As a consequence of Theorems~\ref{tch01.2.4}
and~\ref{tch03.1.2}, an SC-semigroup $(\gamma_t)_{t\ge 0}$ always
consists of infinitely divisible distributions.

Now let us consider a transition semigroup $(Q^\gamma_t)_{t\ge 0}$
defined by \eqref{ch03.1.4} with the SC-semigroup $(\gamma_t)_{t\ge 0}$
given by \eqref{ch03.1.6}. If an immigration process has transition
semigroup $(Q_t^\gamma)_{t\ge 0}$, we say it has \index{branching
mechanism} \textit{branching mechanism} $\phi$ and \index{immigration
mechanism} \textit{immigration mechanism} $\psi$. It is easy to see that
 \beqlb\label{ch03.1.9}
\int_0^\infty \e^{-\lambda y} Q_t^\gamma(x,\d y)
 =
\exp\Big\{-xv_t(\lambda)-\int_0^t\psi(v_s(\lambda))\d s\Big\}.
 \eeqlb
Recall that the function $\psi\in \mcr{I}$ has the representation
 \beqlb\label{ch03.1.10}
\psi(z) = \beta z + \int_0^\infty\big(1-\e^{-zu}\big)n(\d u), \qquad z\ge
0,
 \eeqlb
where $\beta\ge 0$ is a constant and $(1\land u)n(\d u)$ is a finite
measure on $(0,\infty)$. By \eqref{ch03.1.9} and \eqref{ch02.2.4} one can
show
 \beqlb\label{ch03.1.11}
\int_0^\infty y Q_t^\gamma(x,\d y)
 =
x\e^{-bt} - \psi'(0)\int_0^t\e^{-bs}\d s,
 \eeqlb
where
 \beqlb\label{ch03.1.12}
\psi'(0) = \beta + \int_0^\infty u n(\d u).
 \eeqlb

The following theorem gives a necessary and sufficient condition for the
ergodicity of the semigroup $(Q^\gamma_t)_{t\ge 0}$.

 \bgtheorem\label{tch03.1.3}
Suppose that $b\ge 0$ and $\phi(z)\neq 0$ for every $z>0$. Then
$Q_t^\gamma(x,\cdot)$ converges to a probability measure $\eta$ on
$\mbb{R}_+$ as $t\to \infty$ if and only if
 \beqlb\label{ch03.1.13}
\int_0^\lambda \frac{\psi(z)}{\phi(z)}\d z<\infty ~~ \mbox{for some
$\lambda>0$.}
 \eeqlb
If \eqref{ch03.1.13} holds, the Laplace transform of $\eta$ is given by
 \beqlb\label{ch03.1.14}
L_\eta(\lambda)
 =
\exp\Big\{- \int_0^\infty\psi(v_s(\lambda))\d s\Big\}, \qquad \lambda\ge
0.
 \eeqlb
\edtheorem

 \begin{proof}
Since $\phi(z)\ge 0$ for all $z\ge 0$, from \eqref{ch02.2.2} we see
$t\mapsto v_t(\lambda)$ is decreasing. Then \eqref{ch02.2.6} implies
$\lim_{t\to \infty} v_t(\lambda) = 0$. By \eqref{ch03.1.9} we have
 \beqlb\label{ch03.1.15}
\lim_{t\to \infty}\int_0^\infty \e^{-\lambda y} Q_t^\gamma(x,\d y)
 =
\exp\Big\{- \int_0^\infty\psi(v_s(\lambda))\d s\Big\}
 \eeqlb
for every $\lambda\ge 0$. A further application of \eqref{ch02.2.2} gives
 \beqnn
\int_0^t\psi(v_s(\lambda))\d s
 =
\int_{v_t(\lambda)}^\lambda \frac{\psi(z)}{\phi(z)}\d z.
 \eeqnn
It follows that
 \beqnn
\int_0^\infty\psi(v_s(\lambda))\d s
 =
\int_0^\lambda \frac{\psi(z)}{\phi(z)}\d z,
 \eeqnn
which is a continuous function of $\lambda\ge 0$ if and only if
\eqref{ch03.1.13} holds. Then the result follows by \eqref{ch03.1.15} and
Theorem~\ref{tch01.1.2}. \qed\end{proof}

 \bgcorollary\label{tch03.1.4}
Suppose that $b>0$. Then $Q_t^\gamma(x,\cdot)$ converges to a probability
measure $\eta$ on $\mbb{R}_+$ as $t\to \infty$ if and only if
$\int_1^\infty \log u n(\d u)< \infty$. In this case, the Laplace
transform of $\eta$ is given by \eqref{ch03.1.14}. \edcorollary

 \begin{proof}
We have $\phi(z) = bz + o(z)$ as $z\to 0$. Thus \eqref{ch03.1.13} holds
if and only if
 \beqnn
\int_0^\lambda \frac{\psi(z)}{z}\d z<\infty  ~~ \mbox{for some
$\lambda>0$,}
 \eeqnn
which is equivalent to
 \beqnn
\int_0^\lambda \frac{\d z}{z}\int_0^\infty \big(1-\e^{-zu}\big) n(\d u)
 =
\int_0^\infty n(\d u)\int_0^{\lambda u} \frac{1-\e^{-y}}{y}\d y<\infty
 \eeqnn
for some $\lambda>0$. The latter holds if and only if $\int_1^\infty \log
u n(\d u)< \infty$. Then we have the result by Theorem~\ref{tch03.1.3}.
\qed\end{proof}

In the situation of Theorem~\ref{tch03.1.3}, it is easy to show that
$\eta$ is a stationary distribution for $(Q_t^\gamma)_{t\ge 0}$. The fact
that the CBI-process may have a non-trivial stationary distribution makes
it a more interesting model in many respects than the CB-process without
immigration. Note also that the transition semigroup $(Q^b_t)_{t\ge 0}$
given by \eqref{ch02.3.3} is a special case of the one defined by
\eqref{ch03.1.9}.

 \bgtheorem\label{tch03.1.5}
Suppose that $b>0$ and let $q_t^\prime(\lambda)$ be defined by
\eqref{ch02.3.1}. Then for every $\lambda\ge 0$ the limit
$q^\prime(\lambda) := \dlim_{t\to \infty}$ $q_t^\prime(\lambda)$ exists
and is given by
 \beqlb\label{ch03.1.16}
q^\prime(\lambda)
 =
\exp\Big\{- \int_0^\infty\phi_0^\prime(v_s(\lambda))\d s\Big\}, \qquad
\lambda\ge 0.
 \eeqlb
Moreover, we have $q^\prime(0+) = q^\prime(0) = 1$ if and only if
$\int_1^\infty u\log u m(\d u)< \infty$. The last condition is also
equivalent to $q^\prime(\lambda)> 0$ for some and hence all $\lambda>0$.
\edtheorem

 \begin{proof}
The first assertion is easy in view of \eqref{ch02.3.1}. By
Corollary~\ref{tch03.1.4}, we have $\int_1^\infty u\log u m(\d u)<
\infty$ if and only if $\lambda\mapsto q^\prime(\lambda)$ is the Laplace
transform of a probability $\eta$ on $\mbb{R}_+$. Then the other two
assertions hold obviously. \qed\end{proof}

 \bgtheorem\label{tch03.1.6}
Suppose that $b>0$ and $\phi^\prime(z)\to \infty$ as $z\to \infty$. Then
$Q^b_t(x,\cdot)$ converges as $t\to \infty$ to a probability $\eta$ on
$(0,\infty)$ if and only if $\int_1^\infty u\log u m(\d u)< \infty$. If
the condition holds, then $\eta$ has Laplace transform $L_\eta =
q^\prime$ given by \eqref{ch03.1.16}. \edtheorem

 \begin{proof}
By Corollary~\ref{tch03.1.4} and Theorem~\ref{tch03.1.5} we have the
results with $\eta$ being a probability measure on $\mbb{R}_+$. By
Theorem~\ref{tch02.3.1} the measure $Q^b_t(0,\cdot)$ is supported by
$(0,\infty)$, hence $Q^b_t(x,\cdot)$ is supported by $(0,\infty)$ for
every $x\ge 0$. {From} \eqref{ch02.3.5} we have $L_\eta(\infty)\le
q_t^\prime(\infty) = 0$ for $t>0$. That implies $\eta(\{0\}) = 0$.
\qed\end{proof}

 \bexample\label{ech03.1.1}
Suppose that $c>0$, $0<\alpha\le 1$ and $b$ are constants and let
$\phi(z) = cz^{1+\alpha} + bz$ for $z\ge 0$. In this case the cumulant
semigroup $(v_t)_{t\ge 0}$ is given by \eqref{ch02.2.18}. Let $\beta\ge
0$ and let $\psi(z) = \beta z^\alpha$ for $z\ge 0$. We can use
\eqref{ch03.1.9} to define the transition semigroup $(Q^\gamma_t)_{t\ge
0}$. It is easy to show that
 \beqlb\label{ch03.1.17}
\int_0^\infty \e^{-\lambda y} Q_t^\gamma(x,\d y)
 =
\frac{1}{\big[1 + cq^b_\alpha(t) \lambda^\alpha\big]^{\beta/c\alpha}}\,
\e^{-xv_t(\lambda)}, \qquad \lambda\ge 0.
 \eeqlb
\eexample

The concept of SC-semigroup associated with branching processes was
introduced in Li (1995/6, 1996). Theorem~\ref{tch03.1.2} can be regarded
as a special form of main theorem of Li (1995/6). Theorem~\ref{tch03.1.3}
and Corollary~\ref{tch03.1.4} were given in Pinsky (1972). Other results
in this section can be found in Li (2000).

\section{Stationary immigration distributions}\label{sch03.2}

 \setcounter{equation}{0}

In this section, we give a brief discussion of the structures of
stationary distributions of the CBI-processes. The results here were
first given in Li (2002) in the setting of measure-valued processes.
Given two probability measures $\eta_1$ and $\eta_2$ on $\mbb{R}_+$, we
write $\eta_1\preceq \eta_2$ if $\eta_1*\gamma = \eta_2$ for another
probability measure $\gamma$ on $\mbb{R}_+$. Clearly, the measure
$\gamma$ is unique if it exists. Let $(Q_t)_{t\ge 0}$ be the transition
semigroup defined by \eqref{ch02.1.15} and \eqref{ch02.1.20}, where
$(v_t)_{t\ge 0}$ has the representation \eqref{ch02.2.1}. Let
$(Q_t^\circ)_{t\ge 0}$ be the restriction of $(Q_t)_{t\ge 0}$ to
$(0,\infty)$. Let $\mcr{E}^*(Q)$ denote the set of probabilities $\eta$
on $\mbb{R}_+$ satisfying $\eta Q_t\preceq \eta$ for all $t\ge 0$.

 \bgtheorem\label{tch03.2.1}
For each $\eta\in \mcr{E}^*(Q)$ there is a unique SC-semigroup
$(\gamma_t)_{t\ge 0}$ associated with $(Q_t)_{t\ge 0}$ such that $\eta
Q_t*\gamma_t = \eta$ for $t\ge 0$ and $\eta = \lim_{t\to
\infty}\gamma_t$. \edtheorem

 \begin{proof}
By the definition of $\mcr{E}^*(Q)$, for each $t\ge 0$ there is a unique
probability measure $\gamma_t$ on $\mbb{R}_+$ satisfying $\eta = (\eta
Q_t) * \gamma_t$. By the branching property of $(Q_t)_{t\ge 0}$ one can
show $(\mu_1*\mu_2)Q_t = (\mu_1Q_t)*(\mu_2Q_t)$ for any $t\ge 0$ and any
probability measures $\mu_1$ and $\mu_2$ on $\mbb{R}_+$. Then for $r,t\ge
0$ we have
 \beqnn
(\eta Q_{r+t})*\gamma_{r+t}
 =
(\eta Q_t)*\gamma_t
 =
\{[(\eta Q_r)*\gamma_r]Q_t\}*\gamma_t
 =
(\eta Q_{r+t})*(\gamma_rQ_t)*\gamma_t.
 \eeqnn
A cancelation gives \eqref{ch03.1.1}, so $(\gamma_t)_{t\ge 0}$ is an
SC-semigroup associated with $(Q_t)_{t\ge 0}$. Now for every $\lambda\ge
0$ the function $t\mapsto L_{\gamma_t}(\lambda)$ is decreasing. By the
relation $\eta = (\eta Q_t)*\gamma_t$ one can see $t\mapsto L_{\eta
Q_t}(\lambda)$ is increasing. Then there are probability measures
$\eta_i$ and $\eta_p$ on $\mbb{R}_+$ so that $\eta_i*\eta_p = \eta$ and
 \beqnn
L_{\eta_i}(\lambda) = \lim_{t\to \infty}L_{\eta Q_t}(\lambda),
 \quad
L_{\eta_p}(\lambda) = \lim_{t\to \infty}L_{\gamma_t}(\lambda).
 \eeqnn
These imply $\eta_i = \lim_{t\to\infty} \eta Q_t$ and $\eta_p =
\lim_{t\to\infty} \gamma_t$. It follows that $\eta_i$ is a stationary
distribution of $(Q_t)_{t\ge 0}$, so we must have $\eta_i = \delta_0$ and
$\eta_p = \eta$. \qed\end{proof}

In the situation of Theorem~\ref{tch03.2.1}, it is easy to see the
measure $\eta\in \mcr{E}^*(Q)$ is the unique stationary distribution of
the transition semigroup $(Q_t^\gamma)_{t\ge 0}$ defined by
\eqref{ch03.1.4}. Then we can identify $\mcr{E}^*(Q)$ with the set of
stationary distributions of immigration processes associated with
$(Q_t)_{t\ge 0}$. As a consequence of Theorem~\ref{tch03.1.3}
and~\ref{tch03.2.1}, every $\mu\in \mcr{E}^*(Q)$ is infinitely divisible.
Recall that we write $\mu = I(h,l)$ if $\mu$ is an infinitely divisible
probability measure on $\mbb{R}_+$ with $\psi := -\log L_\mu$ given by
\eqref{ch01.2.3}. Let $\mcr{E}(Q^\circ)$ denote the set of excessive
measures $\nu$ for $(Q_t^\circ)_{t\ge 0}$ satisfying
 \beqnn
\int_0^\infty (1\land u) \nu(\d u)< \infty.
 \eeqnn
The following result gives some characterizations of the set
$\mcr{E}^*(Q)$.

 \bgtheorem\label{tch03.2.3}
Let $\eta = I(\beta,\nu)$ be an infinitely divisible probability measure
on $\mbb{R}_+$. Then $\eta\in \mcr{E}^*(Q)$ if and only if $(\beta,\nu)$
satisfy
 \beqlb\label{ch03.2.1}
\beta h_t\le \beta
 \quad\mbox{and}\quad
\beta l_t + \nu Q_t^\circ\le \nu, \qquad t\ge 0.
 \eeqlb
In particular, if $\nu\in \mcr{E}(Q^\circ)$, then $\eta = I(0,\nu) \in
\mcr{E}^*(Q)$. \edtheorem

 \begin{proof}
It is easy to show that $\eta Q_t = I(\beta_t,\nu_t)$, where $\beta_t =
\beta h_t$ and $\nu_t = \beta l_t + \nu Q_t^\circ$. Then $\eta Q_t\preceq
\eta$ holds if and only if \eqref{ch03.2.1} is satisfied. The second
assertion is immediate. \qed\end{proof}

\section{Scaling limits of discrete immigration models}\label{sch03.3}

 \setcounter{equation}{0}

In this section, we prove a limit theorem of rescaled Galton--Watson
branching processes with immigration, which leads to the CBI-processes.
This kind of limit theorems were studied in Aliev and Shchurenkov (1982),
Kawazu and Watanabe (1971) and Li (2006) among many others.

Let $g$ and $h$ be two probability generating functions. Suppose that
$\{\xi_{n,i}: n,i=1,2,\ldots\}$ and $\{\eta_n: n=1,2,\ldots\}$ are
independent families of positive integer-valued i.i.d.\ random variables
with distributions given by $g$ and $h$, respectively. Given another
positive integer-valued random variable $y(0)$ independent of
$\{\xi_{n,i}\}$ and $\{\eta_n\}$, we define inductively
 \beqlb\label{ch03.3.1}
y(n) = \sum_{i=1}^{y(n-1)}\xi_{n,i} + \eta_n, \qquad n=1,2,\ldots.
 \eeqlb
Then $\{y(n): n = 0,1,2,\ldots\}$ is a discrete-time positive
integer-valued Markov chain with transition matrix $Q(i,j)$ determined by
 \beqlb\label{ch03.3.2}
\sum^\infty_{j=0} Q(i,j)z^j = g(z)^ih(z), \qquad |z|\le 1.
 \eeqlb
The random variable $y(n)$ can be thought of as the number of individuals
in generation $n\ge 0$ of an evolving particle system. After one unit
time, each of the $y(n)$ particles splits independently of others into a
random number of offspring according to the distribution given by $g$ and
a random number of immigrants are added to the system according to the
probability law given by $h$. The $n$-step transition matrix $Q^n(i,j)$
of $\{y(n): n = 0,1,2,\ldots\}$ is given by
 \beqlb\label{ch03.3.3}
\sum^\infty_{j=0} Q^n(i,j)z^j = g^n(z)^i\prod_{j=1}^nh(g^{j-1}(z)),
\qquad |z|\le 1,
 \eeqlb
where $g^n(z)$ is defined by $g^n(z) = g(g^{n-1}(z))$ successively with
$g^0(z) = z$. We call any positive integer-valued Markov chain with
transition probabilities given by \eqref{ch03.3.2} or \eqref{ch03.3.3} a
\index{Galton--Watson branching process with immigration}
\textit{Galton--Watson branching process with immigration}
\index{GWI-process} (GWI-process) with parameters $(g,h)$. When $h\equiv
1$, this reduces to the GW-process defined in the first section.

Suppose that for each integer $k\ge 1$ we have a GWI-process $\{y_k(n):
n\ge 0\}$ with parameters $(g_k,h_k)$. Let $z_k(n) = y_k(n)/k$. Then
$\{z_k(n): n\ge 0\}$ is a Markov chain with state space $E_k :=
\{0,1/k,2/k,\ldots\}$ and $n$-step transition probability $Q_k^n(x,\d y)$
determined by
 \beqlb\label{ch03.3.4}
\int_{E_k}\e^{-\lambda y}Q_k^n(x,\d y)
 =
g_k^n(\e^{-\lambda/k})^{kx}\prod_{j=1}^nh(g_k^{j-1}(\e^{-\lambda/k})),
\quad \lambda\ge 0.
 \eeqlb
Suppose that $\{\gamma_k\}$ is a positive real sequence so that
$\gamma_k\to \infty$ increasingly as $k\to \infty$. Let $[\gamma_kt]$
denote the integer part of $\gamma_kt\ge 0$. In view of \eqref{ch03.3.4},
given $z_k(0) = x$ the conditional distribution
$Q_k^{[\gamma_kt]}(x,\cdot)$ of $z_k([\gamma_kt])$ on $E_k$ is determined
by
 \beqlb\label{ch03.3.5}
\ar\ar\int_{E_k}\e^{-\lambda y}Q_k^{[\gamma_kt]}(x,\d y) \cr
 \ar\ar\qquad
= \exp\Big\{-xv_k(t,\lambda) - \int_0^{\frac{[\gamma_kt]}{\gamma_k}}
\bar{H}_k(v_k(s,\lambda))\d s\Big\},
 \eeqlb
where $v_k(t,\lambda)$ is given by \eqref{ch02.1.7} and
 \beqnn
\bar{H}_k(\lambda)
 =
-\gamma_k\log h_k(\e^{-\lambda/k}), \qquad \lambda\ge 0.
 \eeqnn
For any $z\ge 0$ let $G_k(z)$ be defined by \eqref{ch02.1.8} and let
 \beqlb\label{ch03.3.6}
H_k(z)=\gamma_k[1-h_k(\e^{-z/k})].
 \eeqlb

 \bgcondition\label{tch03.3.1}
There is a function $\psi$ on $[0,\infty)$ such that $H_k(z)\to \psi(z)$
uniformly on $[0,a]$ for every $a\ge 0$ as $k\to \infty$. \edcondition

It is simple to see that $H_k\in \mcr{I}$. By Theorem~\ref{tch01.2.2}, if
the above condition is satisfied, the limit function $\psi$ has the
representation \eqref{ch03.1.10}. A different of proof of the following
theorem was given in Li (2006).

 \bgtheorem\label{tch03.3.2}
Suppose that Conditions~\ref{tch02.1.2} and~\ref{tch03.3.1} are
satisfied. Let $\{y(t): t\ge 0\}$ be a CBI-process with transition
semigroup $(Q_t^\gamma)_{t\ge 0}$ defined by \eqref{ch03.1.9}. If
$z_k(0)$ converges to $y(0)$ in distribution, then $\{z_k([\gamma_kt]):
t\ge 0\}$ converges to $\{y(t): t\ge 0\}$ in distribution on
$D([0,\infty), \mbb{R}_+)$. \edtheorem

 \begin{proof}
By Theorem~\ref{tch02.1.6} for every $a\ge 0$ we have $v_k(t,\lambda)\to
v_t(\lambda)$ uniformly on $[0,a]^2$ as $k\to \infty$. For $\lambda>0$
and $x\ge 0$ set $e_\lambda(x) = \e^{-\lambda x}$. In view of
\eqref{ch03.3.5} we have
 \beqnn
\lim_{k\to\infty} \sup_{x\in E_k}\big|Q_k^{[\gamma_kt]}e_\lambda(x) -
Q_t^\gamma e_\lambda(x)\big| = 0
 \eeqnn
for every $t\ge 0$. Then the result follows as in the proof of
Theorem~\ref{tch02.1.9}. \qed\end{proof}

 \bexample\label{ech03.3.1}
In a special case, we can give a characterization for the CBI-process in
terms of a stochastic differential equation. Let $m =
\mbf{E}[\xi_{1,1}]$. From \eqref{ch03.3.1} we have
 \beqnn
y(n) - y(n-1)
 =
\sqrt{y(n-1)}\sum_{i=1}^{y(n-1)}\frac{\xi_{n,i} - m}{\sqrt{y(n-1)}} -
(1-m)y(n-1) + \eta_n.
 \eeqnn
Then it is natural to expect that a typical CBI-process would solve the
stochastic differentia equation
 \beqlb\label{ch03.3.7}
\d y(t)
 =
\sqrt{2cy(t)}\d B(t) - by(t)\d t + \beta\d t, \qquad t\ge 0,
 \eeqlb
where $\{B(t): t\ge 0\}$ is a Brownian motion. The above equation has a
unique positive strong solution; see Ikeda and Watanabe (1989,
pp.235--236). In fact, the solution $\{y(t): t\ge 0\}$ has transition
semigroup given by \eqref{ch03.1.17} with $\alpha=1$. Let
$C^2(\mbb{R}_+)$ denote the set of bounded continuous real functions on
$\mbb{R}_+$ with bounded continuous derivatives up to the second order.
Then $\{y(t): t\ge 0\}$ has generator $A$ given by
 \beqnn
Af(x) = c\frac{\d^2}{\d x^2}f(x) + (\beta - bx)\frac{\d}{\d x}f(x),
\qquad f\in C^2(\mbb{R}_+).
 \eeqnn
In particular, for $\beta=0$ the solution of \eqref{ch03.3.7} is called
\index{Feller's branching diffusion} \textit{Feller's branching
diffusion}. \eexample

\section{A reconstruction of the sample paths}\label{sch03.4}

 \setcounter{equation}{0}

In this section, we give a reconstruction of the sample paths of the
CBI-process using a Poisson random measure. Let us consider an entrance
law $(H_t)_{t>0}$ for $(Q_t^\circ)_{t\ge 0}$ and let $\mbf{Q}_H$ be the
$\sigma$-finite measure on $D(0,\infty)$ determined by \eqref{ch02.4.1}.
Suppose that $N(\d s,\d w)$ is a Poisson random measure on
$(0,\infty)\times D(0,\infty)$ with intensity $\d s\mbf{Q}_H(\d w)$. We
define the process
 \beqlb\label{ch03.4.1}
Y_t = \int_{(0,t)} \int_{D(0,\infty)} w_{t-s} N(\d s,\d w), \qquad t\ge
0.
 \eeqlb
A special form of the following theorem was established in Pitman and Yor
(1982):

 \bgtheorem\label{tch03.4.1}
For $t\ge 0$ let $\mcr{G}_t$ be the $\sigma$-algebra generated by the
collection of random variables $\{N((0,u]\times A): A\in \mcr{A}^0_{t-u},
0\le u< t\}$. Then $\{(Y_t,\mcr{G}_t): t\ge 0\}$ is an immigration
process with transition semigroup $(Q^H_t)_{t\ge 0}$ given by
 \beqlb\label{ch03.4.2}
\int_0^\infty \e^{-\lambda y} Q_t^H(x,\d y)
 =
\exp\Big\{-xv_t(\lambda) - \int_0^t\d s\int_0^\infty(1-\e^{-\lambda y})
H_s(\d y)\Big\}.
 \eeqlb
\edtheorem

 \begin{proof}
It is easy to show that $Y_t$ has distribution $Q_t^H(0,\cdot)$ on
$\mbb{R}_+$. Let $t\ge r>u\ge 0$ and let $h$ be a bounded positive
function on $D(0,\infty)$ measurable relative to $\mcr{A}^0_{r-u}$. For
$\lambda\ge 0$ we can see as in the proof of Theorem~\ref{tch02.4.2} that
 \beqnn
 \ar\ar\mbf{P}\Big[\exp\Big\{-\int_0^u\int_{D(0,\infty)}h(w)N(\d s,\d w)
- \lambda Y_t\Big\}\Big] \cr
 \ar\ar\qquad
=\, \mbf{P}\Big[\exp\Big\{-\int_0^t\int_{D(0,\infty)} \big[h(w) 1_{\{s\le
u\}} + \lambda w_{t-s}\big] N(\d s,\d w)\Big\}\Big] \cr
 \ar\ar\qquad
=\, \exp\Big\{-\int_0^t\d s\int_{D(0,\infty)}\big(1 - \e^{-h(w) 1_{\{s\le
u\}}}\e^{-\lambda w_{t-s}}\big)\mbf{Q}_H(\d w)\Big\} \cr
 \ar\ar\qquad
=\, \exp\Big\{-\int_0^u\d s\int_{D(0,\infty)}\big(1 - \e^{-h(w)}
\e^{-v_{t-r}(\lambda)w_{r-s}}\big)\mbf{Q}_H(\d w)\Big\} \cr
 \ar\ar\qquad\qquad
\cdot\exp\Big\{-\int_u^r\d s\int_{D(0,\infty)}\big(1 -
\e^{-v_{t-r}(\lambda)w_{r-s}}\big)\mbf{Q}_H(\d w)\Big\} \cr
 \ar\ar\qquad\qquad
\cdot\exp\Big\{-\int_r^t\d s\int_{D(0,\infty)}\big(1 - \e^{-\lambda
w_{t-s}}\big) \mbf{Q}_H(\d w)\Big\} \cr
 \ar\ar\qquad
=\, \mbf{P}\Big[\exp\Big\{-\int_0^u\int_{D(0,\infty)} h(w) N(\d s,\d w) -
Y_rv_{t-r}(\lambda)\Big\}\Big] \cr
 \ar\ar\qquad\qquad
\cdot\exp\Big\{-\int_r^t\d s\int_0^\infty(1-\e^{-\lambda y})H_{t-s}(\d
y)\Big\},
 \eeqnn
where we have used the Markov property \eqref{ch02.4.1} for the third
equality. That shows $\{(Y_t,\mcr{G}_t): t\ge 0\}$ is a Markov process in
$\mbb{R}_+$ with transition semigroup $(Q^H_t)_{t\ge 0}$. \qed\end{proof}

Let $\mbf{P}_x(\d w)$ denote the distribution on $D[0,\infty)$ of the
CB-process $\{x(t): t\ge 0\}$ with $x(0)=x$. Suppose that $\psi\in
\mcr{I}$ be given by \eqref{ch03.1.10}. If the condition
\eqref{ch02.2.14} is satisfied, we can define a $\sigma$-finite measure
$\mbf{Q}_H(\d w)$ on $D[0,\infty)$ by
 \beqlb\label{ch03.4.3}
\mbf{Q}_H(\d w) = \beta \mbf{Q}_{(0)}(\d w) + \int_0^\infty n(\d x)
\mbf{P}_x(\d w).
 \eeqlb
This corresponds to the entrance law $(H_t)_{t>0}$ for $(Q_t^\circ)_{t\ge
0}$ defined by
 \beqlb\label{ch03.4.4}
H_t = \beta l_t + \int_0^\infty n(\d x) Q_t(x,\cdot), \qquad t>0.
 \eeqlb
In this case, it is easy to show that
 \beqlb\label{ch03.4.5}
\int_0^\infty(1-\e^{-\lambda y}) H_t(\d y)
 =
\psi(v_t(\lambda)), \quad t>0, \lambda\ge 0.
 \eeqlb
Then from Theorem~\ref{tch03.4.1} we obtain

 \bgcorollary\label{tch03.4.2}
Suppose that \eqref{ch02.2.14} is satisfied and $(H_t)_{t>0}$ is given by
\eqref{ch03.4.3}. Then $\{(Y_t,\mcr{G}_t): t\ge 0\}$ is an immigration
process with transition semigroup $(Q^\gamma_t)_{t\ge 0}$ given by
\eqref{ch03.1.9}. \edcorollary


\chapter{Martingale problems and stochastic equations}\label{ch04}

Martingale problems play a very important role in the study of Markov
processes. In this chapter we prove the equivalence of a number of
martingale problems for CBI-processes. From the martingale problems we
derive some stochastic equations. Using the stochastic equations, we give
simple proofs of Lamperti's transformations on CB-processes and
spectrally positive L\'{e}vy processes.

\section{Martingale problem formulations}\label{sch04.1}

 \setcounter{equation}{0}

In this section we give some formulations of the CBI-process in terms of
martingale problems and prove their equivalence. Suppose that
$(\phi,\psi)$ are given respectively by \eqref{ch02.1.13} and
\eqref{ch03.1.10} with $un(\d u)$ being a finite measure on $(0,\infty)$.
For $f\in C^2(\mbb{R}_+)$ define
 \beqlb\label{ch04.1.1}
Lf(x) \ar=\ar cxf^{\prime\prime}(x) + x\int_0^\infty \big[f(x+z) - f(x) -
zf^\prime(x)\big] m(\d z) \cr
 \ar~\ar
+\, (\beta - bx)f^\prime(x) + \int_0^\infty \big[f(x+z) - f(x)\big] n(\d
z).
 \eeqlb
We shall identify the operator $L$ as the generator of the CBI-process.
For this purpose we need the following:

 \bgproposition\label{tch04.1.1}
Let $(Q_t^\gamma)_{t\ge 0}$ be the transition semigroup defined by
\eqref{ch02.1.20} and \eqref{ch03.1.9}. Then for any $t\ge 0$ and
$\lambda\ge 0$ we have
 \beqlb\label{ch04.1.2}
\int_0^\infty \e^{-\lambda y} Q_t^\gamma(x,\d y)
 =
\e^{-x\lambda} + \int_0^t \d s\int_0^\infty [x\phi(\lambda) -
\psi(\lambda)] \e^{-y\lambda} Q_s^\gamma(x,\d y).
 \eeqlb
\edproposition

 \begin{proof}
Recall that $v_t^\prime(\lambda) = (\partial/\partial \lambda)
v_t(\lambda)$. By differentiating both sides of \eqref{ch03.1.9} we get
 \beqnn
\int_0^\infty y\e^{-y\lambda} Q_t^\gamma(x,\d y)
 \ar=\ar
\int_0^\infty \e^{-y\lambda} Q_s^\gamma(x,\d y)\Big[xv_t^\prime(\lambda)
+ \int_0^t \psi^\prime(v_s(\lambda))v_s^\prime(\lambda)\d s\Big].
 \eeqnn
From this and \eqref{ch02.2.3} it follows that
 \beqnn
{\partial\over\partial t}\int_0^\infty \e^{-y\lambda} Q_t^\gamma(x,\d y)
 \ar=\ar
- \Big[x{\partial\over \partial t}v_t(\lambda) + \psi(v_t(\lambda))\Big]
\int_0^\infty \e^{-y\lambda} Q_s^\gamma(x,\d y) \cr
 \ar=\ar
\Big[x\phi(\lambda) v_t^\prime(\lambda) - \psi(\lambda)\Big]\int_0^\infty
\e^{-y\lambda} Q_s^\gamma(x,\d y) \cr
 \ar\ar
- \int_0^t \psi^\prime(v_s(\lambda)){\partial\over \partial
s}v_s(\lambda) \d s \int_0^\infty \e^{-y\lambda} Q_s^\gamma(x,\d y) \cr
 \ar=\ar
\Big[x\phi(\lambda) v_t^\prime(\lambda) - \psi(\lambda)\Big]\int_0^\infty
\e^{-y\lambda} Q_s^\gamma(x,\d y) \cr
 \ar\ar
+\, \phi(\lambda)\int_0^t \psi^\prime(v_s(\lambda))v_s^\prime(\lambda) \d
s \int_0^\infty \e^{-y\lambda} Q_s^\gamma(x,\d y) \cr
 \ar=\ar
\int_0^\infty [x\phi(\lambda) - \psi(\lambda)] \e^{-y\lambda}
Q_s^\gamma(x,\d y).
 \eeqnn
That gives \eqref{ch04.1.2}. \qed\end{proof}

Suppose that $(\Omega, \mcr{G}, \mcr{G}_t, \mbf{P})$ is a filtered
probability space satisfying the usual hypotheses and $\{y(t): t\ge 0\}$
is a c\`{a}dl\`{a}g process in $\mbb{R}_+$ that is adapted to
$(\mcr{G}_t)_{t\ge 0}$ and satisfies $\mbf{P}[y(0)]< \infty$. Let us
consider the following properties:
 \benumerate

\itm[{\rm(1)}] For every $T\ge 0$ and $\lambda\ge 0$,
 \beqnn
\exp\Big\{-v_{T-t}(\lambda)y(t) - \int_0^{T-t} \psi(v_s(\lambda))\d
s\Big\}, \qquad 0\le t\le T,
 \eeqnn
is a martingale.

\itm[{\rm(2)}] For every $\lambda\ge 0$,
 \beqnn
H_t(\lambda)
 :=
\exp\Big\{-\lambda y(t)+\int_0^t[\psi(\lambda) - y(s)\phi(\lambda)]\d
s\Big\}, \quad t\ge 0,
 \eeqnn
is a local martingale.

\itm[{\rm(3)}] {\rm(a)} The process $\{y(t): t\ge 0\}$ has no negative
jumps. Let $N(\d s,\d z)$ be the optional random measure on
$(0,\infty)^2$ defined by
 \beqnn
N(\d s,\d z) = \sum_{s>0}1_{\{\Delta y(s)\ne 0\}}\delta_{(s,\Delta
y(s))}(\d s,\d z),
 \eeqnn
where $\Delta y(s) = y(s) - y(s-)$, and let $\hat{N}(\d s,\d z)$ denote
the predictable compensator of $N(\d s,\d z)$. Then $\hat{N}(\d s,\d z) =
y(s-)\d sm(\d z) + \d sn(\d z)$.

{\rm(b)} If we let $\tilde{N}(\d s,\d z) = N(\d s,\d z) - \hat{N}(\d s,\d
z)$, then
 \beqnn
y(t) = y(0) + M_t^c + M_t^d + \int_0^t \Big[\beta + \int_0^\infty z n(\d
z) - by(s)\Big]\d s,
 \eeqnn
where $t\mapsto M_t^c$ is a continuous local martingale with quadratic
variation $2cy(t-)\d t$ $=$ $2cy(t)\d t$ and
 \beqnn
t\mapsto M_t^d = \int_0^t\int_0^\infty z \tilde{N}(\d s,\d z)
 \eeqnn
is a purely discontinuous local martingale.

\itm[{\rm(4)}] For every $f\in C^2(\mbb{R}_+)$ we have
 \beqnn
f(y(t)) = f(y(0)) + \int_0^t Lf(y(s))\d s + \mbox{local mart.}
 \eeqnn

 \eenumerate

 \bgtheorem\label{tch04.1.2}
The above properties {\rm(1)}, {\rm(2)}, {\rm(3)} and {\rm(4)} are
equivalent to each other. Those properties hold if and only if
$\{(y(t),\mcr{G}_t): t\ge 0\}$ is a CBI-process with parameters
$(\phi,\psi)$. \edtheorem

 \begin{proof}
Clearly, {\rm(1)} holds if and only if $\{y(t): t\ge 0\}$ is a Markov
process relative to $(\mcr{G}_t)_{t\ge 0}$ with transition semigroup
$(Q_t^\gamma)_{t\ge 0}$ defined by \eqref{ch03.1.9}. Then we only need to
prove the equivalence of the four properties.

{\rm(1)}$\Rightarrow${\rm(2)}: Suppose that {\rm(1)} holds. Then $\{y(t):
t\ge 0\}$ is a CBI-process with transition semigroup $(Q_t^\gamma)_{t\ge
0}$ given by \eqref{ch03.1.9}. By \eqref{ch04.1.2} and the Markov
property it is easy to see that
 \beqnn
Y_t(\lambda) := \e^{-\lambda y(t)} + \int_0^t [\psi(\lambda) -
y(s)\phi(\lambda)] \e^{-\lambda y(s)}\d s
 \eeqnn
is a martingale. By integration by parts applied to
 \beqlb\label{ch04.1.3}
Z_t(\lambda) := \e^{-\lambda y(t)}
 ~~\mbox{and}~~
W_t(\lambda) := \exp\Big\{\int_0^t[\psi(\lambda) - y(s)\phi(\lambda)]\d
s\Big\}
 \eeqlb
we obtain
 \beqnn
\d H_t(\lambda)
 =
\e^{-\lambda y(t-)}\d W_t(\lambda) + W_t(\lambda)\d\e^{-\lambda y(t)}
 =
W_t(\lambda)\d Y_t(\lambda).
 \eeqnn
Then $\{H_t(\lambda)\}$ is a local martingale.

{\rm(2)}$\Rightarrow${\rm(3)}: For any $\lambda\ge 0$ define
$Z_t(\lambda)$ and $W_t(\lambda)$ by \eqref{ch04.1.3}. We have
$Z_t(\lambda) = H_t(\lambda)W_t(\lambda)^{-1}$ and so
 \beqlb\label{ch04.1.4}
\d Z_t(\lambda) = W_t(\lambda)^{-1}\d H_t(\lambda) -
Z_{t-}(\lambda)[\psi(\lambda) - y(t-)\phi(\lambda)]\d t
 \eeqlb
by integration by parts. Then $\{Z_t(\lambda)\}$ is a special
semi-martingale; see, e.g., Dellacherie and Meyer (1982, p.213). By
It\^{o}'s formula we find the $\{y(t)\}$ is also a special
semi-martingale. We define the optional random measure $N(\d s,\d z)$ on
$[0,\infty)\times \mbb{R}$ by
 \beqnn
N(\d s,\d z) = \sum_{s>0}1_{\{\Delta y(s)\ne 0\}} \delta_{(s,\Delta
y(s))}(\d s,\d z),
 \eeqnn
where $\Delta y(s) = y(s) - y(s-)$. Let $\hat{N}(\d s,\d z)$ denote the
predictable compensator of $N(\d s,\d z)$ and let $\tilde{N}(\d s,\d z)$
denote the compensated random measure; see Dellacherie and Meyer (1982,
pp.371--374). It follows that
 \beqlb\label{ch04.1.5}
y(t) = y(0) + U_t + M_t^c + M_t^d,
 \eeqlb
where $\{U_t\}$ is a predictable process with locally bounded variations,
$\{M^c_t\}$ is a continuous local martingale and
 \beqlb\label{ch04.1.6}
M_t^d = \int_0^t\int_{\mbb{R}} z \tilde{N}(\d s,\d z), \quad t\ge 0,
 \eeqlb
is a purely discontinuous local martingale; see Dellacherie and Meyer
(1982, p.353 and p.376) or Jacod and Shiryaev (2003, p.84). Let $\{C_t\}$
denote the quadratic variation process of $\{M^c_t\}$. By It\^{o}'s
formula,
 \beqlb\label{ch04.1.7}
Z_t(\lambda) \ar=\ar Z_0(\lambda) - \lambda\int_0^tZ_{s-}(\lambda)\d U_s
+ \frac{1}{2}\lambda^2\int_0^t Z_{s-}(\lambda)\d C_s \cr
 \ar~\ar
+ \int_0^t\int_{\mbb{R}} Z_{s-}(\lambda) \big(\e^{-z\lambda} - 1 +
z\lambda\big) \hat{N}(\d s,\d z) + \mbox{local mart.}
 \eeqlb
In view of \eqref{ch04.1.4} and \eqref{ch04.1.7} we get
 \beqnn
[y(t)\phi(\lambda) - \psi(\lambda)]\d t
 =
\frac{1}{2}\lambda^2\d C_t - \lambda\d U_t + \int_{\mbb{R}}
\big(\e^{-z\lambda} - 1 + z\lambda\big) \hat{N}(\d t,\d z)
 \eeqnn
by the uniqueness of canonical decompositions of special
semi-martingales; see Dellacherie and Meyer (1982, p.213). By
substituting the representation \eqref{ch02.1.13} of $\phi$ into the
above equation and comparing both sides it is easy to find that
{\rm(3.a)} and {\rm(3.b)} hold.

{\rm(3)}$\Rightarrow${\rm(4)}: This follows by It\^{o}'s formula.

{\rm(4)}$\Rightarrow${\rm(1)}: Let $G = G(t,x)\in C^{1,2}([0,T]\times
\mbb{R}_+)$. For $0\le t\le T$ and $k\ge1$ we have
 {\small\beqnn
G(t,y(t))
 \ar=\ar
G(0,y(0)) + \sum_{j=0}^\infty \big[G(t\land j/k,y(t\land (j+1)/k)) -
G(t\land j/k,y(t\land j/k))\big] \cr
 \ar~\ar
+ \sum_{j=0}^\infty \big[G(t\land (j+1)/k,y(t\land (j+1)/k)) - G(t\land
j/k,y(t\land (j+1)/k))\big],
 \eeqnn\!}
where the summations only consist of finitely many non-trivial terms. By
applying {\rm(4)} term by term we obtain
 \beqnn
G(t,y(t))
 \ar=\ar
G(0,y(0)) + \sum_{j=0}^\infty \int_{t\land j/k}^{t\land (j+1)/k}
\Big\{[\beta - by(s)]G^\prime_y(t\land j/k,y(s)) \cr
 \ar~\ar
+\, cy(s)G^{\prime\prime}_{xx}(t\land j/k,y(s)) + y(s)\int_0^\infty
\Big[G(t\land j/k,y(s) + z) \cr
 \ar~\ar
-\, G(t\land j/k,y(s)) - z G^\prime_y(t\land j/k,y(s))\Big] m(\d z) \cr
 \ar~\ar
+ \int_0^\infty \Big[G(t\land j/k,y(s) + z) - G(t\land j/k,y(s))\Big]
n(\d z)\Big\}\d s \cr
 \ar~\ar
+ \sum_{j=0}^\infty\int_{t\land j/k}^{t\land (j+1)/k}
G^\prime_t(s,y(t\land (j+1)/k))\d s + M_k(t),
 \eeqnn
where $\{M_k(t)\}$ is a local martingale. Since $\{y(t)\}$ is a
c\`{a}dl\`{a}g process, letting $k\to\infty$ in the equation above gives
 \beqnn
G(t,y(t)) \ar=\ar G(0,y(0)) + \int_0^t \Big\{G^\prime_t(s,y(s)) -
by(s)G^\prime_y(s,y(s)) \cr
 \ar~\ar
+\, cy(s)G^{\prime\prime}_{xx}(s,y(s)) + y(s)\int_0^\infty
\Big[G(s,y(s)+z) \cr
 \ar~\ar
-\, G(s,y(s)) - z G^\prime_x(s,y(s))\Big]m(\d z) \cr
 \ar~\ar
+ \int_0^\infty\Big[G(s,y(s)+z) - G(s,y(s))\Big]n(\d z)\Big\}\d s + M(t),
 \eeqnn
where $\{M(t)\}$ is a local martingale. For any $T\ge 0$ and $\lambda\ge
0$ we may apply the above to
 \beqnn
G(t,x) = \exp\Big\{-v_{T-t}(\lambda)x-\int_0^{T-t}\psi(v_s(\lambda))\d
s\Big\}
 \eeqnn
to see $t\mapsto G(t,y(t))$ is a local martingale. \qed\end{proof}

The above property {\rm(4)} implies that the generator of the CBI-process
is the closure of the generator $L$ in the sense of Ethier and Kurtz
(1986). This explicit form of the generator was first given in Kawazu and
Watanabe (1971).

\section{Stochastic equations of CBI-processes}\label{sch04.2}

 \setcounter{equation}{0}

In this section we establish some stochastic equations for the
CBI-processes. The reader may refer to Dawson and Li (2006, 2010) and Fu
and Li (2010) for more results on this topic. Suppose that $(\phi,\psi)$
are given respectively by \eqref{ch02.1.13} and \eqref{ch03.1.10} with
$un(\d u)$ being a finite measure on $(0,\infty)$. Let
$(Q_t^\gamma)_{t\ge 0}$ be the transition semigroup defined by
\eqref{ch02.1.20} and \eqref{ch03.1.9}. In this section, we derive some
stochastic equations for the CBI-processes.

Suppose that $(\Omega, \mcr{F}, \mcr{F}_t, \mbf{P})$ is a filtered
probability space satisfying the usual hypotheses. Let $\{B(t): t\ge 0\}$
be an $(\mcr{F}_t)$-Brownian motion and let $\{p_0(t): t\ge 0\}$ and
$\{p_1(t): t\ge 0\}$ be $(\mcr{F}_t)$-Poisson point processes on
$(0,\infty)^2$ with characteristic measures $m(dz)du$ and $n(dz)du$,
respectively. We assume that the white noise and the Poisson processes
are independent of each other. Let $N_0(ds,dz,du)$ and $N_1(ds,dz,du)$
denote the Poisson random measures on $(0,\infty)^3$ associated with
$\{p_0(t)\}$ and $\{p_1(t)\}$, respectively. Let $\tilde{N}_0(ds,dz,du)$
denote the compensated measure of $N_0(ds,dz,du)$. Let us consider the
stochastic integral equation
 \beqlb\label{ch04.2.1}
y(t)
 \ar=\ar
y(0) + \int_0^t \sqrt{2cy(s)}\d B(s) + \int_0^t\int_0^\infty
\int_0^{y(s-)} z \tilde{N}_0(\d s,\d z,\d u) \cr
 \ar~\ar\qquad
+ \int_0^t (\beta-by(s))\d s + \int_0^t\int_0^\infty z N_1(\d s,\d z),
 \eeqlb
where $\tilde{N}_0(\d s,\d z,\d u) = N_0(\d s,\d z,\d u) - \d sm(\d z)\d
u$. We understand the last term on the right-hand side as an integral
over the set $\{(s,z,u): 0<s\le t, 0<z<\infty, 0<u\le y(s-)\}$ and give
similar interpretations for other integrals with respect to Poisson
random measures in this section.

 \bgtheorem\label{tch04.2.1}
There is a unique positive weak solution to \eqref{ch04.2.1} and the
solution is a CBI-process with transition semigroup $(Q_t^\gamma)_{t\ge
0}$. \edtheorem

 \begin{proof}
Suppose that $\{y(t)\}$ is a c\`{a}dl\`{a}g realization of the
CBI-process with transition semigroup given by \eqref{ch02.1.20} and
\eqref{ch03.1.9}. By Theorem~\ref{tch04.1.2} the process has no negative
jumps and the random measure
 \beqnn
N(\d s,\d z) := \sum_{s>0}1_{\{y(s)\neq y(s-)\}}
\delta_{(s,y(s)-y(s-))}(\d s,\d z)
 \eeqnn
has predictable compensator
 \beqnn
\hat{N}(\d s,\d z) = y(s-)\d sm(\d z) + \d sn(\d z)
 \eeqnn
and
 \beqlb\label{ch04.2.2}
y(t) \ar=\ar y(0) + t\Big[\beta + \int_0^\infty un(\d u)\Big] - \int_0^t
by(s-)\d s \cr
 \ar~\ar
+\, M^c(t) + \int_0^t \int_0^\infty z \tilde{N}(\d s,\d z),
 \eeqlb
where $\tilde{N}(\d s,\d z) = N(\d s,\d z) - \hat{N}(\d s,\d z)$ and
$t\mapsto M^c(t)$ is a continuous local martingale with quadratic
variation $2cy(t-)\d t$. By representation theorems for semimartingales,
we have equation \eqref{ch04.2.1} on an extension of the original
probability space; see, e.g., Ikeda and Watanabe (1989, p.90 and p.93).
That proves the existence of a weak solution to \eqref{ch04.2.1}.
Conversely, if $\{y(t)\}$ is a positive solution to \eqref{ch04.2.1}, one
can use It\^{o}'s formula to see the process is a solution of the
martingale problem associated with the generator $L$ defined by
\eqref{ch04.1.1}. By Theorem~\ref{tch04.2.1} we see $\{y(t)\}$ is a
CBI-process with transition semigroup $(Q_t^\gamma)_{t\ge 0}$. That
implies the weak uniqueness of the solution to \eqref{ch04.2.1}.
\qed\end{proof}

 \bgtheorem\label{tch04.2.2}
Suppose that $m(\d z) = \sigma z^{-1-\alpha}\d z$ for constants
$\sigma\ge 0$ and $1<\alpha<2$. Then the CBI-process with transition
semigroup $(Q_t^\gamma)_{t\ge 0}$ is the unique positive weak solution of
 \beqlb\label{ch04.2.3}
\d y(t) = \sqrt{2cy(t)}\d B(t) + \sqrt[\alpha]{\sigma y(t-)}\d z_0(t) -
by(t)\d t + \d z_1(t),
 \eeqlb
where $\{B(t)\}$ is a Brownian motion, $\{z_0(t)\}$ is a one-sided
$\alpha$-stable process with L\'{e}vy measure $z^{-1-\alpha}\d z$,
$\{z_1(t)\}$ is an increasing L\'{e}vy process defined by $(\beta,n)$,
and $\{B(t)\}$, $\{z_0(t)\}$ and $\{z_1(t)\}$ are independent of each
other. \edtheorem

 \begin{proof}
We assume $\sigma>0$, for otherwise the proof is easier. Let us consider
the CBI-process $\{y(t)\}$ given by \eqref{ch04.2.1} with $\{N_0(\d s,\d
z,\d u)\}$ being a Poisson random measure on $(0,\infty)^3$ with
intensity $\sigma z^{-1-\alpha}\d s\d z\d u$. We define the random
measure $\{N(\d s,\d z)\}$ on $(0,\infty)^2$ by
 \beqnn
N((0,t]\times B)
 \ar=\ar
\int_0^t\int_0^\infty\int_0^{y(s-)} 1_{\{y(s-)>0\}}1_B\Big(\frac{z}
{\sqrt[\alpha]{\sigma y(s-)}}\Big) N_0(\d s,\d z,\d u) \cr
 \ar~\ar
+ \int_0^t\int_0^\infty\int_0^{1/\sigma} 1_{\{y(s-)=0\}}1_B(z) N_0(\d
s,\d z,\d u).
 \eeqnn
It is easy to compute that $\{N(\d s,\d z)\}$ has predictable compensator
 \beqnn
\hat{N}((0,t]\times B)
 \ar=\ar
\int_0^t\int_0^\infty 1_{\{y(s-)>0\}}1_B\Big(\frac{z}
{\sqrt[\alpha]{\sigma y(s-)}}\Big)\frac{\sigma y(s-)\d s\d
z}{z^{1+\alpha}} \cr
 \ar~\ar
+ \int_0^t\int_0^\infty 1_{\{y(s-)=0\}}1_B(z)\frac{\d s\d z}
{z^{1+\alpha}} \cr
 \ar=\ar
\int_0^t\int_0^\infty 1_B(z)\frac{\d s\d z}{z^{1+\alpha}}.
 \eeqnn
Thus $\{N(\d s,\d z)\}$ is a Poisson random measure with intensity
$z^{-1-\alpha}\d s\d z$; see, e.g., Ikeda and Watanabe (1989, p.93). Now
define the L\'{e}vy processes
 \beqnn
z_0(t) = \int_0^t\int_0^\infty z \tilde{N}(\d s,\d z)
 ~~\mbox{and}~~
z_1(t) = \beta t + \int_0^t\int_0^\infty z N_1(\d s,\d z),
 \eeqnn
where $\tilde{N}(\d s,\d z) = N(\d s,\d z) - \hat{N}(\d s,\d z)$. It is
easy to see that
 \beqnn
\int_0^t\sqrt[\alpha]{\sigma y(s-)}\d z_0(s)
 \ar=\ar
\int_0^t\int_0^\infty \sqrt[\alpha]{\sigma y(s-)}\,z \tilde{N}(\d s,\d z)
\cr
 \ar=\ar
\int_0^t\int_0^\infty \int_0^{y(s-)} z \tilde{N}_0(\d s,\d z,\d u).
 \eeqnn
Then we get \eqref{ch04.2.3} from \eqref{ch04.2.1}. Conversely, if
$\{y(t)\}$ is a solution of \eqref{ch04.2.3}, one can use It\^{o}'s
formula to see that $\{y(t)\}$ solves the martingale problem associated
with the generator $L$ defined by \eqref{ch04.1.1} with $m(\d z) = \sigma
z^{-1-\alpha}\d z$. Then $\{y(t)\}$ is a CBI-process with transition
semigroup $(Q_t^\gamma)_{t\ge 0}$ and the solution of \eqref{ch04.2.3} is
unique in law. \qed\end{proof}

 \bgtheorem\label{tch04.2.3}
The pathwise uniqueness holds for positive solutions to \eqref{ch04.2.1}.
\edtheorem

 \begin{proof} For each integer $n\ge 0$ define $a_n = \exp\{-n(n+1)/2\}$. Then $a_n\to 0$ decreasingly as $n\to \infty$ and
 \beqnn
\int_{a_n}^{a_{n-1}}z^{-1} dz = n, \qquad n\ge 1.
 \eeqnn
Let $x\mapsto g_n(x)$ be a positive continuous function supported by
$(a_n,a_{n-1})$ so that
 \beqnn
\int_{a_n}^{a_{n-1}}g_n(x)dx=1
 \eeqnn
and $g_n(x)\le 2(nx)^{-1}$ for every $x>0$. For $n\ge 0$ let
 \beqnn
f_n(z)=\int_0^{|z|}dy\int_0^yg_n(x)dx, \qquad z\in \mbb{R}.
 \eeqnn
It is easy to see that $|f_n^\prime(z)|\le 1$ and
 \beqnn
0\le |z|f_n^{\prime\prime}(z) = |z|g_n(|z|)\le 2n^{-1}, \qquad z\in
\mbb{R}.
 \eeqnn
Moreover, we have $f_n(z)\rightarrow |z|$ increasingly as $n\to \infty$.
Suppose that $\{y(t): t\ge 0\}$ and $\{z(t): t\ge 0\}$ are both positive
solutions of \eqref{ch04.2.1}. Let $\alpha_t = z(t)-y(t)$ for $t\ge 0$.
From \eqref{ch04.2.1} we have
 \beqnn
\alpha_t \ar=\ar \alpha_0 - b\int_0^t \alpha_{s-} ds +
\sqrt{2c}\int_0^t\big(\sqrt{z(s)}-\sqrt{y(s)}\big)dB(s) \cr
 \ar\ar\qquad
+ \int_0^t\int_0^\infty \int_{y(s-)}^{z(s-)} z \tilde{N}_0(ds,dz,du).
 \eeqnn
By this and It\^o's formula,
 \beqlb\label{ch04.2.4}
f_n(\alpha_t) \ar=\ar f_n(\alpha_0) - b\int_0^t f_n^\prime(\alpha_s)
\alpha_s ds + c\int_0^t f_n^{\prime\prime}(\alpha_s) \big(\sqrt{z(s)} -
\sqrt{y(s)}\big)^2ds \cr
 \ar\ar
+ \int_0^t\alpha_s1_{\{\alpha_s>0\}}ds \int_0^\infty [f_n(\alpha_s+z) -
f_n(\alpha_s) - zf_n^\prime(\alpha_s)] m(dz) \cr
 \ar\ar
- \int_0^t\alpha_s1_{\{\alpha_s<0\}}ds \int_0^\infty [f_n(\alpha_s-z) -
f_n(\alpha_s) + zf_n^\prime(\alpha_s)] m(dz) \ccr
 \ar\ar
+\, \mbox{martingale.}
 \eeqlb
It is easy to see that $|f_n(a+x) - f_n(a)|\le |x|$ for any
$a,x\in\mbb{R}$. If $ax\ge 0$, we have
 \beqnn
|f_n(a+x) - f_n(a) - xf_n^\prime(a)|\le (2|ax|)\land (n^{-1}|x|^2).
 \eeqnn
Taking the expectation in both sides of \eqref{ch04.2.4} gives
 \beqnn
\mbf{P}[f_n(\alpha_t)]
 \ar\le\ar
\mbf{P}[f_n(\alpha_0)] + |b|\int_0^t \mbf{P}[|\alpha_s|] ds + c\int_0^t
\mbf{P}[f_n^{\prime\prime}(\alpha_s)|\alpha_s|] ds \cr
 \ar\ar
+ \int_0^tds\int_0^\infty \{(2z\mbf{P}[|\alpha_s|])\land (n^{-1}z^2)\}
m(dz) .
 \eeqnn
Then letting $n\to \infty$ we get
 \beqnn
\mbf{P}[|z(t)-y(t)|]
 \le
\mbf{P}[|z(0)-y(0)|] + |b|\int_0^t \mbf{P}[|z(s)-y(s)|] ds.
 \eeqnn
By this and Gronwall's inequality one can see the pathwise uniqueness
holds for \eqref{ch04.2.1}. \qed\end{proof}

Theorem~\ref{tch04.2.3} was first proved in Dawson and Li (2006). By
Theorems~\ref{tch04.2.1} and~\ref{tch04.2.3} there is a unique positive
strong solution to \eqref{ch04.2.1}; see, e.g., Situ (2005, p.76 and
p.104). The pathwise uniqueness of \eqref{ch04.2.3} was proved in Fu and
Li (2010).

\section{Lamperti's transformations by time changes}\label{sch04.3}

 \setcounter{equation}{0}

The results of Lamperti (1967b) assert that CB-processes are in
one-to-one correspondence with spectrally positive L\'{e}vy processes via
simple random time changes. Caballero et al.\ (2009) recently gave proofs
of those results using the approach of stochastic equations; see also
Helland (1978) and Silverstein (1968). Suppose that $\phi$ is a branching
mechanism given by \eqref{ch02.1.13}. Let $\{x(t): t\ge 0\}$ be a
CB-process with $x(0)=x\ge 0$ and with branching mechanism $\phi$ is
given by \eqref{ch02.1.13}. Let $\{Y_t: t\ge 0\}$ be a spectrally
positive L\'{e}vy process such that $Y_0=x$ and
 \beqlb\label{ch04.3.1}
\log \mbf{P}\exp\{i\lambda(Y_t-Y_r)\}
 =
(t-r)\phi(-i\lambda), \qquad \lambda\in \mbb{R}, t\ge r\ge 0.
 \eeqlb
Let $\tau = \inf\{s\ge 0: Y_s = 0\}$ be its first hitting time at zero
and let $Z_t = Y_{t\land \tau}$ for  $t\ge 0$. The proofs of the
following two theorems were essentially adopted from Caballero et al.\
(2009).

 \bgtheorem\label{tch04.3.1}
For any $t\ge 0$ let $z(t) = x(\kappa(t))$, where
 \beqlb\label{ch04.3.2}
\kappa(t)
 =
\inf\Big\{u\ge 0: \int_0^u x(s-)\d s = \int_0^u x(s)\d s\ge t\Big\}.
 \eeqlb
Then $\{z(t): t\ge 0\}$ is distributed identically on $D[0,\infty)$ with
$\{Z_t: t\ge 0\}$. \edtheorem

 \begin{proof}
By the result of Theorem~\ref{tch04.2.1}, we may assume $\{x(t)\}$ solves
the stochastic integral equation
 \beqlb\label{ch04.3.3}
x(t)
 \ar=\ar
x + \int_0^t \sqrt{2cx(s-)}\d B(s) - \int_0^t bx(s-)\d s \cr
 \ar~\ar
+ \int_0^t\int_0^\infty \int_0^{x(s-)} z \tilde{N}_0(\d s,\d z,\d u),
 \eeqlb
where $\{B(t)\}$ is a Brownian motion, $\{N_0(\d s,\d z,\d u)\}$ is a
Poisson random measures on $(0,\infty)^3$ with intensity $\d sm(\d z)\d
u$ and $\tilde{N}_0(\d s,\d z,\d u) = N_0(\d s,\d z,\d u) - \d sm(\d z)\d
u$. It follows that
 \beqlb\label{ch04.3.4}
z(t)
 \ar=\ar
x + \int_0^{\kappa(t)} \sqrt{2cx(s-)}\d B(s) - \int_0^{\kappa(t)}
bx(s-)\d s \cr
 \ar~\ar
+ \int_0^{\kappa(t)}\int_0^\infty \int_0^{x(s-)} z \tilde{N}_0(\d s,\d
z,\d u) \cr
 \ar=\ar
x + \sqrt{2c}W(t) - b\int_0^t z(s-) \d\kappa(s) \cr
 \ar~\ar
+ \int_0^t\int_0^\infty \int_0^{x(\kappa(s)-)} z \tilde{N}_0(\d
\kappa(s),\d z,\d u),
 \eeqlb
where
 \beqnn
W(t) = \int_0^{\kappa(t)} \sqrt{x(s-)}\d B(s)
 =
\int_0^t \sqrt{z(s-)} \d B(\kappa(s))
 \eeqnn
is a continuous martingale. From \eqref{ch04.3.2} we have
 \beqnn
1_{\{z(s-)>0\}}\d\kappa(s)
 =
1_{\{z(s-)>0\}}z(s-)^{-1}\d s.
 \eeqnn
Let $\tau_0 = \inf\{t\ge 0: z(t) = 0\}$. Since zero is a trap for
$\{z(t)\}$, we have
 \beqnn
\int_0^t z(s-) \d\kappa(s)
 =
\int_0^t 1_{\{z(s-)>0\}}\d s
 =
t\land\tau_0.
 \eeqnn
Then $\{W(t)\}$ has quadratic variation process $\<W\>(t) =
t\land\tau_0$, so it is a Brownian motion stopped at $\tau_0$. It is easy
to extend $\{W(t)\}$ to a Brownian motion with infinite time. Now define
the random measure $\{N(\d s,\d z)\}$ on $(0,\infty)^2$ by
 \beqnn
N((0,t]\times (a,b])
 \ar=\ar
\int_0^t\int_a^b\int_0^{z(s-)} 1_{\{z(s-)>0\}} N_0(\d\kappa(s),\d z,\d
u),
 \eeqnn
where $t\ge 0$ and $b\ge a>0$. It is easy to compute that
$\{N((0,t]\times (a,b]): t\ge 0\}$ has predictable compensator
 \beqnn
\hat{N}((0,t]\times (a,b])
 =
\int_0^t m(a,b]z(s-)\d\kappa(s)
 =
\int_0^t m(a,b]1_{\{s\le\tau_0\}} \d s.
 \eeqnn
Then we can extend $\{N(\d s,\d z)\}$ is a Poisson random measure on
$(0,\infty)^2$ with intensity $\d sm(\d z)$; see, e.g., Ikeda and
Watanabe (1989, p.93). From \eqref{ch04.3.4} we conclude that $\{z(t)\}$
is distributed on $D[0,\infty)$ identically with $\{Z_t: t\ge 0\}$.
\qed\end{proof}

 \bgtheorem\label{tch04.3.2}
For any $t\ge 0$ let $X_t = Z_{\theta(t)}$, where
 \beqlb\label{ch04.3.5}
\theta(t)
 =
\inf\Big\{u\ge 0: \int_0^u Z_{s-}^{-1}\d s = \int_0^u Z_s^{-1}\d s\ge
t\Big\}.
 \eeqlb
Then $\{X_t: t\ge 0\}$ is distributed identically on $D[0,\infty)$ with
$\{x(t): t\ge 0\}$. \edtheorem

 \begin{proof}
By the L\'{e}vy--It\^{o} decomposition, up to an extension of the
original probability space we may assume $\{Y_t\}$ is given by
 \beqnn
Y_t = x + \sqrt{2c}W(t) -  bt + \int_0^t\int_0^\infty \int_0^1 z
\tilde{M}_0(\d s,\d z,\d u),
 \eeqnn
where $\{W(t)\}$ is a Brownian motion, $\{M_0(\d s,\d z,\d u)\}$ is a
Poisson random measures on $(0,\infty)^3$ with intensity $\d sm(\d z)\d
u$ and $\tilde{M}_0(\d s,\d z,\d u) = M_0(\d s,\d z,\d u) - \d sm(\d z)\d
u$. It follows that
 \beqlb\label{ch04.3.6}
X_t = x + \sqrt{2c}W(\theta(t)) - b\theta(t) + \int_0^t \int_0^\infty
\int_0^1 z \tilde{M}_0(\d\theta(s),\d z,\d u).
 \eeqlb
From \eqref{ch04.3.5} we have
 \beqnn
\theta(t) = \int_0^t Z_{\theta(s)} \d s = \int_0^t X_s \d s.
 \eeqnn
Then the continuous martingale $\{W(\theta(t))\}$ has the representation
 \beqnn
W(\theta(t)) = \int_0^t \sqrt{X_s}\d B(s), \qquad t\ge 0
 \eeqnn
for another Brownian motion $\{B(t)\}$. Now we take an independent
Poisson random measure $\{M_1(\d s,\d z,\d u)\}$ on $(0,\infty)^3$ with
intensity $\d sm(\d z)\d u$ and define the random measure
 \beqnn
N_0(\d s,\d z,\d u) = 1_{\{u\le X_{s-}\}} M_0(\d\theta(s),\d
z,X_{s-}^{-1}\d u) + 1_{\{u>X_{s-}\}} M_1(\d s,\d z,\d u).
 \eeqnn
It is easy to see that $\{N_0(\d s,\d z,\d u)\}$ has deterministic
compensator $\d sm(\d z)\d u$, so it is a Poisson random measures. From
\eqref{ch04.3.6} we see that $\{X_t\}$ is a weak solution of
\eqref{ch04.3.3}. That gives the desired result. \qed\end{proof}

\chapter{State-dependent immigration structures}\label{ch05}

In this chapter we investigate the structures of state-dependent
immigration associated with CB-processes. For simplicity, we only
consider interactive immigration rates. The models are defined in terms
of some stochastic integral equations generalizing \eqref{ch04.2.1}. We
prove the existence and pathwise uniqueness of solutions to the
stochastic integral equations. Similar immigration structures were
studied in Li (2011) in the setting of superprocesses by considering
different types of stochastic equations. We shall deal with processes
with c\`{a}dl\`{a}g paths as in Li (2011).

\section{Time-dependent immigration}\label{sch05.1}

 \setcounter{equation}{0}

In this section, we introduce a generalization of the CBI-process. Let
$\phi$ is a branching mechanism given by \eqref{ch02.1.13} and let
$(Q_t)_{t\ge 0}$ be the transition semigroup defined by \eqref{ch02.1.15}
and \eqref{ch02.1.20}. We consider a set of functions $\{\psi_s: s\ge
0\}\subset \mcr{I}$ given by
 \beqlb\label{ch05.1.1}
\psi_s(z) = \beta_s z + \int_0^\infty\big(1-\e^{-zu}\big)n_s(\d u),
\qquad z\ge 0,
 \eeqlb
where $\beta_s\ge 0$ and $(1\land u)n_s(\d u)$ is a finite measure on
$(0,\infty)$. We assume $s\mapsto \psi_s(z)$ is locally bounded and
measurable on $[0,\infty)$ for each $z\ge 0$. By Theorems~\ref{tch01.2.3}
and~\ref{tch01.2.4}, for any $t\ge r\ge 0$ there is an infinitely
divisible probability measure $\gamma_{r,t}$ on $[0,\infty)$ defined by
 \beqlb\label{ch05.1.2}
\int_0^\infty \e^{-\lambda y} \gamma_{r,t}(\d y)
 =
\exp\Big\{-\int_r^t\psi_s(v_{t-s}(\lambda))\rho(s)\d s\Big\}, \qquad
\lambda\ge 0.
 \eeqlb
Then we can define the probability kernels $(Q_{r,t}^\gamma: t\ge r\ge
0)$ by
 \beqlb\label{ch05.1.3}
Q_{r,t}^\gamma(x,\cdot) := Q_{t-r}(x,\cdot)*\gamma_{r,t}(\cdot), \qquad
x\ge 0.
 \eeqlb
It is easily seen that
 \beqlb\label{ch05.1.4}
\int_0^\infty \e^{-\lambda y} Q_{r,t}^\gamma(x,\d y)
 =
\exp\Big\{-xv_{t-r}(\lambda)-\int_r^t\psi_s(v_{t-s}(\lambda))\d s\Big\}.
 \eeqlb
Moreover, the kernels $(Q_{r,t}^\gamma: t\ge r\ge 0)$ form a transition
semigroup on $\mbb{R}_+$. A Markov process with transition semigroup
$(Q_{r,t}^\gamma: t\ge r\ge 0)$ is called a \index{special inhomogeneous
CBI-process} \textit{special inhomogeneous CBI-process} with
\index{branching mechanism} \textit{branching mechanism} $\phi$ and
\index{time-dependent immigration mechanism} \textit{time-dependent
immigration mechanism} $\{\psi_s: s\ge 0\}$. One can see that the
time-space homogeneous transition semigroup associated with
$(Q_{r,t}^\gamma: t\ge r\ge 0)$ is a Feller semigroup. Then
$(Q_{r,t}^\gamma: t\ge r\ge 0)$ has a c\`{a}dl\`{a}g realization $X =
(\Omega, \mcr{F}, \mcr{F}_t, y(t), \mbf{Q}_{r,x}^\gamma)$. In particular,
if $un_s(\d u)$ is a locally bounded kernel from $[0,\infty)$ to
$(0,\infty)$, one can derive from \eqref{ch02.2.4} and \eqref{ch05.1.4}
that
 \beqlb\label{ch05.1.5}
\int_0^\infty y Q_{r,t}^\gamma(x,\d y)
 =
x\e^{-b(t-r)} + \int_r^t\e^{-b(t-s)}\psi'_s(0)\d s.
 \eeqlb
where
 \beqnn
\psi'_s(0) = \beta_s + \int_0^\infty z n_s(\d z).
 \eeqnn
The reader may refer to Li (2002) for the discussions of general
inhomogeneous immigration processes in the setting of measure-valued
processes.

\section{Predictable immigration rates}\label{sch05.2}

\setcounter{equation}{0}

Let $\phi$ be a branching mechanism given by \eqref{ch02.1.13} and $\psi$
an immigration mechanism given by \eqref{ch03.1.10}. In this section, we
give a construction of CBI-processes with random immigration rates given
by predictable processes. Suppose that $(\Omega, \mcr{F}, \mcr{F}_t,
\mbf{P})$ is a filtered probability space satisfying the usual
hypotheses. Let $\{B(t): t\ge 0\}$ be an $(\mcr{F}_t)$-Brownian motion
and let $\{p_0(t): t\ge 0\}$ and $\{p_1(t): t\ge 0\}$ be
$(\mcr{F}_t)$-Poisson point processes on $(0,\infty)^2$ with
characteristic measures $m(\d z)du$ and $n(\d z)du$, respectively. We
assume that the white noise and the Poisson processes are independent of
each other. Let $N_0(\d s,\d z,\d u)$ and $N_1(\d s,\d z,\d u)$ denote
the Poisson random measures on $(0,\infty)^3$ associated with
$\{p_0(t)\}$ and $\{p_1(t)\}$, respectively. Let $\tilde{N}_0(\d s,\d
z,\d u)$ denote the compensated measure of $N_0(\d s,\d z,\d u)$. Suppose
that $\rho = \{\rho(t): t\ge 0\}$ is a positive $(\mcr{F}_t)$-predictable
process such that $t\mapsto \mbf{P}[\rho(t)]$ is locally bounded. We are
interested in positive c\`adl\`ag solutions of the stochastic equation
 \beqlb\label{ch05.2.1}
Y_t \ar=\ar Y_0 + \sigma\int_0^t\sqrt{Y_{s-}} \d B(s) + \int_0^t
\int_0^\infty\int_0^{Y_{s-}} z \tilde{N}_0(\d s,\d z,\d u) \cr
 \ar\ar
+ \int_0^t (\beta\rho(s)-bY_{s-}) \d s + \int_0^t \int_0^\infty
\int_0^{\rho(s)} z N_1(\d s,\d z,\d u).
 \eeqlb

Clearly, the above equation is a generalization of \eqref{ch04.2.1}. For
any positive c\`adl\`ag solution $\{Y_t: t\ge 0\}$ of \eqref{ch05.2.1}
satisfying $\mbf{P}[Y_0]< \infty$, one can use a standard stopping time
argument to show that $t\mapsto \mbf{P}[Y_t]$ is locally bounded and
 \beqlb\label{ch05.2.2}
\mbf{P}[Y_t] = \mbf{P}[Y_0] + \psi'(0)\int_0^t \mbf{P}[\rho(s)] \d s -
b\int_0^t \mbf{P}[Y_s] \d s,
 \eeqlb
where $\psi'(0)$ is defined by \eqref{ch03.1.12}. By It\^o's formula, it
is easy to see that $\{Y_t: t\ge 0\}$ solves the following martingale
problem: For every $f\in C^2(\mbb{R}_+)$,
 \beqlb\label{ch05.2.3}
f(Y_t)
 \ar=\ar
f(Y_0) + \mbox{local mart.} - b\int_0^t f'(Y_s)Y_s\d s + \frac{1}{2}
\sigma^2\int_0^t f''(Y_s)Y_s\d s \cr
 \ar\ar
+ \int_0^tY_s\d s\int_0^\infty [f(Y_s + z) - f(Y_s) - zf'(Y_s)]m(\d z)
\cr
 \ar\ar
+ \int_0^t\rho(s)\Big\{\beta f'(Y_s) + \int_0^\infty [f(Y_s + z) -
f(Y_s)] n(\d z)\Big\}\d s.
 \eeqlb

 \bgproposition\label{tch05.2.1}
Suppose that $\{Y_t: t\ge 0\}$ is a positive c\`adl\`ag solution of
\eqref{ch05.2.1} and $\{Z_t: t\ge 0\}$ is a positive c\`adl\`ag solution
of the equation with $\rho = \{\rho(t): t\ge 0\}$ replaced by $\eta =
\{\eta(t): t\ge 0\}$. Then for any $t\ge 0$ we have
 \beqlb\label{ch05.2.4}
\mbf{P}[|Z_t-Y_t|]
 \le
\e^{|b|t}\Big\{\mbf{P}[|Z_0-Y_0|] + \psi'(0) \int_0^t
\mbf{P}[|\eta(s)-\rho(s)|] \d s\Big\},
 \eeqlb
where $\psi'(0)$ is defined by \eqref{ch03.1.12}. \edproposition

 \begin{proof}
For each integer $n\ge 0$ define $a_n = \exp\{-n(n+1)/2\}$. Then $a_n\to
0$ decreasingly as $n\to \infty$ and
 \beqnn
\int_{a_n}^{a_{n-1}}z^{-1} \d z = n, \qquad n\ge 1.
 \eeqnn
Let $x\mapsto g_n(x)$ be a positive continuous function supported by
$(a_n,a_{n-1})$ so that
 \beqnn
\int_{a_n}^{a_{n-1}}g_n(x)\d x=1
 \eeqnn
and $g_n(x)\le 2(nx)^{-1}$ for every $x>0$. Let
 \beqnn
f_n(z)=\int_0^{|z|}dy\int_0^yg_n(x)\d x, \qquad z\in \mbb{R}.
 \eeqnn
It is easy to see that $|f_n'(z)|\le 1$ and
 \beqnn
0\le |z|f_n''(z) = |z|g_n(|z|)\le 2n^{-1}, \qquad z\in \mbb{R}.
 \eeqnn
Moreover, we have $f_n(z)\rightarrow |z|$ increasingly as $n\to \infty$.
Let $\alpha_t = Z_t-Y_t$ for $t\ge 0$. From \eqref{ch05.2.1} we have
 \beqlb\label{ch05.2.5}
\alpha_t \ar=\ar \alpha_0 + \beta\int_0^t [\eta(s)-\rho(s)] \d s +
\sigma\int_0^t (\sqrt{Z_{s-}} - \sqrt{Y_{s-}}) \d B(s) \cr
 \ar\ar
-\, b\int_0^t \alpha_{s-} \d s + \int_0^t \int_0^\infty
\int_{Y_{s-}}^{Z_{s-}} z \tilde{N}_0(\d s,\d z,\d u) \cr
 \ar\ar
+ \int_0^t \int_0^\infty \int_{\rho(s)}^{\eta(s)} z N_1(\d s,\d z,\d u).
 \eeqlb
By this and It\^o's formula,
 \beqlb\label{ch05.2.6}
f_n(\alpha_t) \ar=\ar f_n(\alpha_0) + \beta\int_0^t f_n'(\alpha_s)
[\eta(s)-\rho(s)] \d s - b\int_0^t f_n'(\alpha_s)\alpha_s \d s \cr
 \ar\ar
+\, {1\over 2}\sigma^2\int_0^t f_n''(\alpha_s)
(\sqrt{Z_{s-}}-\sqrt{Y_{s-}})^2\d s \cr
 \ar\ar
+ \int_0^t\alpha_s1_{\{\alpha_s>0\}}\d s \int_0^\infty [f_n(\alpha_s+z) -
f_n(\alpha_s) - zf_n'(\alpha_s)] m(\d z) \cr
 \ar\ar
- \int_0^t\alpha_s1_{\{\alpha_s<0\}}\d s \int_0^\infty [f_n(\alpha_s-z) -
f_n(\alpha_s) + zf_n'(\alpha_s)] m(\d z) \cr
 \ar\ar
+ \int_0^t[\eta(s)-\rho(s)]1_{\{\eta(s)>\rho(s)\}}\d s \int_0^\infty
[f_n(\alpha_s+z) - f_n(\alpha_s)] n(\d z) \cr
 \ar\ar
- \int_0^t[\rho(s)-\eta(s)]1_{\{\rho(s)>\eta(s)\}}\d s \int_0^\infty
[f_n(\alpha_s-z) - f_n(\alpha_s)] n(\d z) \ccr
 \ar\ar
+\, \mbox{martingale.}
 \eeqlb
It is easy to see that $|f_n(a+x) - f_n(a)|\le |x|$ for any
$a,x\in\mbb{R}$. If $ax\ge 0$, we have
 \beqnn
|f_n(a+x) - f_n(a) - xf_n'(a)|\le (2|ax|)\land (n^{-1}|x|^2).
 \eeqnn
Taking the expectation in both sides of \eqref{ch05.2.6} gives
 \beqnn
\mbf{P}[f_n(\alpha_t)]
 \ar\le\ar
\mbf{P}[f_n(\alpha_0)] + \beta\int_0^t \mbf{P}[|\eta(s)-\rho(s)|] \d s +
|b|\int_0^t \mbf{P}[|\alpha_s|] \d s \cr
 \ar\ar
+\, \int_0^t \mbf{P}[|\eta(s)-\rho(s)|]\d s \int_0^\infty z n(\d z) +
n^{-1}\sigma^2t \cr
 \ar\ar
+\, \int_0^t\d s\int_0^\infty \{(2z\mbf{P}[|\alpha_s|]) \land
(n^{-1}z^2)\} m(\d z).
 \eeqnn
By letting $n\to \infty$ we get
 \beqlb\label{ch05.2.7}
\mbf{P}[|Z_t-Y_t|]
 \ar\le\ar
\mbf{P}[|Z_0-Y_0|] + |b|\int_0^t \mbf{P}[|Z_s-Y_s|] \d s \cr
 \ar\ar
+\, \psi'(0) \int_0^t \mbf{P}[|\eta(s)-\rho(s)|] \d s.
 \eeqlb
Then we get the desired estimate follows by Gronwall's inequality.
\qed\end{proof}

 \bgproposition\label{tch05.2.2}
Suppose that $\{Y_t: t\ge 0\}$ is a positive c\`adl\`ag solution of
\eqref{ch05.2.1} and $\{Z_t: t\ge 0\}$ is a positive c\`adl\`ag solution
of the equation with $(b,\rho)$ replaced by $(c,\eta)$. Then for any
$t\ge 0$ we have
 \beqnn
\mbf{P}\Big[\sup_{0\le s\le t}|Z_s-Y_s|\Big]
 \ar\le\ar
\mbf{P}[|Z_0-Y_0|] + \psi'(0) \int_0^t \mbf{P}[|\eta(s)-\rho(s)|] \d s
\cr
 \ar\ar
+ \Big(|b|+2\int_1^\infty z m(\d z)\Big)\int_0^t \mbf{P}[|Z_s-Y_s|] \d s
\cr
 \ar\ar
+\, 2\sigma \Big(\int_0^t\mbf{P}[|Z_s-Y_s|]\d s\Big)^{1\over 2} \cr
 \ar\ar
+\, 2\Big(\int_0^t \mbf{P}[|Z_s-Y_s|]\d s \int_0^1 z^2 m(\d
z)\Big)^{1\over 2},
 \eeqnn
where $\psi'(0)$ is defined by \eqref{ch03.1.12}. \edproposition

 \begin{proof}
This follows by applying Doob's martingale inequality to
\eqref{ch05.2.5}. \qed\end{proof}

 \bgtheorem\label{tch05.2.3}
For any $Y_0\ge 0$ there is a pathwise unique positive c\`adl\`ag
solution $\{Y_t: t\ge 0\}$ of \eqref{ch05.2.1}. \edtheorem

 \begin{proof}
The pathwise uniqueness of the solution follows by
Proposition~\ref{tch05.2.1} and Gronwall's inequality. Without loss of
generality, we may assume $Y_0\ge 0$ is deterministic in proving the
existence of the solution. We give the proof in two steps.

\noindent\textit{Step~1.} Let $0=r_0<r_1<r_2<\cdots$ be an increasing
sequence. For each $i\ge 1$ let $\eta_i$ be a positive integrable random
variable measurable with respect to $\mcr{F}_{r_{i-1}}$. Let $\rho =
\{\rho(t): t\ge 0\}$ be the positive $(\mcr{F}_t)$-predictable step
process given by
 \beqnn
\rho(t) = \sum_{i=1}^\infty \eta_i 1_{(r_{i-1},r_i]}(t), \qquad t\ge 0.
 \eeqnn
By Theorem~\ref{tch04.2.1}, on each interval $(r_{i-1},r_i]$ there is a
pathwise unique solution $\{Y_t: r_{i-1}<t\le r_i\}$ to
 \beqnn
Y_t \ar=\ar Y_{r_{i-1}} + \sigma\int_{r_{i-1}}^t\sqrt{Y_{s-}} \d B(s) +
\int_{r_{i-1}}^t \int_0^\infty\int_0^{Y_{s-}} z \tilde{N}_0(\d s,\d z,\d
u) \cr
 \ar\ar\qquad
+ \int_{r_{i-1}}^t (\beta\eta_i-bY_{s-}) \d s + \int_{r_{i-1}}^t
\int_0^\infty \int_0^{\eta_i} z N_1(\d s,\d z,\d u).
 \eeqnn
Then $\{Y_t: t\ge 0\}$ is a solution to \eqref{ch05.2.1}.

\noindent\textit{Step~2.} Suppose that $\rho = \{\rho(t): t\ge 0\}$ is
general positive $(\mcr{F}_t)$-predictable process such that $t\mapsto
\mbf{P}[\rho(t)]$ is locally bounded. Take a sequence of positive
predictable step processes $\rho_k = \{\rho_k(t): t\ge 0\}$ so that
 \beqlb\label{ch05.2.8}
\mbf{P}\Big[\int_0^t |\rho_k(s) - \rho(s)| \d s\Big] \to 0
 \eeqlb
for every $t\ge 0$ as $k\to \infty$. Let $\{Y_k(t): t\ge 0\}$ be the
solution to \eqref{ch05.2.1} with $\rho=\rho_k$. By
Proposition~\ref{tch05.2.1}, Gronwall's inequality and \eqref{ch05.2.8}
one sees
 \beqnn
\sup_{0\le s\le t}\mbf{P}[|Y_k(s)-Y_i(s)|]\to 0
 \eeqnn
for every $t\ge 0$ as $i,k\to \infty$. Then Proposition~\ref{tch05.2.2}
implies
 \beqnn
\mbf{P}\Big[\sup_{0\le s\le t}|Y_k(s)-Y_i(s)|\Big]\to 0
 \eeqnn
for every $t\ge 0$ as $i,k\to \infty$. Thus there is a subsequence
$\{k_i\}\subset \{k\}$ and a c\`adl\`ag process $\{Y_t: t\ge 0\}$ so that
 \beqnn
\sup_{0\le s\le t}|Y_{k_i}(s)-Y_s|\to 0
 \eeqnn
almost surely for every $t\ge 0$ as $i\to \infty$. It is routine to show
that $\{Y_t: t\ge 0\}$ is a solution to \eqref{ch05.2.1}. \qed\end{proof}

 \bgtheorem\label{tch05.2.4}
If $\rho = \{\rho(t): t\ge 0\}$ is a deterministic locally bounded
positive Borel function, the solution $\{Y_t: t\ge 0\}$ of
\eqref{ch05.2.1} is a special inhomogeneous CBI-process with branching
mechanism $\phi$ and time-dependent immigration mechanisms
$\{\rho(t)\psi: t\ge 0\}$. \edtheorem

 \begin{proof}
By Theorem~\ref{tch04.2.1}, when $\rho(t) = \rho$ is a deterministic
constant function, the process $\{Y_t: t\ge 0\}$ is a CBI-process with
branching mechanism $\phi$ and immigration mechanisms $\rho\psi$. If
$\rho = \{\rho(t): t\ge 0\}$ is a general deterministic locally bounded
positive Borel function, we can take each step function $\rho_k =
\{\rho_k(t): t\ge 0\}$ in the last proof to be deterministic. Then the
solution $\{Y_k(t): t\ge 0\}$ of \eqref{ch05.2.1} with $\rho=\rho_k$ is a
special inhomogeneous CBI-process with branching mechanism $\phi$ and
time-dependent immigration mechanisms $\{\rho_k(t)\psi: t\ge 0\}$. In
other words, for any $\lambda\ge 0$, $t\ge r\ge 0$ and $G\in \mcr{F}_r$
we have
 \beqnn
\mbf{P}[1_G\e^{-\lambda Y_k(t)}]
 =
\mbf{P}\Big[1_G\exp\Big\{-Y_k(r)v_{t-r}(\lambda) -
\int_r^t\rho_k(s)\psi(v_{t-s}(\lambda))\d s\Big\}\Big].
 \eeqnn
Letting $k\to \infty$ along the sequence $\{k_i\}$ mentioned in the last
proof gives
 \beqnn
\mbf{P}[1_G\e^{-\lambda Y_t}]
 =
\mbf{P}\Big[1_G\exp\Big\{-Y_rv_{t-r}(\lambda) -
\int_r^t\rho(s)\psi(v_{t-s}(\lambda))\d s\Big\}\Big].
 \eeqnn
Then $\{Y_t: t\ge 0\}$ is a CBI-process with immigration rate $\rho =
\{\rho(t): t\ge 0\}$. \qed\end{proof}

In view of the result of Theorem~\ref{tch05.2.4}, the solution $\{Y_t:
t\ge 0\}$ to \eqref{ch05.2.1} can be called an inhomogeneous CBI-process
with branching mechanism $\phi$, immigration mechanism $\psi$ and
\index{predictable immigration rate} \textit{predictable immigration
rate} $\rho = \{\rho(t): t\ge 0\}$. The results in this section are
slight modifications of those in Li (2011+), where some path-valued
branching processes were introduced.

\section{Interactive immigration rates}\label{sch05.3}

 \setcounter{equation}{0}

In this section, we give a construction of CBI-processes with interactive
immigration rates. We shall use the set up of the second section. Suppose
that $z\mapsto q(z)$ is a positive Lipschitz function on $[0,\infty)$. We
consider the stochastic equation
 \beqlb\label{ch05.3.1}
Y_t \ar=\ar Y_0 + \sigma\int_0^t\sqrt{Y_{s-}} \d B(s) + \int_0^t
\int_0^\infty\int_0^{Y_{s-}} z \tilde{N}_0(\d s,\d z,\d u) \cr
 \ar\ar
+ \int_0^t [\beta q(Y_{s-})-bY_{s-}] \d s + \int_0^t \int_0^\infty
\int_0^{q(Y_{s-})} z N_1(\d s,\d z,\d u).
 \eeqlb
This reduces to \eqref{ch04.2.1} when $q$ is a constant function. We may
interpret the solution $\{Y_t: t\ge 0\}$ of \eqref{ch05.3.1} as a
CBI-process with \index{interactive immigration rate} \textit{interactive
immigration rate} given by the process $s\mapsto q(Y_{s-})$.

 \bgtheorem\label{tch05.3.1}
There is a pathwise unique solution $\{Y_t: t\ge 0\}$ of
\eqref{ch05.3.1}. \edtheorem

 \begin{proof}
Suppose that $\{Y_t: t\ge 0\}$ and $\{Z_t: t\ge 0\}$ are two solutions to
this equation. Let $K\ge 0$ be a Lipschitz constant for the function
$z\mapsto q(z)$. From \eqref{ch05.2.4} we get
 \beqnn
\mbf{P}[|Z_t-Y_t|]
 \le
\psi'(0)\e^{|b|t}\int_0^t \mbf{P}[|q(Z_s)-q(Y_s)|] \d s
 \le
K\psi'(0)\e^{|b|t} \int_0^t \mbf{P}[|Z_s-Y_s|] \d s.
 \eeqnn
Then the pathwise uniqueness for \eqref{ch05.3.1} follows by Gronwall's
inequality. We next prove the existence of the solution using an
approximating argument. Let $Y_0(t)\equiv 0$. By Theorem~\ref{tch05.2.3}
we can define inductively the sequence of processes $\{Y_k(t): t\ge 0\}$,
$k=1,2,\dots$ as pathwise unique solutions of the stochastic equations
 \beqlb\label{ch05.3.2}
Y_k(t) \ar=\ar Y_0 - b\int_0^t Y_k(s-) \d s + \sigma\int_0^t
\sqrt{Y_k(s-)} \d B(s) \cr
 \ar\ar
+\, \beta\int_0^t q(Y_{k-1}(s-)) \d s + \int_0^t\int_0^\infty
\int_0^{Y_k(s-)} z \tilde{N}_0(\d s,\d z,\d u)  \cr
 \ar\ar
+ \int_0^t \int_0^\infty \int_0^{q(Y_{k-1}(s-))} z N_1(\d s,\d z,\d u).
 \eeqlb
Let $Z_k(t) = Y_k(t) - Y_{k-1}(t)$. By \eqref{ch05.2.4} we have
 \beqnn
\mbf{P}[|Z_k(t)|]
 \ar\le\ar
\psi'(0)\e^{|b|t}\int_0^t \mbf{P}[|q(Y_{k-1}(s))-q(Y_{k-2}(s))|] \d s \cr
 \ar\le\ar
K\psi'(0)\e^{|b|t} \int_0^t \mbf{P}[|Z_{k-1}(s)|] \d s.
 \eeqnn
By \eqref{ch05.3.2} one sees that $\{Z_1(t): t\ge 0\}$ is a CBI-process
with branching mechanism $\phi$ and immigration mechanism $q(0)\psi$. In
view of \eqref{ch03.1.11}, we have
 \beqnn
\mbf{P}[|Z_1(t)|]
 =
\e^{-bt}\mbf{P}[Y_0] + \psi'(0)\int_0^t \e^{-bs} \d s.
 \eeqnn
By a standard argument, one shows
 \beqnn
\sum_{k=1}^\infty\sup_{0\le s\le t}\mbf{P}[|Y_k(s) - Y_{k-1}(s)|]<
\infty,
 \eeqnn
so the Lipschitz property of $z\mapsto q(z)$ implies
 \beqnn
\sum_{k=1}^\infty\sup_{0\le s\le t}\mbf{P}[|q(Y_k(s)) - q(Y_{k-1}(s))|]<
\infty.
 \eeqnn
It follows that
 \beqnn
\lim_{k,l\to \infty}\int_0^t \mbf{P}[|q(Y_k(s)) - q(Y_l(s))|]\d s = 0.
 \eeqnn
Then there exists a predictable process $\rho = \{\rho(s): s\ge 0\}$ so
that
 \beqlb\label{ch05.3.3}
\lim_{k\to \infty}\int_0^t \mbf{P}[|q(Y_k(s)) - \rho(s)|]\d s = 0.
 \eeqlb
Let $\{Y_t: t\ge 0\}$ be the positive c\`{a}dl\`{a}g process defined by
\eqref{ch05.2.1}. By Proposition~\ref{tch05.2.2}, there is a subsequence
$\{k_n\}\subset \{k\}$ so that a.s.
 \beqnn
\lim_{n\to\infty} \sup_{0\le s\le t}|Y_{k_n}(s) - Y_s| = 0, \qquad t\ge
0.
 \eeqnn
By the continuity of $z\mapsto q(z)$ we get a.s.
 \beqnn
\lim_{n\to\infty} q(Y_{k_n}(s-)) = q(Y(s-)), \qquad t\ge 0.
 \eeqnn
This and \eqref{ch05.3.3} imply that
 \beqnn
\int_0^t \mbf{P}[|q(Y(s-)) - \rho(s)|]\d s = 0.
 \eeqnn
Then letting $k\to \infty$ along $\{k_n+1\}$ in \eqref{ch05.3.2} we see
$\{Y_t: t\ge 0\}$ is a solution of \eqref{ch05.3.1}.  \qed\end{proof}

By It\^o's formula, it is easy to see that the solution $\{Y_t: t\ge 0\}$
of \eqref{ch05.3.1} solves the following martingale problem: For every
$f\in C^2(\mbb{R}_+)$,
 \beqlb\label{ch05.3.4}
f(Y_t)
 \ar=\ar
f(Y_0) + \mbox{local mart.} - b\int_0^t f'(Y_s)Y_s\d s + \frac{1}{2}
\sigma^2\int_0^t f''(Y_s)Y_s\d s \cr
 \ar\ar
+ \int_0^tY_s\d s\int_0^\infty [f(Y_s + z) - f(Y_s) - zf'(Y_s)]m(\d z)
\cr
 \ar\ar
+ \int_0^tq(Y_s)\Big\{\beta f'(Y_s) + \int_0^\infty [f(Y_s + z) - f(Y_s)]
n(\d z)\Big\}\d s.
 \eeqlb
By Theorem~\ref{tch05.3.1}, the solution is a strong Markov process with
generator given by
 \beqlb\label{ch05.3.5}
Af(x)
 \ar=\ar
\frac{1}{2}\sigma^2xf''(x) + x\int_0^\infty [f(x+z)-f(x)-zf'(x)]m(\d z)
\cr
 \ar\ar
-\, bxf'(x) + q(x)\Big\{\beta f'(x) + \int_0^\infty [f(x+z)-f(x)] n(\d
z)\Big\}.
 \eeqlb

We can also consider two Lipschitz functions $z\mapsto q_1(z)$ and
$z\mapsto q_2(z)$ on $[0,\infty)$. By slightly modifying the arguments,
one can show there is a pathwise unique solution to
 \beqlb\label{ch05.3.6}
Y_t \ar=\ar Y_0 + \sigma\int_0^t\sqrt{Y_{s-}} \d B(s) + \int_0^t
\int_0^\infty\int_0^{Y_{s-}} z \tilde{N}_0(\d s,\d z,\d u) \cr
 \ar\ar
+ \int_0^t [\beta q_1(Y_{s-})-bY_{s-}] \d s + \int_0^t \int_0^\infty
\int_0^{q_2(Y_{s-})} z N_1(\d s,\d z,\d u).
 \eeqlb
The solution of this equation can be understood as a CBI-process with
\textit{interactive immigration rates} given by the processes $s\mapsto
q_1(Y_{s-})$ and $s\mapsto q_2(Y_{s-})$. This type of immigration
structures were studied in Li (2011) in the setting of superprocesses by
considering a different type of stochastic equations.

\backmatter

 \addcontentsline{toc}{chapter}{Index}
 \printindex

 \end{document}